\documentclass[12pt]{article}
\usepackage{latexsym}
\usepackage[T1]{fontenc}
\usepackage{amsmath,amsfonts,amssymb}
\usepackage{theorem}
\theoremstyle{plain}
\theoremheaderfont{\scshape}

\newtheorem{definition}{Definition}[section]

\newtheorem{proposition}{Proposition}[section]
\newtheorem{theorem}{Theorem}[section]

\theorembodyfont{\rmfamily}
\newtheorem{remark}{Remark}[section]

\def\cA{{\cal A}}          
\def\cB{{\cal B}}          
\def\cC{{\cal C}}
\def\cD{{\cal D}}                                     
\def\cF{{\cal F}}
\def\cH{{\cal H}}          
\def\cI{{\cal I}}
\def\cK{{\cal K}}
\def\cN{{\cal N}}
\def\cR{{\cal R}}
\def\cU{{\cal U}}
\def\cX{{\cal X}}
\def\cZ{{\cal Z}}

\def\fg{{\mathfrak g}}
\def\fA{{\mathfrak A}}
\def\fB{{\mathfrak B}}
\def\fH{{\mathfrak H}}

\def\CC{{\mathbb C}}
\def\HH{{\mathbb H}}
\def\II{{\mathbb I}}
\def\NN{{\mathbb N}}
\def\ZZ{{\mathbb Z}}

\def\Ad{\mathop{\rm Ad}\nolimits}
\def\ad{\mathop{\rm ad}\nolimits}
\def\id{\mathop{\rm id}\nolimits}
\def\op{\mathop{\rm op}\nolimits}
\def\qdet{\text{q-}\!\det}
\def\qqquad{\quad\quad\quad}

\newcommand{\elpa}[1]{{{\cA}_{q,p}(\widehat{sl}_{#1})}}
\newcommand{\elpb}[1]{{{\cB}_{q,p,\lambda}(\widehat{sl}_{#1})}}
\newcommand{\dy}[1]{{{\cD}Y(#1)}}
\newcommand{\ddy}[2]{{{\cD}Y_{#2}(#1)}}
\newcommand{\sfrac}[2]{{\textstyle{\frac{#1}{#2}}}}
\newcommand{\half}{{\scriptstyle{\frac{1}{2}}}}
\newcommand{\car}[2]{\genfrac{[}{]}{0pt}{}{#1}{#2}}
\newcommand{\hypergeom}[4]{{_{2}F_{1}\left(\begin{array}{c}{#1}\quad{#2}\\
                           {#3}\end{array};{#4}\right)}}
\newcommand{\finproof}{{\hfill \rule{5pt}{5pt}}}
\newcommand{\finrmq}{{\hfill $\square$}}

\begin{document}
\baselineskip=15pt
\pagestyle{empty}

\begin{center}
{\LARGE \textbf{\textsf{Quantum elliptic algebras and double Yangians}}} 
\\
\vspace*{24pt}
{\large \textbf{\textsf{Luc Frappat}}}
\\
\vspace*{24pt} 
\emph{Laboratoire d'Annecy-le-Vieux de Physique Th{\'e}orique (LAPTH)} 
\\
\emph{UMR 5108, CNRS-Universit{\'e} de Savoie}
\\
\emph{B.P. 110, F-74941 Annecy-le-Vieux Cedex, France}
\\
\vspace*{8pt} 
\emph{Membre de l'Institut Universitaire de France} 
\\
\vspace*{24pt}
Lectures given at the First French-Moroccan School on Non-Commutative 
Geometry
\\
Marseille (France), 17-21 December 2001
\end{center}

\vfill

\begin{abstract}
Quantum universal enveloping algebras, quantum elliptic algebras and
double (deformed) Yangians provide fundamental algebraic structures
relevant for many integrable systems. They are described in the FRT
formalism by R-matrices which are solutions of elliptic, trigonometric
or rational type of the Yang--Baxter equation with spectral parameter
or its generalization known as the Gervais--Neveu--Felder equation.
While quantum groups and double Yangians appear as quasi-triangular
Hopf algebras, this is no longer the case for elliptic algebras and
the various deformations of Yangian type algebras. These structures
are dealt with the framework of quasi-Hopf algebras. These algebras
can be obtained from Hopf algebras through particular Drinfel'd twists
satisfying the so-called shifted cocycle condition. We review these
different structures and the pattern of connections between them.
\end{abstract}

\vfill


\vspace*{24pt}

\rightline{LAPTH-Conf-898/02}
\rightline{math.QA/0201245}
\rightline{January 2002}

\newpage

\def\contentsname{}
\vspace*{-60pt}
\baselineskip=14pt
\tableofcontents


\clearpage
\baselineskip=18pt
\pagestyle{myheadings}
\markright{Quantum elliptic algebras and double Yangians} 
\pagenumbering{arabic}
\setcounter{page}{1}
\parindent=0pt

\section{Introduction}
\setcounter{equation}{0}

\parskip=6pt

The quantum generalization of the inverse scattering method by
E.K.~Sklyanin, L.D.~Faddeev and L.A.~Takhtajan gives a unified
framework for the resolution of integrable systems and puts in
evidence underlying algebraic structures \cite{SF78,SFT80}. These
structures appear to be a very powerful tool both for resolution and
construction of two-dimensional models in statistical mechanics and in
quantum field theory. They also lead to the formulation of the
classical inverse scattering problem in an algebraic way with the
formalism of the classical $r$-matrix. The quantization of classical
$r$-matrices yields quantum R-matrices \cite{Skl80} satisfying the
so-called star-triangle relation, also known as the Yang--Baxter
equation \cite{Bax72,Yang67}.

A particularly interesting case is that of R-matrices $R \in
\text{End}(V \otimes V)$ depending on a complex parameter, called the
spectral parameter, where $V$ is a representation of an
infinite-dimensional algebra. In \cite{KuSkl82} a first attempt to
give a systemization of the (at this time) known solutions of the
Yang--Baxter equation with spectral parameter was undertaken. However,
the problem of classifying such solutions is not solved yet, and since
this pioneering work no new type of solution has been found. In fact,
the known solutions of the Yang--Baxter equation with spectral
parameter fall into three classes: the elliptic ones, the
trigonometric ones and the rational ones.

The trigonometric solutions have been used in solvable lattice model,
see for example \cite{IIJMNT,Jim86,JKMO88}. They appear also to be the
R-matrices used to define the quantum affine algebras in the FRT
formalism. The notion of quantum affine algebra was introduced in
\cite{KuRe81,Skl85} for the $\widehat{sl}_{2}$ case and in
\cite{Dri85,Jim85,Jim86} for any simple Kac--Moody affine algebra.

Elliptic R-matrices have yielded many algebraic structures relevant to
certain integrable systems in quantum mechanics and statistical
mechanics, e.g. the $XYZ$ model \cite{JKKMW}, RSOS models
\cite{ABF,Kon98} and Sine--Gordon theory \cite{Kon97,KLCP}. In the
elliptic case, the R-matrix associated to the eight vertex model was
first derived by R.J.~Baxter \cite{Bax82} and generalized to the
$\ZZ_{N}$-symmetric vertex model by A.A.~Belavin \cite{Bel81}. The
vertex-type elliptic R-matrix for $sl_{2}$ was first used by
E.K.~Sklyanin \cite{Skl82,Skl83} to construct a two-parameter
deformation of the enveloping algebra $\cU(sl_{2})$. The affine
version, including the central extension $c$, was introduced in
\cite{FIJKMY}, provided the Yang--Baxter relations take the form
$R_{12}L_{1}L_{2} = L_{2}L_{1}R_{12}^*$ where $R_{12}^*$ denotes the
R-matrix $R_{12}$ with a modified elliptic nome $p^* = pq^{-2c}$. At
this level, the coalgebra structure remains undetermined: indeed the
fact that the R-matrices $R_{12}$ and $R_{12}^*$ are different prevent
the usual coproduct formula $\Delta L = L \dot\otimes L$ from
applying. Moreover a free field realization of $\elpa{N}$ is still
missing. The connection between the quantum elliptic algebras
$\elpa{N}$ and the $q$-deformed Virasoro and $W_N$ algebras
\cite{SKAO96,AKOS96,FF96,FR96} was established in \cite{AFRS3,AFRS5}.

Other elliptic solutions, known as the face-type R-matrices, depending
on extra parameters $\lambda$ belonging to the dual of the Cartan
algebra in the underlying algebra, were first used by G.~Felder
\cite{Fel94} to define the algebra $\elpb{N}$ in the R-matrix
approach. This class of elliptic solutions to the Yang--Baxter
equation is associated with the face-type statistical models
\cite{ABF,DJMO,JMO88}. B.~Enriquez and G.~Felder \cite{EF97} and
H.~Konno \cite{Kon98} introduced a current representation, although
differences arise in the treatment of the central extension. A
slightly different structure, also based upon face-type R-matrices but
incorporating extra Heisenberg algebra generators, was introduced as
$\cU_{q,p}(\widehat{sl}(2))$ \cite{JKOS2,Kon98}. This structure is
relevant to the resolution of the quantum Calogero--Moser and
Ruijsenaars--Schneider models \cite{ABB96,BBB96,JS97}. It is
worthwhile to notice that the RLL type relations take in this case a
particular form, like the Yang--Baxter equation, where the Lax
matrices also depend on the extra parameters $\lambda$, and where some
shifts along the Cartan generators occur. As in the vertex case, the
coalgebra structure was undetermined.

This problem was solved in \cite{JKOS} both for the vertex-type and
for the face-type elliptic algebras. O.~Babelon \cite{Bab91} was the
first to suggest that the notion of quasi-Hopf algebra was the right
framework to deal with such structures. Originating with the
proposition of \cite{BBB96} on face-type algebras, the construction of
the twist operators was undertaken in both cases by C.~Fr{\o}nsdal
\cite{Fro97b,Fro97a} and finally achieved at the level of formal
universal twists in \cite{ABRR,JKOS}. In \cite{ABRR}, the universal
twist is obtained by solving a linear equation introduced in
\cite{BR99}, this equation playing a fundamental role for complex
continuation of $6j$ symbols. Moreover in the case of finite
(super)algebras, the convergence of the infinite products defining the
twists was also proved in \cite{ABRR}. This led to a formal
construction of universal R-matrices for the elliptic algebras
$\elpa{N}$ and $\elpb{N}$ in \cite{JKOS}, of which the Belavin--Baxter
and the Andrews--Baxter--Forrester $4 \!\times\! 4$ matrices are
respectively (spin $1/2$) evaluation representations.

\bigskip

Particularly interesting limits of elliptic or trigonometric
R-matrices are the so-called degeneration or scaling limits. The
resulting R-matrices are rational solutions of the Yang--Baxter
equation. The scaling limit is defined by taking $z=q^{u}$ (spectral
parameter), $p=q^{2r}$ (elliptic nome) and $w = q^{2s}$ (dynamical
parameter) with $q \to 1$, keeping $u$, $r$ and $s$ fixed. The
degeneracy limits of elliptic R-matrices was considered either in the
vertex case \cite{JKM,Kon97} (R-matrix formulation) and \cite{KLP}
(current algebra formulation) or in the face case \cite{Cla}. They
give rise to algebraic structures which have been variously
characterized as scaled elliptic algebras denoted
$\cA_{\hbar,\eta}(\widehat{sl}_{2})$, where $\eta \equiv 1/r$ and $q
\simeq e^{\epsilon\hbar}$ with $\epsilon \to 0$ \cite{JKM,KLP}, or
double Yangian type algebras \cite{Cla,BL93,KT96}. The algebra
$\cA_{\hbar,\eta}(\widehat{sl}_{2})$ is relevant to the study of the
$XXZ$ model in its gapless regime \cite{JKM}. It admits a further
limit $r \to \infty$ ($\eta \to 0$) where its R-matrix becomes
identical to the R-matrix defining the double Yangian $\dy{\fg}$
(centrally extended), defined in \cite{Kho95,KT96}.

One must however be careful in these identifications in terms of
R-matrix structure since the generating functionals (Lax matrices) of
these algebras admit different interpretations in terms of modes
(generators of the enveloping algebra). Although represented by
formally identical RLL type relations, these two classes of objects
differ fundamentally in their structure: in the context of scaled
elliptic algebras, the expansion is done in terms of continuous-index
Fourier modes of the spectral parameter \cite{KLP,Kon97}, while in the
context of double Yangian type algebras the expansion is done in terms
of powers of the spectral parameter \cite{Kho95,KT96}. Hence these
structures must be considered separately.

It appears clearly that the universal algebraic structures associated
with any limit of evaluated R-matrices may not be taken for granted,
but must be explicitly constructed. This will be achieved by
identification of these particular limits as evaluation
representations of universal R-matrices for the deformations by
particular Drinfel'd twists \cite{Dri90}, known as ``shifted-cocycle''
twists \cite{Bab91,BBB96}, of Hopf algebra structures. This
construction systematically endows these deformations with a
Gervais--Neveu--Felder type quasi-triangular quasi-Hopf algebra
structure. It is characterized by a particular form of the universal
Yang--Baxter equation, called Gervais--Neveu--Felder or dynamical
Yang--Baxter equation.

\bigskip

The plan of the paper is as follows. We first introduce in section 2
the basic tools necessary for the description of the various algebraic
structures we will deal with, especially the notions of Hopf algebra,
quasi-Hopf algebra and Drinfel'd twist. We also describe the structure
of quantum affine algebras, in particular the construction of their
universal R-matrix and the connections between their various
presentations (Cartan--Weyl basis, Drinfel'd realization, FRT
formalism). The section 3 is devoted to the quantum affine elliptic
algebras of vertex type $\elpa{N}$ and of face type $\elpb{N}$ and
some of their ``degenerations''. Emphasis is given on the structure of
these algebras, in particular their quasi-Hopf structure. The
Drinfel'd twists that connect them to the quantum group
$\cU_{q}(\widehat{sl}_{N})$ are explicitly constructed following the
lines of \cite{JKOS}. In section 4, we deal with (double) Yangians and
related structures. The degeneration limit of $\elpa{N}$, interpreted
as a deformed double Yangian, is shown to be obtained from the double
Yangian $\dy{sl_{N}}$ by a Drinfel'd twist of Gervais--Neveu--Felder
type, promoting the deformed double Yangian as a quasi-triangular
quasi-Hopf algebra. A similar degeneration limit of the face type
elliptic algebra $\elpb{N}$ and removing the elliptic dependence lead
to a dynamical double Yangian algebra, which also exhibits the
structure of quasi-triangular quasi-Hopf algebra. Finally the nature
of the degeneration limit of $\elpb{N}$ without removing the elliptic
dependence is discussed (``dynamical deformed double Yangian'').

\parskip=0pt


\section{Hopf algebras, quasi-Hopf algebras and Drinfel'd twists}
\setcounter{equation}{0}

\subsection{Hopf algebras}

\subsubsection{Definitions}

\begin{definition}
A unital associative algebra $\fA$ over $\CC$ is called a Hopf algebra 
if it is endowed with a coalgebra structure: the coproduct $\Delta : 
\fA \to \fA \otimes \fA$ and counit $\epsilon : \fA \to \CC$ are 
algebra homomorphisms, the product $m : \fA \otimes \fA \to \fA$ and 
unit $\iota : \CC \to \fA$ are coalgebra homomorphisms~%
\footnote{$\varphi : \fA \to \fB$ is a coalgebra homomorphism if 
$\Delta_{\fB} \circ \varphi = (\varphi \otimes \varphi) \circ 
\Delta_{\fA}$ and $\epsilon_{\fB} \circ \varphi = \epsilon_{\fA}$.}, 
and $\fA$ is equipped with an antihomomorphism $S : \fA \to \fA$ 
(antipode), with
\begin{align*}
& (\id \otimes \Delta) (\Delta(x)) = (\Delta \otimes \id) 
(\Delta(x)) \qquad (\forall x \in \fA) \qquad 
\text{(coassociativity)} \\
& (\id \otimes \epsilon) \circ \Delta = (\epsilon \otimes \id) 
\circ \Delta = \id \\
& m \circ (S \otimes \id) \circ \Delta = m \circ (\id \otimes S) 
\circ \Delta = \iota \circ \epsilon
\end{align*}
\end{definition}

If $\sigma$ denotes the flip map, $\sigma(x \otimes y) = y \otimes x$ for 
$x,y \in \fA$, $\Delta^{\op} = \sigma \circ \Delta$ is the opposite 
coproduct. Whenever $\Delta = \Delta^{\op}$, the Hopf algebra $\fA$ is said 
cocommutative.

\begin{definition}
A Hopf algebra $\fA$ is said to be quasi-triangular if it exists an 
invertible element $\cR \in \fA \otimes \fA$, called the universal 
R-matrix, such that
\begin{align*}
& \Delta^{\op}(x) = \cR \Delta(x) \cR^{-1} \qquad (\forall x \in 
\fA) \qquad \text{(almost cocommutativity)} \\
& (\Delta \otimes \id) (\cR) = \cR_{13} \, \cR_{23} \\
& (\id \otimes \Delta) (\cR) = \cR_{13} \, \cR_{12}
\end{align*}
\end{definition}

It follows that $\cR$ satisfies the \emph{Yang--Baxter equation} (in $\fA 
\otimes \fA \otimes \fA$):
\begin{equation}
\cR_{12} \, \cR_{13} \, \cR_{23} = \cR_{23} \, \cR_{13} \, \cR_{12}
\end{equation}
and that
\begin{equation}
\label{eq:Reps}
\begin{split}
& (\epsilon \otimes \id) (\cR) = (\id \otimes \epsilon) (\cR) = 1 
\\
& (S \otimes \id) (\cR) = (\id \otimes S^{-1}) (\cR) = \cR^{-1}
\end{split}
\end{equation}

\subsubsection{Quantum affine algebras $\cU_{q}(\widehat\fg)$}
\label{sect:212}

Let $\fg$ be a finite-dimensional complex simple Lie algebra and
$\CC[z,z^{-1}]$ be the ring of Laurent polynomials in the
indeterminate $z$. The Kac--Moody affine algebra $\widehat\fg$ is
defined by $\widehat\fg = \fg \otimes \CC[z,z^{-1}] \oplus \CC c
\oplus \CC d$, $c$ is the central extension and $d$ is the derivation.
Let $(A_{ij})$ and $(A^{sym}_{ij})$ be the Cartan matrix and the
symmetrized Cartan matrix of $\widehat\fg$ respectively, where
$A^{sym}_{ij} = d_{i} A_{ij}$ and $d_{i}$ are relatively coprime
integers such that $d_{i} A_{ij} = d_{j} A_{ji}$. If $\widehat\Pi^0 =
\{ \alpha_{0}, \ldots, \alpha_{r} \}$ is a simple root system of
$\widehat\fg$, one has $(A^{sym}_{ij}) = (\alpha_{i},\alpha_{j})$
where $(\,\cdot\,,\,\cdot\,)$ denotes a scalar product on the root
space.

\medskip

The quantum universal enveloping Kac--Moody affine algebra -- or quantum 
affine algebra -- $\cU_{q}(\widehat\fg)$ is the unital algebra over 
$\CC(q)$, the ring of rational functions in the indeterminate $q$, with 
generators $k_{i}^{\pm}$, $e_{i}$ and $f_{i}$ ($0 \le i \le r$) such that
\begin{align}
& [k_{i}^{\pm} , k_{j}^{\pm}] = 0 && k_{i}^{+} k_{i}^{-} = k_{i}^{-} 
k_{i}^{+} = 1 && [e_{i} , f_{j}] = \delta_{ij} \; \frac{k_{i}^{+} - 
k_{i}^{-}}{q-q^{-1Ñ}} \\
& k_{i}^{\pm} e_{j} = q^{\pm A^{sym}_{ij}} e_{j} k_{i}^{\pm} && 
k_{i}^{\pm} f_{j} = q^{\mp A^{sym}_{ij}} f_{j} k_{i}^{\pm} && \\
\intertext{and for $i \ne j$}
& (\ad_{q} e_{i})^{1-A_{ij}} (e_{j}) = 0 && (\ad_{q} f_{i})^{1-A_{ij}} 
(f_{j}) = 0 &&
\end{align}
where by definition $(\ad_{q} e_{i}) (e_{j}) = e_{i} e_{j} - 
q^{A^{sym}_{ij}} e_{j} e_{i}$. \\
The Hopf algebra structure of $\cU_{q}(\widehat\fg)$ is given by
\begin{align}
& \Delta(k_{i}^{\pm}) = k_{i}^{\pm} \otimes k_{i}^{\pm} && 
\Delta(e_{i}) = e_{i} \otimes 1 + k_{i}^{+} \otimes e_{i} && 
\Delta(f_{i}) = f_{i} \otimes k_{i}^{-} + 1 \otimes f_{i} \nonumber \\
& S(k_{i}^{\pm}) = k_{i}^{\mp} && S(e_{i}) = -k_{i}^{-} e_{i} && 
S(f_{i}) = f_{i} k_{i}^{+} \\
& \epsilon(k_{i}^{\pm}) = \epsilon(e_{i}) = \epsilon(f_{i}) = 0 && 
\epsilon(1) = 1 && \nonumber
\end{align}

\medskip

Let $\widehat\Pi^+$ be the set of positive roots of $\widehat\fg$. One 
defines a so-called normal ordering $\prec$ among the roots of 
$\widehat\Pi^+$ in the following way: if $\alpha, \beta, \alpha+\beta \in 
\widehat\Pi^+$ and $[\alpha,\beta]$ is a minimal segment containing 
$\alpha+\beta$, then one writes $\alpha \prec \alpha+\beta \prec \beta$. 
Once the generators $e_{\alpha_{i}} \equiv e_{i}$, $e_{-\alpha_{i}} \equiv 
f_{i}$ corresponding to the simple roots are given, the generators 
corresponding to all positive roots are constructed by induction as follows 
\cite{KT94} (see also \cite{KT91,KT93} for the finite case):
\begin{equation}
\begin{split}
& e_{\alpha+\beta} = [e_{\alpha} , e_{\beta}]_{q} = e_{\alpha} 
e_{\beta} - q^{(\alpha,\beta)} e_{\beta} e_{\alpha} \\
& f_{\alpha+\beta} = [f_{\beta} , f_{\alpha}]_{q^{-1}} = f_{\beta} 
f_{\alpha} - q^{-(\alpha,\beta)} f_{\alpha} f_{\beta}
\end{split}
\end{equation}
if $\alpha, \beta, \alpha+\beta \in \widehat\Pi^+$ are such that $\alpha 
\prec \alpha+\beta \prec \beta$.

For any root $\gamma = \sum_{i=0}^{r} n_{i} \alpha_{i}$, $\gamma \in 
\Pi^{+}$, one sets $k_{\gamma}^{\pm} = \prod_{i=0}^{r} 
{k_{i}^{\pm}}^{n_{i}}$. If $\gamma$ is a real root, one gets
\begin{equation}
[e_{\gamma},f_{\gamma}] = \eta_{\gamma} \; 
\frac{k_{\gamma}^{+}-k_{\gamma}^{-}}{q-q^{-1}}
\end{equation}
but for the imaginary roots $n\delta$, one has (the superscript $(i)$ 
denotes the multiplicity of the imaginary root $n\delta$)
\begin{equation}
[e_{n\delta}^{(i)},f_{m\delta}^{(j)}] \ne a_{ij}(n) \; 
\frac{{k_{\delta}^{+}}^n-{k_{\delta}^{-}}^n}{q-q^{-1}} \; 
\delta_{m+n,0}
\end{equation}
It is necessary instead to introduce new generators 
$\check{e}_{n\delta}^{(i)}$, $\check{f}_{n\delta}^{(j)}$ such that
\begin{equation}
[\check{e}_{n\delta}^{(i)},\check{f}_{m\delta}^{(j)}] = a_{ij}(n) \; 
\frac{{k_{\delta}^{+}}^n-{k_{\delta}^{-}}^n}{q-q^{-1}} \; 
\delta_{m+n,0}
\end{equation}
where $\displaystyle a_{ij}(n) = 
\frac{q^{nA^{sym}_{ij}}-q^{-nA^{sym}_{ij}}}{n(q-q^{-1})}$. \\
The generators $e_{n\delta}^{(i)}$ are expressed in terms of the 
$\check{e}_{n\delta}^{(i)}$ generators by means of Schur relations (see 
\cite{KT94}).

\medskip

The universal R-matrix of the quantum affine algebra $\cU_{q}(\widehat\fg)$ 
has been constructed in ref. \cite{KT92,KT94}. It has the following 
structure:
\begin{equation}
\cR[\cU_{q}(\widehat{\fg})] = \Big( \prod_{\gamma \in 
\widehat\Pi^{+}}^{\rightarrow} \widehat{\cR}_{\gamma} \Big) \cK
\end{equation}
where the arrow means that the product has to be done with respect to the 
normal ordering $\prec$ defined on $\widehat\Pi^+$ and the factors 
$\widehat{\cR}_{\gamma}$ and $\cK$ are given by
\begin{equation}
\label{eq:rhat1}
\widehat{\cR}_{\gamma} = \exp_{q^{-(\gamma,\gamma)}} \Big( (q-q^{-1}) 
\eta_{\gamma}^{-1} \, e_{\gamma} \otimes f_{\gamma} \Big)
\end{equation}
for the real roots $\gamma \in \widehat\Pi^{+}$,
\begin{equation}
\widehat{\cR}_{n\delta} = \exp \Big( (q-q^{-1}) \sum_{ij} c_{ij}(n) \, 
\check{e}_{n\delta}^{(i)} \otimes \check{f}_{n\delta}^{j)} \Big)
\end{equation}
for the imaginary roots $n\delta \in \widehat\Pi^{+}$ where the matrix 
$\big(c_{ij}(n)\big)$ is the inverse of the matrix $\big(a_{ij}(n)\big)$, 
and
\begin{equation}
{\cal K} = q^{\sum_{ij} d_{ij} h_{i} \otimes h_{j}}
\end{equation}
where $h_{i}$ is such that $k_{i}^{\pm} \equiv q^{\pm h_{i}}$ and the 
matrix $(d_{ij})$ is the inverse of a non-degenerated extension of the 
symmetrized Cartan matrix $(A^{sym}_{ij})$. \\
In (\ref{eq:rhat1}) the $q$-exponential is defined by
\begin{equation}
\exp_{q}(x) \equiv \sum_{n \in \NN} \frac{x^n}{(n)_{q}!} \qquad 
\text{where} \qquad (n)_{q}! \equiv (1)_{q} (2)_{q} \ldots (n)_{q} \;\; 
\text{and} \;\; (k)_{q} \equiv \frac{1-q^k}{1-q}
\end{equation}

\subsubsection{The Drinfel'd realization of $\cU_{q}(\widehat\fg)$}
\label{sect:drinfeld}

There exists another realization of the quantum affine algebra 
$\cU_{q}(\widehat\fg)$ found by Drinfel'd \cite{Dri86,Dri88}. 
$\cU_{q}(\widehat\fg)$ is isomorphic to the algebra with generators 
$\cK_{i,n}$, $\cX^{\pm}_{i,n}$ and $C$ (with $1 \le i \le r$ and $n \in 
\ZZ$) such that
\begin{equation}
\label{eq:reldri}
\begin{split}
& \big[ \cK_{i,m},\cK_{j,n} \big] = \delta_{n+m,0} \; 
\frac{q^{mA^{sym}_{ij}}-q^{-mA^{sym}_{ij}}}{m(q-q^{-1})} \; 
\frac{C^m-C^{-m}}{q-q^{-1}} \\
& \cK_{i,0} \; \cX^{\pm}_{j,n} = q^{\pm A^{sym}_{ij}} \; 
\cX^{\pm}_{j,n} \; \cK_{i,0} \\
& \big[ \cK_{i,m},\cX^{\pm}_{j,n} \big] = \pm 
\frac{q^{mA^{sym}_{ij}}-q^{-mA^{sym}_{ij}}}{m(q-q^{-1})} \; 
\cX^{\pm}_{j,m+n} \; C^{\mp |m|/2} \\
& \cX^{\pm}_{i,m+1}\cX^{\pm}_{j,n} - q^{\pm A^{sym}_{ij}} 
\cX^{\pm}_{j,n} \cX^{\pm}_{i,m+1} = q^{\pm A^{sym}_{ij}} 
\cX^{\pm}_{i,m}\cX^{\pm}_{j,n+1} - \cX^{\pm}_{j,n+1} 
\cX^{\pm}_{i,m} \\
& \big[ \cX^+_{i,m},\cX^-_{j,n} \big] = \delta_{ij} \; 
\frac{\Psi^+_{i,m+n} C^{(m-n)/2} - \Psi^-_{i,m+n} 
C^{(n-m)/2}}{q-q^{-1}} \\
& \big[ C,\cK_{i,m} \big] = \big[ C,\cX^{\pm}_{i,m} \big] = 0
\end{split}
\end{equation}
and for $i \ne j$ with $n_{ij} = 1-A_{ij}$
\begin{equation}
\sum_{\sigma \in \mathfrak{S}_{n_{ij}}} \sum_{k=0}^{n_{ij}} (-1)^k \; 
\frac{[n_{ij}]_{q}!}{[k]_{q}! [n_{ij}-k]_{q}!} \; 
\cX^{\pm}_{i,m_{\sigma(1)}} \, \ldots \, \cX^{\pm}_{i,m_{\sigma(k)}} \, 
\cX^{\pm}_{j,n} \, \cX^{\pm}_{i,m_{\sigma(k+1)}} \, \ldots \, 
\cX^{\pm}_{i,m_{\sigma(n_{ij})}} = 0
\end{equation}
The generators $\Psi^{\pm}_{i,n}$ are determined by
\begin{equation}
\sum_{n \ge 0} \Psi^{\pm}_{i,\pm n} \, z^{\mp n} = \cK^{\pm 1}_{i,0} 
\exp\Big( \pm(q-q^{-1}) \sum_{n \ge 1} \cK_{i,\pm n} \, z^{\mp n} \Big)
\end{equation}
and the $q$-factorial is defined by $[n]_{q}! \equiv [1]_{q} \, [2]_{q} \, 
\ldots \, [n]_{q}$ and $\displaystyle [k]_{q} \equiv 
\frac{q^{k}-q^{-k}}{q-q^{-1}}$.

\subsubsection{FRT formalism of $\cU_{q}(\widehat{sl}_{2})$}

Consider the two-dimensional evaluation representation of 
$\widehat{sl}_{2}$ with evaluation parameter $z$, in the homogeneous 
gradation given by
\begin{equation}
\label{eq:evalrepr}
\begin{array}{ll}
\pi_z(e_1) = e_{12} \qqquad \pi_z(f_1) = e_{21} & \qqquad 
\pi_z(e_0) = z e_{21} \qqquad \pi_z(f_0) = z^{-1} e_{12} \\
\Big. \pi_z(h_1) = e_{11} - e_{22} & \qqquad \pi_z(h_0) = e_{22} - 
e_{11}
\end{array}
\end{equation}
Then the R-matrix of $\cU_{q}(\widehat{sl}_{2})$ in the fundamental 
representation $R(z_{1}/z_{2}) = (\pi_{z_{1}} \otimes \pi_{z_{2}}) \cR$ 
reads
\begin{equation}
\label{eq:ruq}
R[\cU_{q}(\widehat{sl}_{2})](z) = \rho(z) \left(
\begin{array}{cccc}
1 & 0 & 0 & 0 \\
0 & \displaystyle \frac{q(1-z)}{1-q^2z} & \displaystyle 
\frac{1-q^2}{1-q^2z} & 0 \\[12pt] 0 & \displaystyle 
\frac{z(1-q^2)}{1-q^2z} & \displaystyle \frac{q(1-z)}{1-q^2z} & 0 
\\
0 & 0 & 0 & 1 \\
\end{array}
\right)
\end{equation}
the normalization factor being
\begin{equation}
\label{eq:rhouq}
\rho(z) = q^{-1/2} \; \frac{(q^2z;q^4)_{\infty}^2}{(z;q^4)_{\infty} \; 
(q^4z;q^4)_{\infty}}
\end{equation}
The infinite multiple products are defined by $(z;a)_\infty = \prod_{n \ge 
0} (1-z a^n)$.

\medskip

In the FRT formalism \cite{FRT}, the algebra $\cU_{q}(\widehat{sl}_{2})$ can be 
alternatively defined as an algebra with generators $L^{\pm}_{ij}(z) = 
\sum_{k \ge 0} L_{ij}^{\pm}(\mp k) \, z^{\pm k}$, encapsulated into the two 
$2 \!\times\! 2$ independent matrices
\begin{equation}
L^{\pm}(z) = \left(
\begin{array}{cc}
L_{11}^{\pm}(z) & L_{12}^{\pm}(z) \\
L_{21}^{\pm}(z) & L_{22}^{\pm}(z) \\
\end{array}
\right)
\end{equation}
subject to the relations
\begin{equation}
\begin{split}
R_{12}(z_1/z_2) \, L^{\pm}_1(z_1) \, L^{\pm}_2(z_2) &= 
L^{\pm}_2(z_2) \, L^{\pm}_1(z_1) \, R_{12}(z_1/z_2) \\
R_{12}(q^{c} z_1/z_2) \, L^+_1(z_1) \, L^-_2(z_2) &= L^-_2(z_2) \, 
L^+_1(z_1) \, R_{12}(q^{-c} z_1/z_2)
\end{split}
\end{equation}
and
\begin{equation}
\qdet L^{\pm}(z) \equiv L^{\pm}_{11}(q^{-1}z) L^{\pm}_{22}(z) - 
L^{\pm}_{21}(q^{-1}z) L^{\pm}_{12}(z) = 1
\end{equation}
where $L^{\pm}_1(z) = L^{\pm}(z) \otimes 1$, $L^{\pm}_2(z) = 1 \otimes 
L^{\pm}(z)$.

The Hopf algebra structure is given by
\begin{equation}
\Delta L^{\pm}(z) = L^{\pm}(zq^{\pm(1 \otimes c/2)}) \,\dot\otimes\, 
L^{\pm}(zq^{\mp c/2 \otimes 1})
\end{equation}
i.e. explicitly
\begin{equation}
\Delta L_{ij}^{\pm}(z) = \sum_{k} L_{ik}^{\pm}(zq^{\pm(1 \otimes c/2)}) 
\otimes L_{kj}^{\pm}(zq^{\mp c/2 \otimes 1})
\end{equation}
for the coproduct, $S(L^{\pm}(z)) = {(L^{\pm}(z))}^{-1}$ for the antipode 
and $\epsilon(L^{\pm}) = 1$ for the counit.

Using a Gauss decomposition for the Lax matrices $L^{\pm}(z)$,
\begin{align}
L^{\pm}(z) &= \left(
\begin{array}{cc}
1 & 0 \\
e^{\pm}(z) & 1 \\
\end{array}
\right) \left(
\begin{array}{cc}
k_1^{\pm}(z) & 0 \\
0 & k_2^{\pm}(z) \\
\end{array}
\right) \left(
\begin{array}{cc}
1 & f^{\pm}(z^\mp) \\
0 & 1 \\
\end{array}
\right) \nonumber \\[6pt] &= \left(
\begin{array}{cc}
k_1^{\pm}(z) & k_1^{\pm}(z) f^{\pm}(z^\mp) \\
e^{\pm}(z) k_1^{\pm}(z) & k_2^{\pm}(z) + e^{\pm}(z) k_1^{\pm}(z) 
f^{\pm}(z^\mp) \\
\end{array}
\right)
\end{align}
the modes of the generating functions
\begin{align}
& \cX^+(z) = (q-q^{-1})^{-1} \big( e^+(qz_{-}) - e^-(qz_{+}) \big) = 
\sum_{n \in \ZZ} \cX^+_{n} \, z^{-n} \nonumber \\
& \cX^-(z) = (q-q^{-1})^{-1} \big( f^+(qz_{+}) - f^-(qz_{-}) \big) = 
\sum_{n \in \ZZ} \cX^-_{n} \, z^{-n} \\
& \Psi^{\pm}(z) = k_2^\mp(zq) \, k_1^\mp(zq)^{-1} = \sum_{n \in \NN} 
\Psi^{\pm}_{n} \, z^{\mp n} \nonumber
\end{align}
with $z_{\pm} = zq^{\pm c/2}$ satisfy the Drinfel'd relations 
(\ref{eq:reldri}) of $\cU_{q}(\widehat{sl}_{2})$ \cite{DF93}.

\subsection{Quasi-Hopf algebras}

\begin{definition}
A unital associative algebra $\fA$ over $\CC$ is called a quasi-Hopf 
algebra if it is endowed with a coalgebra structure: the coproduct 
$\Delta : \fA \to \fA \otimes \fA$ and counit $\epsilon : \fA \to \CC$ 
are algebra homomorphisms, the product $m : \fA \otimes \fA \to \fA$ 
and unit $\iota : \CC \to \fA$ are coalgebra homomorphisms, and $\fA$ 
is equipped with an antihomomorphism $S : \fA \to \fA$ (antipode) and 
elements $\alpha$, $\beta \in \fA$, and an invertible element $\Phi \in 
\fA \otimes \fA \otimes \fA$ (coassociator), with
\begin{align*}
& (\id \otimes \Delta) (\Delta(x)) = \Phi (\Delta \otimes \id) 
(\Delta(x)) \Phi^{-1} \qquad (\forall x \in \fA) \\
& (\id \otimes \epsilon) \circ \Delta = (\epsilon \otimes \id) 
\circ \Delta = \id \\
& (\id \otimes \id \otimes \Delta) (\Phi) \cdot (\Delta \otimes \id 
\otimes \id) (\Phi) = (1 \otimes \Phi) \cdot (\id \otimes \Delta 
\otimes \id) (\Phi) \cdot (\Phi \otimes 1) \\
& (\id \otimes \epsilon \otimes \id) (\Phi) = 1
\end{align*}
and for the antipode
\begin{align*}
& \sum_{i} S(x_{i}^{(1)}) \alpha x_{i}^{(2)} = \epsilon(x) \alpha 
&& \sum_{i} x_{i}^{(1)} \beta S(x_{i}^{(2)}) = \epsilon(x) \beta \\
& \sum_{i} \varphi_{i}^{(1)} \beta S(\varphi_{i}^{(2)}) \alpha 
\varphi_{i}^{(3)} = 1 && \sum_{i} S(\psi_{i}^{(1)}) \alpha 
\psi_{i}^{(2)} \beta S(\psi_{i}^{(3)}) = 1
\end{align*}
where $x \in \fA$ with $\Delta(x) = \sum_{i} x_{i}^{(1)} \otimes 
x_{i}^{(2)}$ and
\begin{equation*}
\Phi = \sum_{i} \varphi_{i}^{(1)} \otimes \varphi_{i}^{(2)} \otimes 
\varphi_{i}^{(3)} \;, \qquad \Phi^{-1} = \sum_{i} \psi_{i}^{(1)} 
\otimes \psi_{i}^{(2)} \otimes \psi_{i}^{(3)}
\end{equation*}
\end{definition}
The element $\Phi$ measures the lack of coassociativity of the coproduct.

\begin{definition}
A quasi-Hopf algebra $\fA$ is said to be quasi-triangular if it exists 
an invertible element $\cR \in \fA \otimes \fA$, called the universal 
R-matrix, such that
\begin{align*}
& \Delta^{\op}(x) = \cR \Delta(x) \cR^{-1} \qquad (\forall x \in 
\fA) \\
& (\Delta \otimes \id) (\cR) = \Phi^{(312)} \, \cR_{13} \, 
{\Phi^{(132)}}^{-1} \, \cR_{23} \, \Phi^{(123)} \\
& (\id \otimes \Delta) (\cR) = {\Phi^{(231)}}^{-1} \, \cR_{13} \, 
\Phi^{(213)} \, \cR_{12} \, {\Phi^{(123)}}^{-1}
\end{align*}
\end{definition}
It follows that $\cR$ satisfies the generalized Yang--Baxter equation (in 
$\fA \otimes \fA \otimes \fA$):
\begin{equation}
\cR_{12} \, \Phi^{(312)} \, \cR_{13} \, {\Phi^{(132)}}^{-1} \, \cR_{23} 
\, \Phi^{(123)} = \Phi^{(321)} \, \cR_{23} \, {\Phi^{(231)}}^{-1} \, 
\cR_{13} \, \Phi^{(213)} \, \cR_{12}
\end{equation}
The notation $\Phi^{(312)}$ means that if $\Phi^{(123)} = \sum_{i} 
\varphi_{i}^{(1)} \otimes \varphi_{i}^{(2)} \otimes \varphi_{i}^{(3)}$, 
then $\Phi^{(312)} = \sum_{i} \varphi_{i}^{(3)} \otimes \varphi_{i}^{(1)} 
\otimes \varphi_{i}^{(2)}$, and so on.

\smallskip

Obviously, a quasi-Hopf algebra with $\Phi = 1$, $\alpha = \beta = 1$ is a 
Hopf algebra.

\subsection{Drinfel'd twist}

The notion of Drinfel'd twist allows one to associate to a given 
quasi-triangular quasi-Hopf algebra another quasi-triangular quasi-Hopf 
algebra in the following way. Consider an invertible element $\cF \in \fA 
\otimes \fA$ such that $(\id \otimes \epsilon) \cF = (\epsilon \otimes \id) 
\cF = 1$ (when $\fA$ is a quantum universal enveloping algebra, this means 
that the ``leading'' term in $\cF$ is $1 \otimes 1$). One sets
\begin{align}
& \widetilde\Delta(x) = \cF_{12} \, \Delta(x) \, \cF_{12}^{-1} \qquad 
(\forall x \in \fA) \\
& \widetilde\cR = \cF_{21} \, \cR_{12} \, \cF_{12}^{-1} \\
& \widetilde\Phi = \big( \cF_{23} (\id \otimes \Delta)(\cF) \big) \, 
\Phi \, \big( \cF_{12} (\Delta \otimes \id)(\cF) \big)^{-1}
\label{eq:deftwistc}
\\
& \widetilde\alpha = \sum_{i} S(w_{i}^{(1)}) \, \alpha \, w_{i}^{(2)} 
\qquad \text{and} \qquad \widetilde\beta = \sum_{i} v_{i}^{(1)} \, 
\beta \, S(v_{i}^{(2)})
\end{align}
where
\begin{equation}
\cF_{12} = \sum_{i} v_{i}^{(1)} \otimes v_{i}^{(2)} \qquad \text{and} 
\qquad \cF_{12}^{-1} = \sum_{i} w_{i}^{(1)} \otimes w_{i}^{(2)}
\end{equation}
\begin{proposition}[Drinfel'd]
If $(\fA$, $\Phi$, $\Delta$, $\epsilon$, $S$, $\alpha$, $\beta$, $\cR)$ 
is a quasi-triangular quasi-Hopf algebra (QTQHA), then $(\fA$, 
$\widetilde\Phi$, $\widetilde\Delta$, $\epsilon$, $S$, 
$\widetilde\alpha$, $\widetilde\beta$, $\widetilde\cR)$ is also a QTQHA 
(denoted $\widetilde\fA$ for short).
\end{proposition}
$\cF$ is called a Drinfel'd twist.

\medskip

In the following, we will mainly be concerned with twists of Hopf algebras. 
{From} now on, we consider the case where $\fA$ is a Hopf algebra ($\Phi = 
1$, $\alpha = \beta = 1$) and $\cF$ depends on parameters $\lambda \in \fH$ 
where $\fH$ is an Abelian subalgebra of $\fA$.

\begin{definition}
A Drinfel'd twist $\cF$ satisfying the so-called \emph{shifted cocycle 
condition} ($h \in \fH$)
\begin{equation}
\cF_{12}(\lambda) \, (\Delta \otimes \id) (\cF(\lambda)) = 
\cF_{23}(\lambda + h^{(1)}) \, (\id \otimes \Delta) (\cF(\lambda))
\end{equation}
is called a \emph{Gervais--Neveu--Felder} (GNF) twist.
\end{definition}

In the case of a GNF twist, the coassociator $\widetilde\Phi$ is given by 
$\widetilde\Phi^{(123)} = \cF_{23}(\lambda) \, \cF_{23}(\lambda + 
h^{(1)})^{-1}$ (from eq. (\ref{eq:deftwistc})) and the universal R-matrix 
$\widetilde{\cR}$ satisfies the so-called \emph{Gervais--Neveu--Felder} or 
\emph{dynamical Yang--Baxter equation} \cite{Fel94,GN84}:
\begin{equation}
\label{eq:DYBE}
\widetilde{\cR}_{12}(\lambda + h^{(3)}) \, 
\widetilde{\cR}_{13}(\lambda) \, \widetilde{\cR}_{23}(\lambda + 
h^{(1)}) = \widetilde{\cR}_{23}(\lambda) \, 
\widetilde{\cR}_{13}(\lambda + h^{(2)}) \, 
\widetilde{\cR}_{12}(\lambda)
\end{equation}

Denoting by $\pi_{ev}(x)$ an evaluation representation of $\fA$ with 
evaluation parameter $x$, the Lax matrix $\widetilde{L}(x) = (\pi_{ev}(x) 
\otimes \id) \, \widetilde{\cR}$ realizes an FRT-type formalism of 
$\widetilde\fA$ with an evaluated R-matrix defined by 
$\widetilde{R}(x_1-x_2) = (\pi_{ev}(x_{1}) \otimes \pi_{ev}(x_{2})) \, 
\widetilde{\cR}$. The RLL relations take the form
\begin{equation}
\begin{split}
\widetilde{R}_{12}(x_1-x_2,\lambda+h) \, & 
\widetilde{L}_1(x_1,\lambda) \, 
\widetilde{L}_2(x_2,\lambda+h^{(1)}) = \\
& \widetilde{L}_2(x_2,\lambda) \, 
\widetilde{L}_1(x_1,\lambda+h^{(2)}) \, 
\widetilde{R}_{12}(x_1-x_2,\lambda)
\end{split}
\end{equation}
Denoting by $\{h_{i}\}$ a basis of $\fH$ and $\{h_{i}^{\vee}\}$ the dual 
basis, the notation $\lambda+h^{(k)}$ means $\sum_{i} (\lambda_{i} + 
{h_{i}}^{(k)}) h_{i}^{\vee}$ if $\lambda = \sum_{i} \lambda_{i} \, 
h_{i}^{\vee}$ and the superscript $(k)$ labels the space where $h$ acts.

\section{Elliptic algebras}
\setcounter{equation}{0}

\subsection{The quantum affine elliptic algebras of vertex type}

The quantum affine elliptic algebra of vertex type $\elpa{2}$ was 
introduced in \cite{FIJKMY,FIJKMY95} and generalized to the $sl_{N}$ case 
in \cite{AFRS3,JKOS}. The status of this algebra was elucidated in 
\cite{JKOS}, where it was shown that $\elpa{N}$ was a QTQHA obtained from 
$\cU_{q}(\widehat{sl}_{N})$ by a suitable Drinfel'd twist.

\subsubsection{R-matrix of $\elpa{2}$}

The quantum affine elliptic algebra $\elpa{2}$ is defined in the RLL 
formalism. The entries of the R-matrix of the quantum affine elliptic 
algebra $\elpa{2}$ are given by the Boltzmann weights of the eight vertex 
model, solved by Baxter \cite{Bax72,Bax82}. Explicitly, this R-matrix reads
\begin{equation}
\label{eq:RAqp}
R[\elpa{2}](z,q,p) = \rho(z,q,p) \left(
\begin{array}{cccc}
a(z) & 0 & 0 & d(z) \\
0 & b(z) & c(z) & 0 \\
0 & c(z) & b(z) & 0 \\
d(z) & 0 & 0 & a(z) \\
\end{array}
\right)
\end{equation}
where
\begin{equation}
\label{eq:abcdAqp}
\begin{array}{lcl}
\displaystyle a(z) = z^{-1} \; \frac{\Theta_{p^2}(q^2z^2) \; 
\Theta_{p^2}(pq^2)} {\Theta_{p^2}(pq^2z^2) \; \Theta_{p^2}(q^2)} 
&\qquad& \displaystyle d(z) = -\frac{p^{1/2}}{q z^{2}} \; 
\frac{\Theta_{p^2}(z^2) \; \Theta_{p^2}(q^2z^2)} 
{\Theta_{p^2}(pz^2) \; \Theta_{p^2}(pq^2z^2)} \\[5mm] \displaystyle 
b(z) = qz^{-1} \; \frac{\Theta_{p^2}(z^2) \; \Theta_{p^2}(pq^2)} 
{\Theta_{p^2}(pz^2) \; \Theta_{p^2}(q^2)} &\qquad& \displaystyle 
c(z) = 1 \\
\end{array}
\end{equation}
with $\Theta_{a}(x) = (x;a)_{\infty}\, (ax^{-1};a)_{\infty} \, 
(a;a)_{\infty}$. \\
The normalization factor $\mu(x)$ is chosen as follows \cite{JKM}:
\begin{equation}
\label{eq:normAqp}
{\rho(z)} = \frac{(p^2;p^2)_\infty} {(p;p)_\infty^2} \; 
\frac{\Theta_{p^2}(pz^2)\Theta_{p^2}(q^2)} {\Theta_{p^2}(q^2z^2)} \; 
\frac{(q^4z^{-2};p,q^4)_\infty \; (q^2z^2;p,q^4)_\infty \; 
(pz^{-2};p,q^4)_\infty \; (pq^2z^2;p,q^4)_\infty} 
{(q^4z^2;p,q^4)_\infty \; (q^2z^{-2};p,q^4)_\infty \; 
(pz^2;p,q^4)_\infty \; (pq^2z^{-2};p,q^4)_\infty}
\end{equation}

\begin{proposition}
The R-matrix (\ref{eq:RAqp}) has the following properties:
\begin{align}
& \text{Yang--Baxter equation:} && R_{12}(z) \, R_{13}(w) \, 
R_{23}(w/z) = R_{23}(w/z) \, R_{13}(w) \, R_{12}(z) &&&&&& 
\\[6pt] 
& \text{unitarity:} && R_{12}(z) \, R_{21}(z^{-1}) = 1 &&&&&& 
\\[6pt] 
& \text{crossing symmetry:} && R_{21}(z^{-1})^{t_1} = (\sigma^1 
\otimes 1) R_{12}(-q^{-1}z) (\sigma^1 \otimes 1) &&&&&& \\[6pt] 
& \text{antisymmetry:} && R_{12}(-z) = - (\sigma^3 \otimes 1) 
R_{12}(z) (\sigma^3 \otimes 1) &&&&&&
\end{align}
where $\sigma^1,\sigma^2,\sigma^3$ are the $2 \!\times\! 2$ Pauli 
matrices and $t_i$ denotes the transposition in the space $i$.
\end{proposition}
The proof is straightforward by direct calculation.

\subsubsection{RLL relations for $\elpa{2}$}

For the definition of the quantum affine elliptic algebra $\elpa{2}$, one 
needs to use a slightly modified R-matrix $\widetilde{R}_{12}(z)$, which 
differs from (\ref{eq:RAqp}) by a suitable normalization factor:
\begin{equation}
\widetilde{R}_{12}(z) = \tau(q^{1/2}z^{-1}) R_{12}(z)
\end{equation}
where the factor $\tau(z)$ is given by
\begin{equation}
\tau(z) = z^{-1} \frac{(qz^2;q^4)_\infty \; (q^3z^{-2};q^4)_\infty} 
{(qz^{-2};q^4)_\infty \; (q^3z^2;q^4)_\infty}
\end{equation}
The function $\tau$ is periodic, $\tau(z) = \tau(zq^2)$, and 
$\widetilde{R}_{12}(z)$ obeys a quasi-periodicity property:
\begin{equation}
\widetilde{R}_{12}(-p^{\frac {1}{2}}z)= (\sigma^1 \otimes 1) \left 
(\widetilde{R}_{21}(z^{-1})\right )^{-1} (\sigma^1 \otimes 1)
\end{equation}
The crossing symmetry and the unitarity properties of $R_{12}$ then allow 
one to exchange inversion and transposition for the matrix 
$\widetilde{R}_{12}$ as:
\begin{equation}
\Big( \widetilde{R}_{12}(z)^{t_2} \Big)^{-1} = \Big( 
\widetilde{R}_{12}(q^2z)^{-1} \Big)^{t_2}
\end{equation}

\medskip

The quantum affine elliptic algebra $\elpa{2}$ is defined as an algebra of 
operators
\begin{equation}
L_{ij}(z) = \sum_{n\in\ZZ} L_{ij}(n) \, z^n = \left(
\begin{array}{cc}
L_{++}(z) & L_{+-}(z) \\
L_{-+}(z) & L_{--}(z) \\
\end{array}
\right)
\end{equation}
where $i,j \in \ZZ_2 \equiv \ZZ/2\ZZ$, encapsulated into a $2 \!\times\! 2$ 
matrix, and the functions the $L_{++}$ and $L_{--}$ are even while $L_{+-}$ 
and $L_{-+}$ are odd in the variable $z$. One first defines 
$\cA_{q,p}(\widehat{gl}_{2})$ by imposing the following relations on 
$L(z)$:
\begin{equation}
\label{eq:RLLAqp}
\widetilde{R}_{12}(z/w,q,p) \, L_{1}(z) \, L_{2}(w) = L_{2}(w) \, 
L_{1}(z) \, \widetilde{R}^{*}_{12}(z/w,q,p)
\end{equation}
with $L_1(z) = L(z) \otimes 1$, $L_2(z) = 1 \otimes L(z)$ and 
$\widetilde{R}^{*}_{12}(z,q,p) = \widetilde{R}_{12}(z,q,p^*=pq^{-2c})$. \\
The quantum determinant of $L(z)$ given by
\begin{equation}
\qdet L(z) = L_{++}(q^{-1}z) L_{--}(z) - L_{-+}(q^{-1}z) L_{+-}(z)
\end{equation}
is in the center of ${\cal A}_{q,p}(\widehat{gl}_{2})$. It can be factored 
out, and set to the value $q^{\frac c2}$ ($c$ being the central charge) so 
as to get
\begin{equation}
\elpa{2} = \cA_{q,p}(\widehat{gl}_{2})/ \big\langle \qdet L - 
q^{\frac{c}{2}} \big\rangle
\end{equation}
Note that it may be useful to introduce the following two matrices:
\begin{equation}
L^+(z) \equiv L(q^{\frac c2}z) \qquad \text{and} \qquad L^-(z) \equiv 
\sigma^{1} L(-p^{\frac 12}z) \sigma^{1}
\end{equation}
They obey coupled exchange relations following from (\ref{eq:RLLAqp}) and 
periodicity/unitarity properties of the matrices $\widetilde{R}_{12}$ and 
$\widetilde{R}^{*}_{12}$:
\begin{equation}
\begin{split}
\widetilde{R}_{12}(z/w) \, L^{\pm}_1(z) \, L^{\pm}_2(w) &= 
L^{\pm}_2(w) \, L^{\pm}_1(z) \, \widetilde{R}^{*}_{12}(z/w) \\
\widetilde{R}_{12}(q^{\frac c2}z/w) \, L^+_1(z) \, L^-_2(w) &= 
L^-_2(w) \, L^+_1(z) \, \widetilde{R}^{*}_{12}(q^{-\frac c2}z/w)
\end{split}
\end{equation}

\begin{remark}
A further renormalization of the modes ${\bar L}_n = \big( -p^\half 
\big)^{\max(n,0)} L_n$ is required in order to get the trigonometric 
quantum enveloping algebra $\cU_q(\widehat{sl}(2)_c)$ from $\elpa{2}$ 
when $p \to 0$. It ensures the vanishing of half of the degrees of 
freedom in $L^{\pm}(z)$ and decouples completely $L^+(z)$ from 
$L^-(z)$, thereby keeping the same overall dimensionality \cite{JKM}. 
It is therefore a non-trivial, ``discontinuous'' procedure. \finrmq
\end{remark}

Let us stress that at this level the Hopf structure of $\elpa{2}$ remains 
undetermined. Indeed, the coproduct of the $L(z)$ generators \emph{cannot} 
be defined by $\Delta(L) = L \dot\otimes L$, since the R-matrices appearing 
in the l.h.s. and in the r.h.s. of eq. (\ref{eq:RLLAqp}) differ, the values 
of the elliptic nomes $p$ and $p^*$ of the elliptic functions 
(\ref{eq:abcdAqp}) entering in the definition of the entries of 
(\ref{eq:RAqp}) and of the normalization factor (\ref{eq:normAqp}) being 
different. The answer to this fundamental question will be given in section 
\ref{sect:quasiAqp}.

\subsubsection{Generalization to $\elpa{N}$}

The R-matrix of the quantum affine elliptic algebra $\elpa{N}$, associated 
to the $\ZZ_{N}$-vertex model, is defined as follows \cite{Bel81,ChCh}:
\begin{equation}
\label{eq:RAqpN}
R[\elpa{N}](z,q,p) = \rho(z,q,p) 
\sum_{(\alpha_1,\alpha_2)\in\ZZ_N\times\ZZ_N} 
W_{(\alpha_1,\alpha_2)}(\xi,\zeta,\tau) \,\, I_{(\alpha_1,\alpha_2)} 
\otimes I_{(\alpha_1,\alpha_2)}^{-1}
\end{equation}
where the variables $z,q,p$ are related to the variables $\xi,\zeta,\tau$ 
by
\begin{equation*}
z=e^{i\pi\xi} \,,\qquad q=e^{i\pi\zeta} \,,\qquad p=e^{2i\pi\tau}
\end{equation*}
The functions $W_{(\alpha_1,\alpha_2)}$ are given by
\begin{equation}
W_{(\alpha_1,\alpha_2)}(\xi,\zeta,\tau) = 
\frac{\vartheta\car{\sfrac{1}{2}+\alpha_1/N} 
{\half+\alpha_2/N}(\xi+\zeta/N,\tau)} 
{N\vartheta\car{\sfrac{1}{2}+\alpha_1/N} 
{\half+\alpha_2/N}(\zeta/N,\tau)}
\end{equation}
where the Jacobi theta functions with rational characteristics 
$(\gamma_1,\gamma_2) \in \sfrac{1}{N} \ZZ \times \sfrac{1}{N} \ZZ$ are 
defined by
\begin{equation}
\vartheta\car{\gamma_1}{\gamma_2}(\xi,\tau) = \sum_{m \in \ZZ} 
\exp\Big(i\pi(m+\gamma_1)^2\tau + 2i\pi(m+\gamma_1)(\xi+\gamma_2) \Big)
\end{equation}
and the matrices $I_{(\alpha_1,\alpha_2)}$ by
\begin{equation}
I_{(\alpha_1,\alpha_2)} = g^{\alpha_2} \, h^{\alpha_1}
\end{equation}
the $N \!\times\! N$ matrices $g$ and $h$ being given by $g_{ij} = 
\omega^i\delta_{ij}$ and $h_{ij} = \delta_{i+1,j}$, the addition of indices 
being understood modulo $N$ and $\omega = e^{2i\pi/N}$. \\
Finally, the normalization factor $\rho(z,q,p)$ is given by
\begin{equation}
\rho(z,q,p) = \frac{(q^{2N}z^{-2};p,q^{2N})_\infty \; 
(q^2z^2;p,q^{2N})_\infty \, (pz^{-2};p,q^{2N})_\infty \; 
(pq^{2N-2}z^2;p,q^{2N})_\infty} {(q^{2N}z^2;p,q^{2N})_\infty \; 
(q^2z^{-2};p,q^{2N})_\infty \, (pz^2;p,q^{2N})_\infty \; 
(pq^{2N-2}z^{-2};p,q^{2N})_\infty}
\end{equation}
The matrix (\ref{eq:RAqpN}) is $\ZZ_N$-symmetric, that is
\begin{equation}
R_{a+s\,,\,b+s}^{c+s\,,\,d+s} = R_{a\,,\,b}^{c\,,\,d}
\end{equation}
for any indices $a,b,c,d,s \in \ZZ_N$.

\medskip

\begin{proposition}
The R-matrix $\widehat{R}(z,q,p)$ similar~%
\footnote{Two R-matrices $R(z)$ and $R'(z)$ are said to be similar if 
$R'(z) = (A \otimes A) R(z) (A \otimes A)^{-1}$ where $A$ is a 
non-degenerate operator in the representation space $V$.}
to $R(z,q,p)$ given by
\begin{equation}
\widehat{R}(z,q,p) = (g^{\frac 12} \otimes g^{\frac 12}) R(z,q,p) 
(g^{-\frac 12} \otimes g^{-\frac 12})
\end{equation}
satisfies the following properties:
\begin{align}
& \text{Yang--Baxter equation:} && \widehat{R}_{12}(z) \, 
\widehat{R}_{13}(w) \, \widehat{R}_{23}(w/z) = 
\widehat{R}_{23}(w/z) \, \widehat{R}_{13}(w) \, \widehat{R}_{12}(z) 
&&&&&& \\[6pt] 
& \text{unitarity:} && \widehat{R}_{12}(z) \, 
\widehat{R}_{21}(z^{-1}) = 1 &&&&&& \\[6pt] 
& \text{crossing symmetry:} && \widehat{R}_{12}(z)^{t_2} \, 
\widehat{R}_{21}(q^{-N}z^{-1})^{t_2} = 1 &&&&&& \\[6pt] 
& \text{antisymmetry:} && \widehat{R}_{12}(-z) = \omega \, (g^{-1} 
\otimes 1) \, \widehat{R}_{12}(z) \, (g \otimes 1) &&&&&&
\end{align}
\end{proposition}
One also introduces a modified R-matrix $\widetilde{R}_{12}$
\begin{equation}
\widetilde{R}_{12}(z,q,p) = \tau_{N}(q^{1/2}z^{-1}) \, 
\widehat{R}_{12}(z,q,p)
\end{equation}
where
\begin{equation}
\tau_N(z) = z^{\frac{2}{N}-2} \, 
\frac{\Theta_{q^{2N}}(qz^2)}{\Theta_{q^{2N}}(qz^{-2})}
\end{equation}
The R-matrix $\widetilde{R}_{12}$ obeys a quasi-periodicity property:
\begin{equation}
\widetilde{R}_{12}(-p^{\frac{1}{2}}z) = 
(g^{\frac{1}{2}}hg^{\frac{1}{2}} \otimes 1)^{-1} \left 
(\widetilde{R}_{21}(z^{-1})\right )^{-1} 
(g^{\frac{1}{2}}hg^{\frac{1}{2}} \otimes 1)
\end{equation}
The function $\tau_N(z)$ is periodic with period $q^N$, $\tau_N(q^Nz) = 
\tau_N(z)$, and satisfies $\tau_N(z^{-1}) = \tau_N(z)^{-1}$.

\medskip

Consider the algebra over $\CC$ generated by the operators $L_{ij}(z) = 
\sum_{n \in \ZZ} L_{ij}(n) \, z^{n}$ where $i,j \in \ZZ_{N}$, encapsulated 
into a $N \!\times\! N$ matrix:
\begin{equation}
L(z) = \left(
\begin{array}{ccc}
L_{11}(z) & \cdots & L_{1N}(z) \\
\vdots && \vdots \\
L_{N1}(z) & \cdots & L_{NN}(z) \\
\end{array}
\right)
\end{equation}
The algebra ${\cal A}_{q,p}(\widehat{gl}(N))$ is defined by imposing the 
following relations among the $L(z)$ generators:
\begin{equation}
\widetilde{R}_{12}(z/w,q,p) \, L_1(z) \, L_2(w) = L_2(w) \, L_1(z) \, 
\widetilde{R}_{12}^{*}(z/w,q,p)
\end{equation}
with $L_1(z) = L(z) \otimes 1$, $L_2(z) = 1 \otimes L(z)$ and 
$\widetilde{R}^{*}_{12}(z,q,p) = \widetilde{R}_{12}(z,q,p^*=pq^{-2c})$. \\
The quantum determinant of $L(z)$ given by
\begin{equation}
\qdet L(z) = \sum_{\sigma\in{\mathfrak S}_N} \varepsilon(\sigma) 
\prod_{i=1}^N L_{i,\sigma(i)}(z q^{i-N-1})
\end{equation}
($\varepsilon(\sigma)$ being the signature of the permutation $\sigma$) 
lies in the center of ${\cal A}_{q,p}(\widehat{gl}_{N})$. It can be 
factored out, and set to the value $q^{\frac c2}$ ($c$ being the central 
charge) so as to get
\begin{equation}
\elpa{N} = \cA_{q,p}(\widehat{gl}_{N})/ \big\langle \qdet L - 
q^{\frac{c}{2}} \big\rangle
\end{equation}

\subsubsection{Quasi-Hopf algebra structure for $\elpa{N}$}
\label{sect:quasiAqp}

Consider the affine Lie algebra $\widehat{sl}_{N}$ with Cartan basis $\fH = 
\{h_{0},\ldots,h_{N-1},d\}$ and dual basis $\fH^* = 
\{\Lambda_{0},\ldots,\Lambda_{N-1},c\}$ where $\Lambda_{i}$ are the 
fundamental weights, $d$ the derivation and $c$ the central charge. Let 
$\tau$ be the automorphism of order $N$ related to the cyclic symmetry of 
the extended Dynkin diagram of $\widehat{sl}_{N}$: $\tau(x_{i}) = x_{i+1 
\mod N}$ where the $x_{i}$'s denote the Chevalley generators of 
$\widehat{sl}_{N}$. One defines
\begin{equation}
\varphi = \tau \circ \Ad \Big( q^{2(r+c)\rho/N} \Big)
\end{equation}
where $\rho = \sum_{i=0}^{N-1} \Lambda_{i}$ gives the principal grading and 
$r$ is some complex number. \\
One defines
\begin{equation}
\label{eq:twistAqp}
\cF(r) = \prod_{k \ge 1}^{\leftarrow} \cF_{k}(r) = \prod_{k \ge 
1}^{\leftarrow} (\varphi^k \otimes \id) (\widehat{\cR}^{-1}) \equiv 
\lim_{n \to \infty} \prod_{nN \ge k \ge 1}^{\leftarrow} (\varphi^k 
\otimes \id) (\widehat{\cR}^{-1})
\end{equation}
where
\begin{equation}
\widehat{\cR} = q^{T} \cR[\elpa{N}] \quad \text{with} \quad T = 
\frac{1}{N} \, \Big( \rho \otimes c + c \otimes \rho - 
\frac{N^2-1}{12}\, c \otimes c \Big)
\end{equation}
The arrow means that the product has to be done to the left, i.e. $\ldots 
\cF_{3}(r) \, \cF_{2}(r) \, \cF_{1}(r)$.

In ref. \cite{JKOS} the following theorem was proved:
\begin{theorem}
The vertex-type Drinfel'd twist $\displaystyle \cF(r) = \prod_{k \ge 
1}^{\leftarrow} (\varphi^k \otimes \id) (\widehat{\cR}^{-1})$ satisfies 
the shifted cocycle condition (hence it is a GNF twist)
\begin{equation*}
\cF_{12}(r) \, (\Delta \otimes \id) \cF(r) = \cF_{23}(r + c^{(1)}) 
\, (\id \otimes \Delta) \cF(r)
\end{equation*}
Moreover, one has $(\id \otimes \epsilon) \cF(r) = (\epsilon \otimes 
\id) \cF(r) = 1$. \\
The quantum affine elliptic algebra $\elpa{N}$ (where $p=q^{2r}$) is 
a quasi-triangular quasi-Hopf algebra with the universal R-matrix 
$\widetilde{\cR}(r) = \cF_{21}(r) \, \cR_{12} \, \cF_{12}^{-1}(r)$, 
where $\cR$ is the universal R-matrix of $\cU_{q}(\widehat{sl}_{N})$. 
$\widetilde{\cR}(r)$ satisfies the shifted Yang--Baxter equation:
\begin{equation*}
\widetilde{\cR}_{12}(r + c^{(3)}) \, \widetilde{\cR}_{13}(r) \, 
\widetilde{\cR}_{23}(r + c^{(1)}) = \widetilde{\cR}_{23}(r) \, 
\widetilde{\cR}_{13}(r + c^{(2)}) \, \widetilde{\cR}_{12}(r)
\end{equation*}
\end{theorem}
Sketch of the proof (see \cite{JKOS} for the complete proof). 
$\cF_k(r)$ satisfy the following properties
\begin{align}
\label{eq:DeltaidF}
(\Delta\otimes\id)(\cF_k(r)) &= \cF_k^{(23)}(r+c^{(1)}) 
\cF_k^{(13)}\Big(r+c^{(2)}-\frac{1}{2k}\;c^{(2)}\Big) \\
(\id\otimes\Delta)(\cF_k(r)) &= \cF_k^{(12)}(r) 
\cF_k^{(13)}\Big(r+\frac{1}{2k}\;c^{(2)}\Big)
\end{align}
and
\begin{equation}
\label{eq:YBEFk}
\cF_k^{(12)}(r) \cF_{k+l}^{(13)}\Big(r+\frac{l+\frac12}{k+l}\, 
c^{(2)}\Big) \cF_l^{(23)}(r+c^{(1)}) = \cF_l^{(23)}(r+c^{(1)}) 
\cF_{k+l}^{(13)}\Big(r+\frac{l-\frac12}{k+l}\, c^{(2)}\Big) 
\cF_k^{(12)}(r)
\end{equation}
Using equation (\ref{eq:YBEFk}), one can prove by induction the
following relation
\begin{align}
\label{eq:cocinduc}
\prod_{l\ge k\ge 1}^{\leftarrow} \cF_k^{(23)}(r+c^{(1)}) \; 
(\id\otimes\Delta)\big(\cF(r)\big) &= \prod_{k\ge 1}^{\leftarrow} 
\cF_k^{(12)}(r) \; 
\cF_{k+l}^{(13)}\Big(r+\frac{l+\frac12}{k+l}\,c^{(2)}\Big) \nonumber \\
&\times \prod_{l\ge k\ge 1}^{\leftarrow} \cF_k^{(23)}(r+c^{(1)}) \; 
\cF_{k}^{(13)}\Big(r+c^{(2)}-\frac{1}{2k}\, c^{(2)}\Big)
\end{align}
Letting then $l \to \infty$ and taking into account (\ref{eq:DeltaidF}), 
one recovers the shifted cocycle condition. \\
Finally, $(\id \otimes \epsilon) \cF(r) = (\epsilon \otimes \id) \cF(r) = 
1$ follows obviously from (\ref{eq:Reps}). 
\finproof

\medskip

Let us notice that the vertex-type Drinfel'd twist is the (unique) solution 
of the linear difference equation
\begin{equation}
\cF_{12}(r) = \Ad (\varphi^{-1} \otimes \id) (\cF_{12}(r)) \cdot 
\widehat{\cR}^{-1}
\end{equation}
such that $(\id \otimes \epsilon) \cF(r) = (\epsilon \otimes \id) \cF(r) = 
1$.

\medskip

Example. Consider the two-dimensional evaluation representation of 
$\widehat{sl}_2$ with evaluation parameter $z$, in the principal gradation
\begin{equation}
\begin{split}
& \pi_z(e_1) = z e_{12} \;, \qquad \pi_z(f_1) = z^{-1} e_{21} \;, 
\qquad \pi_z(h_1) = e_{11} - e_{22} \\
& \pi_z(e_0) = z e_{21} \;, \qquad \pi_z(f_0) = z^{-1} e_{12} \;, 
\qquad \pi_z(h_0) = e_{22} - e_{11}
\end{split}
\end{equation}
The R-matrix of $\cU_{q}(\widehat{sl}_2)$ is this representation is given 
by [compare with (\ref{eq:ruq})]
\begin{equation}
\label{eq:ruqa}
R[\cU_{q}(\widehat{sl}_2)](z) = \rho(z^2) \left(
\begin{array}{cccc}
1 & 0 & 0 & 0 \\
0 & \displaystyle \frac{q(1-z^2)}{1-q^2z^2} & \displaystyle 
\frac{z(1-q^2)}{1-q^2z^2} & 0 \\[12pt] 0 & \displaystyle 
\frac{z(1-q^2)}{1-q^2z^2} & \displaystyle \frac{q(1-z^2)}{1-q^2z^2} 
& 0 \\
0 & 0 & 0 & 1 \\
\end{array}
\right)
\end{equation}
where the normalization factor $\rho(z)$ is given by (\ref{eq:rhouq}). \\
One sets $F(z;p) = \Ad(z^\rho \otimes \id) \cF(r)$ ($p$ and $r$ being 
related as above). Using the expression (\ref{eq:ruqa}) of the R-matrix of 
$\cU_{q}(\widehat{sl}_2)$ and the definition (\ref{eq:twistAqp}) of the 
Drinfel'd twist, one gets \cite{Fro97b,Fro97a}
\begin{equation}
F(z;p) = \rho_{F}(z;p) \left(
\begin{array}{cccc}
a_{F}(z) & 0 & 0 & d_{F}(z) \\
0 & b_{F}(z) & c_{F}(z) & 0 \\
0 & c_{F}(z) & b_{F}(z) & 0 \\
d_{F}(z) & 0 & 0 & a_{F}(z) \\
\end{array}
\right)
\end{equation}
where
\begin{align}
& a_{F}(z) \pm d_{F}(z) = \frac{(\mp p^{1/2}qz;p)_{\infty}} {(\mp 
p^{1/2}q^{-1}z;p)_{\infty}} \\
& b_{F}(z) \pm c_{F}(z) = \frac{(\mp pqz;p)_{\infty}} {(\mp 
pq^{-1}z;p)_{\infty}}
\end{align}
and the normalization factor is
\begin{equation}
\rho_{F}(z;p) = \frac{(pz^2;p,q^4)_{\infty} \; 
(pq^4z^2;p,q^4)_{\infty}} {(pq^2z^2;p,q^4)_{\infty}^2}
\end{equation}
Now, computing the product $F_{21}(z^{-1};p) \, 
R[\cU_{q}(\widehat{sl}_2)](z) \, F_{12}(z;p)^{-1}$, one recovers the 
expression (\ref{eq:RAqp}) of the R-matrix of $\elpa{2}$.

\subsection{The quantum affine elliptic algebras of face type}

\subsubsection{R-matrix of the affine elliptic algebra $\elpb{2}$}

The face-type R-matrices depend on the extra parameters $\lambda$ belonging 
to the dual of the Cartan algebra of the underlying algebra. Let 
$\{h,c,d\}$ be a basis of the Cartan subalgebra of $\widehat{sl}_{2}$. If 
$r,s,s'$ are complex numbers, we set $\lambda = \half (s+1)h + s'c + 
(r+2)d$. The elliptic parameter $p$ and the dynamical parameter $w$ are 
related to the deformation parameter $q$ by $p=q^{2r}$, $w=q^{2s}$. The 
R-matrix of $\elpb{2}$ is \cite{Fel94,JKOS}
\begin{equation}
\label{eq:Relpb}
R[\elpb{2}](z,\lambda) = \rho(z;p) \left(
\begin{array}{cccc}
1 & 0 & 0 & 0 \\
0 & b(z) & c(z) & 0 \\
0 & {\bar c}(z) & {\bar b}(z) & 0 \\
0 & 0 & 0 & 1 \\
\end{array}
\right)
\end{equation}
where
\begin{equation}
\begin{split}
b(z) &= q \frac{(pw^{-1}q^2;p)_{\infty} \; 
(pw^{-1}q^{-2};p)_{\infty}} {(pw^{-1};p)_{\infty}^2} \; 
\frac{\Theta_{p}(z)}{\Theta_{p}(q^2z)} \\
{\bar b}(z) &= q \frac{(wq^2;p)_{\infty} \; (wq^{-2};p)_{\infty}} 
{(w;p)_{\infty}^2} \; \frac{\Theta_{p}(z)}{\Theta_{p}(q^2z)} \\
c(z) &= \frac{\Theta_{p}(q^2)}{\Theta_{p}(w)} \; 
\frac{\Theta_{p}(wz)}{\Theta_{p}(q^2z)} \\
{\bar c}(z) &= \frac{\Theta_{p}(q^2)}{\Theta_{p}(w^{-1})} \; 
\frac{\Theta_{p}(w^{-1}z)}{\Theta_{p}(q^2z)}
\end{split}
\end{equation}
The normalization factor is
\begin{equation}
\rho(z;p) = q^{-1/2} \frac{(q^2z;p,q^4)_{\infty}^2} {(z;p,q^4)_{\infty} 
\; (q^4z;p,q^4)_{\infty}} \; \frac{(pz^{-1};p,q^4)_{\infty} \; 
(pq^4z^{-1};p,q^4)_{\infty}} {(pq^2z^{-1};p,q^4)_{\infty}^2}
\end{equation}
The elliptic algebra $\elpb{2}$ is then defined by 
\begin{equation}
R_{12}(z_1/z_2,\lambda+h) \, L_1(z_1,\lambda) \, 
L_2(z_2,\lambda+h^{(1)}) = L_2(z_2,\lambda) \, L_1(z_1,\lambda+h^{(2)}) 
\, R_{12}(z_1/z_2,\lambda)
\end{equation}

\subsubsection{Generalization to $\elpb{N}$}

Let $\fH$ be the Cartan subalgebra of the affine Lie algebra 
$\widehat{sl}_{N}$, with basis $\{ h_{i} \}$ and dual basis 
$\{h_{i}^{\vee}\}$, and $\Pi^0 = \{\alpha_{0},\ldots,\alpha_{N-1}\}$ be the 
corresponding simple root system of $\widehat{sl}_{N}$. We set $\lambda = 
\rho + s_{1} h_{1} + \ldots + s_{N-1} h_{N-1} + (r+N)d + s'c$ where $\rho 
\in \fH$ is such that $(\rho,\alpha_{i})=1$. The R-matrix of the quantum 
affine elliptic algebra $\elpb{N}$ takes the following form ($1 \le a,b \le 
N$) \cite{Fel94}:
\begin{align}
R[\elpb{N}](z,\lambda) = \rho(z;p) &\Big( \sum_{a} E_{aa} \otimes 
E_{aa} + \sum_{a \ne b} q^2 
\frac{\Theta_{p}(q^{-2}w_{ab})}{\Theta_{p}(w_{ab})} \; 
\frac{\Theta_{p}(z)}{\Theta_{p}(q^2z)} \; E_{aa} \otimes E_{bb} 
\nonumber \\
&+ \sum_{a \ne b} \frac{\Theta_{p}(q^2)}{\Theta_{p}(w_{ab})} \; 
\frac{\Theta_{p}(w_{ab}z)}{\Theta_{p}(q^2z)} \; E_{ab} \otimes E_{ba} 
\Big)
\label{eq:RBqpln}
\end{align}
where the $E_{ab}$ are the elementary $N \!\times\! N$ matrices with 
entries ${(E_{ab})}_{i}^{j} = \delta_{ai} \delta_{jb}$, $w_{ab} = 
q^{x_{a}-x_{b}}$, $x_{a} = 2s_{a} - 2s_{a-1}$ with by convention $s_{0} = 
s_{N} = 0$. \\
The normalization factor $\rho(z)$ is given by
\begin{equation}
\rho(z;p) = q^{-\frac{N-1}{N}} \; \frac{(q^2z;q^{2N},p)_{\infty} \; 
(q^{2N-2}z;q^{2N},p)_{\infty}} {(z;q^{2N},p)_{\infty} \; 
(q^{2N}z;q^{2N},p)_{\infty}} \; \frac{(pz^{-1};q^{2N},p)_{\infty} \; 
(pq^{2N}z^{-1};q^{2N},p)_{\infty}} {(pq^2z^{-1};q^{2N},p)_{\infty} \; 
(pq^{2N-2}z^{-1};q^{2N},p)_{\infty}}
\end{equation}
The R-matrix of $\elpb{N}$ used in \cite{JKOS} is obtained from 
(\ref{eq:RBqpln}) by a similarity transformation. It reads
\begin{align}
R[\elpb{N}](z,\lambda) = \rho(z;p) &\Big( \sum_{a} E_{aa} \otimes 
E_{aa} + \sum_{a \ne b} \frac{\Theta_{p}(q^{2})}{\Theta_{p}(w_{ab})} \; 
\frac{\Theta_{p}(w_{ab}z)}{\Theta_{p}(q^2z)} \; E_{ab} \otimes E_{ba} 
\nonumber \\
&+ \sum_{a < b} q \; \frac{(pw_{ab}^{-1}q^2;p)_{\infty} \; 
(pw_{ab}^{-1}q^{-2};p)_{\infty}} {(pw_{ab}^{-1};p)_{\infty}^2} \; 
\frac{\Theta_{p}(z)}{\Theta_{p}(q^2z)} \; E_{aa} \otimes E_{bb} 
\nonumber \\
&+ \sum_{a > b} q \; \frac{(w_{ab}^{-1}q^2;p)_{\infty} \; 
(w_{ab}^{-1}q^{-2};p)_{\infty}} {(w_{ab}^{-1};p)_{\infty}^2} \; 
\frac{\Theta_{p}(z)}{\Theta_{p}(q^2z)} \; E_{aa} \otimes E_{bb} \Big)
\end{align}

\subsubsection{Quasi-Hopf algebra structure for $\elpb{N}$}

One considers the following automorphism of $\cU_{q}(\widehat{sl}_{N})$, 
where $\lambda \in \fH$
\begin{equation}
\varphi = \Ad \Big( q^{\sum_{i} h_{i} h_{i}^{\vee} + 2(\lambda-\rho)} 
\Big)
\end{equation}
One defines the face-type Drinfel'd twist by
\begin{equation}
\label{eq:twistBqpl}
\cF(\lambda) = \prod_{k \ge 1}^{\leftarrow} \cF_{k}(\lambda) = \prod_{k 
\ge 1}^{\leftarrow} (\varphi^k \otimes \id) (\widehat{\cR}^{-1})
\end{equation}
where $\widehat{\cR} = q^{T} \cR$ and $T = \sum_{i} h_{i} \otimes 
h_{i}^{\vee}$.

In ref. \cite{JKOS} the following theorem was proved:
\begin{theorem}
The face-type Drinfel'd twist $\displaystyle \cF(\lambda) = \prod_{k 
\ge 1}^{\leftarrow} (\varphi^k \otimes \id) (\widehat{\cR}^{-1})$ 
satisfies the shifted cocycle condition (hence it is a GNF twist)
\begin{equation*}
\cF_{12}(\lambda) \, (\Delta \otimes \id) \cF(\lambda) = 
\cF_{23}(\lambda + h^{(1)}) \, (\id \otimes \Delta) \cF(\lambda)
\end{equation*}
Moreover, one has $(\id \otimes \epsilon) \cF(\lambda) = (\epsilon 
\otimes \id) \cF(\lambda) = 1$. \\
The quantum affine elliptic algebra $\elpb{N}$ is a quasi-triangular 
quasi-Hopf algebra with the universal R-matrix 
$\widetilde{\cR}(\lambda) = \cF_{21}(\lambda) \, \cR_{12} \, 
\cF_{12}^{-1}(\lambda)$, where $\cR$ is the universal R-matrix of 
$\cU_{q}(\widehat{sl}_{N})$. $\widetilde{\cR}(\lambda)$ satisfies the 
dynamical Yang--Baxter equation:
\begin{equation*}
\widetilde{\cR}_{12}(\lambda + h^{(3)}) \, 
\widetilde{\cR}_{13}(\lambda) \, \widetilde{\cR}_{23}(\lambda + 
h^{(1)}) = \widetilde{\cR}_{23}(\lambda) \, 
\widetilde{\cR}_{13}(\lambda + h^{(2)}) \, 
\widetilde{\cR}_{12}(\lambda)
\end{equation*}
\end{theorem}
The proof follows the same lines as for the vertex case. In particular the 
equations (\ref{eq:DeltaidF}) to (\ref{eq:cocinduc}) have the same form 
with the replacement $r \to \lambda$ and $c \to h$.

\medskip

Let us emphasize that the face-type Drinfel'd twist is the (unique) 
solution of the linear difference equation
\begin{equation}
\label{eq:eqlinB}
\cF_{12}(\lambda) = \Ad (\varphi^{-1} \otimes \id) (\cF_{12}(\lambda)) 
\cdot \widehat{\cR}^{-1}
\end{equation}
such that $(\id \otimes \epsilon) \cF(\lambda) = (\epsilon \otimes \id) 
\cF(\lambda) = 1$. The importance of this equation will be pointed out in 
the example just below.

\medskip

Example. In the two-dimensional evaluation representation for 
$\widehat{sl}_2$ in the homogeneous gradation eq. (\ref{eq:evalrepr}), one 
sets $F(z;p,w) = \Ad(z^\rho \otimes \id) \cF(\lambda)$. Unlike the vertex 
case, using the expression (\ref{eq:ruq}) of the R-matrix of 
$\cU_{q}(\widehat{sl}_2)$, the infinite product (\ref{eq:twistBqpl}) 
defining the face-type Drinfel'd twist cannot be evaluated directly. In 
fact, equation (\ref{eq:eqlinB}) is the only way to obtain the solution. In 
the two-dimensional evaluation representation under consideration, the 
linear difference equation becomes
\begin{equation}
F_{12}(pz;p,w) = \text{diag}(1,1,w,w) \, F_{12}(z;p,w) \, 
\text{diag}(1,1,w^{-1},w^{-1}) \, 
\text{diag}(q^{\half},q^{-\half},q^{-\half},q^{\half}) \, 
R_{12}(pz;p,w)
\label{eq:}
\end{equation}
Therefore, each entry of $F_{12}(z;p,w)$ satisfies a difference equation of 
hypergeometric type. One gets finally
\begin{equation}
F(z;p,w) = \rho_{F}(z;p,w) \left(
\begin{array}{cccc}
1 & 0 & 0 & 0 \\
0 & X_{11}(z) & X_{12}(z) & 0 \\
0 & X_{21}(z) & X_{22}(z) & 0 \\
0 & 0 & 0 & 1 \\
\end{array}
\right)
\end{equation}
where
\begin{equation}
\begin{split}
X_{11}(z) &= \hypergeom{wq^2}{q^2}{w}{p,pq^{-2}z} \\
X_{12}(z) &= \frac{w(q-q^{-1})}{1-w} \; 
\hypergeom{wq^2}{pq^2}{pw}{p,pq^{-2}z} \\
X_{21}(z) &= z \; \frac{pw^{-1}(q-q^{-1})}{1-pw^{-1}} \; 
\hypergeom{pw^{-1}q^2}{pq^2}{p^2w^{-1}}{p,pq^{-2}z} \\
X_{22}(z) &= \hypergeom{pw^{-1}q^2}{q^2}{pw^{-1}}{p,pq^{-2}z}
\end{split}
\end{equation}
and the normalization factor is
\begin{equation}
\rho_{F}(z;p) = \frac{(pz;p,q^4)_{\infty} \; (pq^4z;p,q^4)_{\infty}} 
{(pq^2z;p,q^4)_{\infty}^2}
\end{equation}
The $q$-hypergeometric function $\hypergeom{q^a}{q^b}{q^c}{q,z}$ is defined 
by
\begin{equation}
\hypergeom{q^a}{q^b}{q^c}{q,z} = \sum_{n=0}^\infty 
\frac{(q^a;q)_{n}(q^b;q)_{n}}{(q^c;q)_{n}(q;q)_{n}} \; z^n
\end{equation}
Computing the product $F_{21}(z^{-1};p,w) \, R(z) \, F_{12}(z;p,w)^{-1}$, 
where $R(z)$ is the R-matrix of $\cU_{q}(\widehat{sl}_{2})$ given by 
(\ref{eq:ruq}), one recovers the expression (\ref{eq:Relpb}) of the 
R-matrix of $\elpb{2}$.

\medskip

This example shows that the linear equation (\ref{eq:eqlinB}) indeed plays 
a crucial role. This kind of linear equation was first introduced in 
\cite{BR99} in the context of complex continuation of $6j$ symbols. It was 
then exploited in \cite{ABRR} in the determination and the proof of 
convergence of the dynamical Drinfel'd twists for quantum groups based on 
finite-dimensional Lie algebras.

\begin{remark}
In the face case, the whole construction generalizes without difficulty 
to any Kac--Moody affine algebra $\widehat{\fg}$. $\fH$ being a Cartan 
subalgebra of $\widehat{\fg}$, one defines the Drinfel'd twist 
$\cF(\lambda)$ by formula (\ref{eq:twistBqpl}), where $\cR$ is the 
universal R-matrix of the quantum affine algebra 
$\cU_{q}(\widehat{\fg})$. The proof of the cocycle condition follows 
the same lines as in the $\elpb{N}$ case. Hence one can define a 
quantum affine elliptic algebra ${\cB}_{q,p,\lambda}(\widehat{\fg})$ 
for any Kac--Moody affine algebra $\widehat{\fg}$. \\
One can notice that the situation is rather different in the vertex 
case. Indeed, a cyclic automorphism of the simple root system of 
$\widehat{\fg}$ of order $r = \text{rank}\,\widehat{\fg}$ arises in the 
definition of the vertex-type Drinfel'd twist (\ref{eq:twistAqp}), 
which exists only in the $\widehat{sl}_{N}$ case. 
\finrmq
\end{remark}

\subsection{Degenerations of the quantum affine elliptic algebras}

\subsubsection{Scaling limit procedures}

The so-called scaling limit of an algebra will be understood as the algebra 
defined by the scaling limit of the R-matrix of the initial structure. It 
is obtained by setting in the R-matrix $p=q^{2r}$ (elliptic nome), $w=q^s$ 
(dynamical parameter) and $z=q^u$ (spectral parameter) with $q \to 1$, and 
$r$, $s$, $u$ being kept fixed. The spectral parameter in the Lax operator 
is now to be taken as $u$ (which becomes additive). The R-matrices obtained 
from the scaling limits of the R-matrices of the quantum affine elliptic 
algebras of vertex or face type will be discussed in section 4. It will be 
shown that the corresponding algebraic structures appear to be deformations 
of double Yangians, and that these deformed double Yangians are indeed 
quasi-triangular quasi-Hopf algebras.

\subsubsection{The Sklyanin algebra}

The Sklyanin algebra \cite{Skl82,Skl83} is constructed from $\elpa{2}$ 
taken at $c=0$. The R-matrix (\ref{eq:RAqp}) can be written as
\begin{equation}
R(z) = 1 \otimes 1 + \sum_{\alpha =1}^3 W_\alpha(z) \sigma_\alpha 
\otimes \sigma_\alpha
\end{equation}
where $\sigma_\alpha$ are the Pauli matrices and $W_\alpha(z)$ are 
expressed in terms of the Jacobi elliptic functions. A particular 
$z$-dependence of the $L(z)$ operators is chosen, leading to a 
factorization of the $z$-dependence in the $RLL$ relations. Indeed, setting
\begin{equation}
L(z)= S_0 + \sum_{\alpha=1}^3 W_\alpha(z) S_\alpha \sigma_\alpha
\end{equation}
one obtains an algebra with four generators $S^\alpha$ 
($\alpha=0,\ldots,3$) and commutation relations
\begin{equation}
\begin{split}
& [S_0, S_\alpha ] = -i J_{\beta \gamma } (S_\beta S_\gamma + 
S_\gamma S_\beta ) \\
& [S_\alpha ,S_\beta ] = i (S_0 S_\gamma + S_\gamma S_0)
\end{split}
\end{equation}
where $J_{\alpha \beta} = \displaystyle\frac{W_\alpha^2 - 
W_\beta^2}{W_\gamma^2 - 1}$ and $\alpha$, $\beta$, $\gamma$ are cyclic 
permutations of 1, 2, 3. The structure functions $J_{\alpha \beta}$ are 
actually independent of $z$. Hence we get an algebra where the 
$z$-dependence has been dropped out.

\subsubsection{The dynamical algebra $\cB_{q,\lambda}(\fg)$}

Let $\fg$ be a finite-dimensional complex simple Lie algebra, with 
symmetrized Cartan matrix $(A^{sym}_{ij})$ and inverse ${A^{sym}}^{-1} = 
(d_{ij})$. Let $\fH$ be a Cartan subalgebra of $\fg$ with basis $\{ h_{i} 
\}$ and dual basis $\{ h^\vee_{i} \}$. The positive root system $\Pi^{+}$ 
of $\fg$ is endowed with a normal ordering $\prec$, i.e. if $\alpha, \beta, 
\alpha+\beta \in \Pi^+$ and $[\alpha,\beta]$ is a minimal segment 
containing $\alpha+\beta$, one sets $\alpha \prec \alpha+\beta \prec 
\beta$. The universal R-matrix of the quantum universal enveloping algebra 
$\cU_{q}(\fg)$ is given by
\begin{equation}
\label{eq:univsln}
\cR[\cU_{q}(\fg)] = \Big( \prod_{\gamma \in \Pi^{+}}^{\rightarrow} 
\widehat{\cR}_{\gamma} \Big) \cK
\end{equation}
where the arrow means that the product has to be done with respect to the 
normal ordering $\prec$ defined on $\widehat\Pi^+$ and the factors 
$\widehat{\cR}_{\gamma}$ and $\cK$ are given by
\begin{equation}
\widehat{\cR}_{\gamma} = \exp_{q^{-(\gamma,\gamma)}} \Big( -(q-q^{-1}) 
e_{\gamma} \otimes f_{\gamma} \Big)
\end{equation}
$e_{\gamma}$, $f_{\gamma}$ are the root generators associated to the root 
$\gamma \in \Pi^{+}$ and ${\cal K} = q^{-\sum_{ij} d_{ij} h_{i} \otimes 
h_{j}}$.

\medskip

Let us enounce the following theorem.
\begin{theorem}[\textit{see refs.} \textnormal{\cite{ABRR,JKOS}}]
Let $\phi = q^{\sum_{ij} d_{ij} h_i h_j + 2 \sum_i s_i h_i}$ where 
$s_{i}$ are complex numbers. The linear equation in $\fA \otimes \fA$
\begin{equation}
\label{eq:eqlinBql}
\cF = \Ad(\phi^{-1} \otimes 1) (\cF) \; \cK^{-1} \widehat{\cR} \cK
\end{equation}
has a unique solution in $(\cU_q(\fB^+)\otimes \cU_q(\fB^-))^c$, with 
projection $1 \otimes 1$ on $(\cU_q(\fH)^{\otimes 2})^c$, where the 
superscript $c$ denotes a suitable completion. It is expressed as
\begin{equation}
\label{eq:twistbql}
\cF = \cK^{-1} \widehat{\cF} \cK \qquad \text{where} \qquad 
\widehat{\cF} = \prod_{k \ge 1}^{\leftarrow} \Ad(\phi \otimes 1)^k 
\Big( \widehat{\cR}^{-1} \Big)
\end{equation}
This solution satisfies the shifted cocycle relation
\begin{equation}
\cF_{12}(w) (\Delta \otimes 1) (\cF(w)) = 
\cF_{23}(wq^{{h^\vee}^{(1)}}) (1 \otimes \Delta) (\cF(w))
\end{equation}
with $w=(w_1,\dots,w_{r_\fg}) = q^s = (q^{s_1},\dots,q^{s_{r_\fg}})\in 
\CC^{r_\fg}$, $w q^{h^\vee}=(w_1 q^{h^\vee_1},\dots,w_{r_\fg} 
q^{h^\vee_{r_\fg}})$ and $h^\vee_i = \sum_j d_{ij} h_j$. Hence $\cF(w)$ 
is a GNF twist.
\end{theorem}
The twist $\cF$ leads to the R-matrix $\widetilde\cR = \cF_{21} \cR_{12} 
\cF_{12}^{-1}$, which defines the algebra denoted $\cB_{q,\lambda}(\fg)$.

\medskip

Elements of proof. Expanding the product formula (\ref{eq:univsln}) with 
respect to a Poincar\'{e}--Birkhoff--Witt basis ordered with $\prec$, 
$\widehat\cR$ reads
\begin{equation}
\widehat\cR = \cR \cK^{-1} = 1 \otimes 1 + \sum_{m\in \cZ^*} \sigma_m 
\; \mathbf{e}^m \otimes \mathbf{f}^m
\end{equation}
where $\cZ=\text{Map}(\Pi^+,\NN)$ and $\cZ^*=\cZ\setminus\{(0,\dots,0)\}$. 
The term $\mathbf{e}^m$ (resp. $\mathbf{f}^m$) denotes an element of the 
PBW basis of the deformed enveloping nilpotent subalgebra $\cU_{q}(\cN^+)$ 
(resp. $\cU_{q}(\cN^-)$). Under the assumptions of the theorem,
\begin{equation}
\widehat{\cF} = 1 \otimes 1 + \sum_{\{p,r\}\in (\cZ^*)^2} 
\varphi_{pr}(w) \; \mathbf{e}^p \otimes \mathbf{f}^r
\end{equation}
where the $\varphi_{pr}(w)$ belong to 
$\CC[[s_1,..,s_{r_\fg},s_1^{-1},..,s_{r_\fg}^{-1},\hbar]] \otimes 
(\cU_q(\fH)^{\otimes 2})^c$. They are defined recursively, using 
(\ref{eq:eqlinBql}), by
\begin{equation}
\left( 1 - q^{(-2 {h^\vee}^{(1)} + \gamma_p - s,\gamma_p)} \right) 
\varphi_{pr}(w) = \sum_{\substack{k+m=p \\ l+m=r \\ m \ne 0}} 
(-1)^{[l][m]} a_p^{km} b_r^{lm} \sigma_m q^{(-2 {h^\vee}^{(1)} + 
\gamma_k - s,\gamma_k)} \;\varphi_{kl}(w)
\end{equation}
In the above equation, $\gamma_p$ is the element of the root lattice 
associated to $\mathbf{e}^p$. The scalar product $(\,\cdot\,,\,\cdot\,)$ is 
given by $(x,y)\equiv \sum_{i,j} A^{sym}_{ij} x_i y_j$. The numbers 
$a_p^{km}$ and $b_r^{lm}$ are defined by
\begin{equation}
\mathbf{e}^k \mathbf{e}^m \ =\ \sum_{p\in\cZ} a_{p}^{km} \mathbf{e}^p 
\quad \text{and} \quad \mathbf{f}^l \mathbf{f}^m \ =\ \sum_{r\in\cZ} 
b_{r}^{lm} \mathbf{f}^r
\end{equation}
\finproof

\medskip

Example. In the fundamental representation for $\fg=sl_2$, we get for the 
expression of the twist
\begin{equation}
F(w) = 1 \otimes 1 + \frac{w(q-q^{-1})}{1-w} \; E_{12} \otimes E_{21}
\end{equation}
The R-matrix of $\cB_{q,\lambda}(sl_{2})$ is then given by
\begin{equation}
R[\cB_{q,\lambda}(sl_{2})](w) = q^{-1/2} \left(
\begin{array}{cccc}
1 & 0 & 0 & 0 \\
0 & q & \displaystyle \frac{1-q^2}{1-w} & 0 \\
0 & \displaystyle -\frac{w(1-q^2)}{1-w}\; & \displaystyle 
\;\frac{q(1-wq^2)(1-wq^{-2})}{(1-w)^2} & 0 \\
0 & 0 & 0 & 1 \\
\end{array}
\right)
\end{equation}

\subsubsection{Non-elliptic limits}

Starting from the R-matrix of $\elpb{2}$, and taking the limit $p \to 0$, 
one gets a R-matrix with still a dynamical dependence, which defines the 
algebra $\cU_{q,\lambda}(\widehat{sl}_{2})$. It reads
\begin{equation}
R[\cU_{q,\lambda}(\widehat{sl}_{2})](z,w) = \rho(z) \left(
\begin{array}{cccc}
1 & 0 & 0 & 0 \\
0 & \displaystyle \frac{q(1-z)}{1-q^2z} & \displaystyle 
\frac{(1-q^2)(1-wz)}{(1-q^2z)(1-w)} & 0 \\[12pt] 0 & \displaystyle 
\frac{(1-q^2)(z-w)}{(1-q^2z)(1-w)} & \displaystyle 
\frac{q(1-z)}{(1-q^2z)} \; \frac{(1-wq^2)(1-wq^{-2})}{(1-w)^2} & 0 
\\
0 & 0 & 0 & 1 \\
\end{array} \right)
\end{equation}
The normalization factor is
\begin{equation}
\rho(z) = q^{-1/2} \; \frac{(q^2z;q^4)_{\infty}^2}{(z;q^4)_{\infty} \; 
(q^4z;q^4)_{\infty}}
\end{equation}
This algebra can be seen as a dynamical version of 
$\cU_{q}(\widehat{sl}_{2})$. This can be easily generalized to any 
Kac--Moody (untwisted) affine algebra $\widehat\fg$. The Hopf structure of 
these algebras $\cU_{q}(\widehat\fg)$ is determined by the following 
proposition.
\begin{proposition}
The twist (\ref{eq:twistbql}) applied to the universal R-matrix of 
$\cU_{q}(\widehat\fg)$ leads to the R-matrix of a QTQHA denoted 
$\cU_{q,\lambda}(\widehat\fg)$.
\end{proposition}
Proof. Using the fact that $\cU_{q}(\fg)$ is a Hopf subalgebra of 
$\cU_{q}(\widehat\fg)$, the twist (\ref{eq:twistbql}) can be used to 
construct the dynamical algebra $\cU_{q,\lambda}(\widehat\fg)$. Indeed, the 
twist $\cF$ (\ref{eq:twistbql}) seen as an element of 
$\cU_{q}(\widehat\fg)^{\otimes 2}$ satisfies the shifted cocycle condition, 
yielding a dynamical R-matrix
\begin{equation}
\cR_{U_{q,\lambda}(\widehat\fg)}(w) = \cF_{21}(w) \; 
\cR_{U_{q}(\widehat\fg)} \; \cF_{12}^{-1}(w)
\end{equation}
\finproof

\begin{proposition}
The evaluation representation of $\cR_{U_{q,\lambda}(\widehat\fg)}(w)$ 
for $\fg = sl_N$, is identified with the $p \to 0$ limit of the 
evaluation representation of the elliptic $\elpb{N}$ R-matrix 
(\ref{eq:RBqpln}).
\end{proposition}
Proof. Direct computation of the matrix elements 
$R_{i_{1}i_{2}}^{j_{1}j_{2}}$ of the R-matrix of 
$\cU_{q,\lambda}(\widehat{sl}_{N})$ ($1 \le a,b \le N$) gives:
\begin{align}
R &= \rho(z) \Big( \sum_{a} E_{aa} \otimes E_{aa} + \sum_{a,b} 
\frac{(1-q^2)(1-w_{ab}z)} {(1-q^2z)(1-w_{ab})} E_{ab} \otimes E_{ba} 
\nonumber \\
&+ \sum_{a<b} \frac{q(1-z)} {1-q^2z} E_{aa} \otimes E_{bb} + \sum_{a>b} 
\frac{q(1-z)} {(1-q^2z)} \; \frac{(1-w_{ab}q^2)(1-w_{ab}q^{-2})} 
{(1-w_{ab})^2} E_{aa} \otimes E_{bb} \Big)
\label{eq:ruqln}
\end{align}
the normalization factor being given by
\begin{equation}
\rho(z) = q^{-\frac{N-1}{N}} \; \frac{(q^{2}z;q^{2N})_{\infty} \; 
(q^{2N-2}z;q^{2N})_{\infty}}{(z;q^{2N})_{\infty} \; 
(q^{2N}z;q^{2N})_{\infty}}
\end{equation}
We recognize the limit $p \to 0$ of the R-matrix (\ref{eq:RBqpln}). The 
R-matrix (\ref{eq:ruqln}) satisfies the dynamical Yang--Baxter equation 
(\ref{eq:DYBE}). \finproof

\section{Double Yangians and related structures}
\setcounter{equation}{0}

\subsection{Yangians}

\subsubsection{Definition of the Yangians}

Let $\fg$ be a finite-dimensional complex simple Lie algebra and consider 
$\fg[u] = \fg \otimes \CC[u]$, where $\CC[u]$ is the ring of polynomials in 
the indeterminate $u$ (by misuse of language, we will call $\fg[u]$ the 
half-loop algebra of $\fg$). $\fg[u]$ is endowed with its standard 
bialgebra structure $\delta : \fg[u] \to \fg[u] \otimes \fg[u]$ defined by 
(note that $\fg[u] \otimes \fg[u]$ is isomorphic to $(\fg \otimes 
\fg)[u,v]$ where $v$ is a second indeterminate)
\begin{equation}
\delta(f)(u,v) = \left[ f(u) \otimes 1 + 1 \otimes f(v) \,,\, 
\frac{C}{u-v} \right]
\end{equation}
where $C$ is the second order tensorial Casimir element of $\fg$ associated 
to a given invariant bilinear form on $\fg$ (for example the Killing form). 
\\
Let $\CC[[\hbar]]$ be the ring of formal power series in the inderterminate 
$\hbar$. Then there exists a unique quantization $\cU_{\hbar}(\fg)$ of 
$(\fg[u],\delta)$ which is a graded Hopf algebra over $\CC[[\hbar]]$, the 
gradation being defined by setting $\deg\hbar = 1$ \cite{Dri86}, i.e.
\begin{align}
& U_{\hbar}(\fg[u])/\hbar U_{\hbar}(\fg[u]) \simeq U(\fg[u]) \; 
\text{as graded algebras over $\CC$}
\intertext{and}
& \frac{1}{\hbar} \big( \Delta - \Delta^{\op} \big)(x) 
\big\vert_{\!\!\mod\hbar} = \delta \big( x \big\vert_{\!\!\mod\hbar} 
\big) \quad \text{for} \; x \in U_{\hbar}(\fg[u])
\end{align}

$\cU_{\hbar=1}(\fg)$ is a Hopf algebra over $\CC$, which is called the 
Yangian of $\fg$ and is denoted by $Y(\fg)$. It has been introduced by 
Drinfel'd in ref. \cite{Dri85}.

\medskip

There exists for the Yangian $Y(\fg)$ three different realizations, due to 
Drinfel'd \cite{Dri85,Dri86,Dri88}. In the first realization the Yangian is 
generated by the elements $J^a_0$ of the Lie algebra and a set of other 
generators $J^a_1$ in one-to-one correspondance with $J^a_0$ realizing a 
representation space thereof. The second realization is given in terms of 
generators and relations similar to the description of the loop algebra as 
a space of maps. However in this realization no explicit formula for the 
comultiplication is known in general. The third realization is obtained in 
the FRT formalism.

\subsubsection{FRT formalism for the Yangians}

The Yangian $Y(\fg)$ can be constructed in the FRT formalism as follows. \\
Let $\cU(R)$ be the algebra generated by the operators $T^{ij}_{(n)}$, for $1 
\le i,j \le N$, $n \in \NN$, encapsulated into a $N \!\times\! N$ matrix 
($E_{ij} \in \text{End}(\CC^{N})$ are the standard elementary matrices)
\begin{equation}
T(u) = \sum_{n \in \NN} T_{(n)} \, u^{-n} = \sum_{i,j=1}^{N} \sum_{n 
\in \NN} T^{ij}_{(n)} \, u^{-n} \, E_{ij} = \sum_{i,j=1}^{N} T^{ij}(u) 
\, E_{ij}
\end{equation}
and $T^{ij}_{(0)} = \delta_{ij}$, imposing the following constraints on 
$T(u)$
\begin{equation}
R_{12}(u-v) \, T_1(u) \, T_2(v) = T_2(v) \, T_1(u) \, R_{12}(u-v)
\end{equation}
where $R_{12}$ is a $N \!\times\! N$ matrix which is a rational solution of 
the Yang--Baxter equation. \\
The Hopf algebra structure of $\cU(R)$ is given by \cite{FRT}
\begin{equation}
\label{eq:GopfYang}
\Delta \big( T^{ij}(u) \big) = \sum_{k=1}^{N} T^{ik}(u) \otimes 
T^{kj}(u) \quad ; \quad S(T(u)) = T(u)^{-1} \quad ; \quad 
\epsilon(T(u)) = \II_{N}
\end{equation}

\medskip

In the case of $\fg = sl(N)$, the matrix $R_{12}(u)$ is given by
\begin{equation}
R(u) = \II_{N} + \frac{P}{u}
\end{equation}
where $\II_{N}$ is the $N \!\times\! N$ unit matrix and $P = 
\sum_{i,j=1}^{N} E_{ij} \otimes E_{ji}$ is the permutation matrix. The 
explicit commutation relations between the generators $T^{ij}_{(n)}$ read 
($m,n \ge 0$)
\begin{equation}
\big[ T^{ij}_{(m+1)} \,,\, T^{kl}_{(n)} \big] - \big[ T^{ij}_{(m)} 
\,,\, T^{kl}_{(n+1)} \big] = T^{kj}_{(n)} \, T^{il}_{(m)} - 
T^{kj}_{(m)} \, T^{il}_{(n)}
\end{equation}
The quantum determinant of $T(u)$ defined by
\begin{equation}
\qdet T(u) = \sum_{\sigma\in{\mathfrak S}_N} \varepsilon(\sigma) 
\prod_{i=1}^N T^{i,\sigma(i)}(u-N+i) = \sum_{n \in \NN} c_{n} u^{-n}
\end{equation}
where $\varepsilon(\sigma)$ is the signature of the permutation $\sigma$, 
lies in the center of $\cU(R)$ \cite{MNO96}. Moreover, the center of 
$\cU(R)$ is generated by the coefficients $c_{n}$ and one has $\Delta(\qdet 
T(u)) = \qdet T(u) \otimes \qdet T(u)$. The Yangian $Y(sl(N))$ is then 
identified as the quotient algebra $\cU(R)/\big\langle \qdet T(u) - 1 
\big\rangle$.

\medskip

In the case of $\fg = so(N)$ or $sp(N)$, the matrix $R_{12}(u)$ take the 
form \cite{AACFR,Dri85,KuSkl82}
\begin{equation}
R(u) = \II_{N} + \frac{P}{u} - \frac{K}{u+\kappa}
\end{equation}
where $\II_{N}$ and $P$ are defined as above and $K = \sum_{i,j=1}^{N} 
\epsilon_{i} \epsilon_{j} E_{\bar{\jmath}\bar{\imath}} \otimes E_{ji}$ with 
$\bar{\imath} = N+1-i$. For $so(N)$, $\epsilon_{i} = 1$ for all $i$, while 
for $sp(N)$ with $N=2n$, $\epsilon_{i} = 1$ if $1 \le i \le n$ and 
$\epsilon_{i} = -1$ if $n+1 \le i \le N$. The commutation relations between 
the generators $T^{ij}_{(n)}$ become now
\begin{align}
\big[ T^{ij}_{(m+2)} , T^{kl}_{(n)} \big] &- \big[ T^{ij}_{(m)} , 
T^{kl}_{(n+2)} \big] = 2 \big[ T^{ij}_{(m+1)} , T^{kl}_{(n+1)} \big] - 
\kappa \big[ T^{ij}_{(m+1)} , T^{kl}_{(n)} \big] + \kappa \big[ 
T^{ij}_{(m)} , T^{kl}_{(n+1)} \big] \nonumber \\
&+ T^{kj}_{(n)} T^{il}_{(m+1)} - T^{kj}_{(m+1)} T^{il}_{(n)} - 
T^{kj}_{(n+1)} T^{il}_{(m)} + T^{kj}_{(m)} T^{il}_{(n+1)} + \kappa 
T^{kj}_{(n)} T^{il}_{(m)} - \kappa T^{kj}_{(m)} T^{il}_{(n)} \nonumber 
\\
&+ \sum_{r} \Big( \delta_{i\bar{k}} \, \epsilon_{\bar{\imath}} 
\epsilon_{\bar{r}} \, \big( T^{rj}_{(m+1)} T^{\bar{r}l}_{(n)} - 
T^{rj}_{(m)} T^{\bar{r}l}_{(n+1)} \big) - \delta_{j\bar{l}} \, 
\epsilon_{\bar{r}} \epsilon_{\bar{\jmath}} \, \big( T^{k\bar{r}}_{(n)} 
T^{ir}_{(m+1)} - T^{ir}_{(m+1)} T^{k\bar{r}}_{(m)} \big) \Big)
\end{align}
where $m,n \ge -2$ and by convention $T^{ij}_{(n)} = 0$ for $n<0$. \\
The operators generated by $C^{ij}(u) = \sum_{k} \epsilon_{i} \epsilon_{k} 
T^{\bar{k}\bar{\imath}}(u-\kappa) \, T^{kj}(u)$ are such that $C^{ij}(u) = 
\delta_{ij} \, c(u)$. The element $c(u)$ lies in the center of $\cU(R)$ and 
satisfies $\Delta(c(u)) = c(u) \otimes c(u)$. It generates a Hopf ideal. 
The Yangian $Y(\fg)$ is then given by the quotient algebra 
$\cU(R)/\big\langle c(u) - 1 \big\rangle$.

\subsubsection{Drinfel'd second realization of the Yangians}

The Yangian $Y(\fg)$ is isomorphic to the associative algebra over $\CC$ 
with generators $e_{i,n}$, $f_{i,n}$ and $h_{i,n}$ where $i = 1,\ldots,r$ 
($r$ is the rank of $\fg$) and $n \in \NN$, and defining relations
\begin{equation}
\label{eq:reldriyang}
\begin{split}
& [h_{i,m} , h_{j,n}] = 0 \hspace*{72pt} [e_{i,m} , f_{j,n}] = 
\delta_{ij} h_{i,m+n} \\
& [h_{i,0} , e_{j,n}] = 2\alpha_{ij} e_{j,n} \hspace*{45pt} 
[h_{i,0} , f_{j,n}] = -2\alpha_{ij} f_{j,n} \\
& [h_{i,m+1} , e_{j,n}] - [h_{i,m} , e_{j,n+1}] = \alpha_{ij} \{ 
h_{i,m} , e_{j,n} \} \\
& [h_{i,m+1} , f_{j,n}] - [h_{i,m} , f_{j,n+1}] = -\alpha_{ij} \{ 
h_{i,m} , f_{j,n} \} \\
& [e_{i,m+1} , e_{j,n}] - [e_{i,m} , e_{j,n+1}] = \alpha_{ij} \{ 
e_{i,m} , e_{j,n} \} \\
& [f_{i,m+1} , f_{j,n}] - [f_{i,m} , f_{j,n+1}] = -\alpha_{ij} \{ 
f_{i,m} , f_{j,n} \}
\end{split}
\end{equation}
and for $i \ne j$ with $n_{ij} = 1-A_{ij}$
\begin{equation}
\begin{split}
& \sum_{\sigma \in {\mathfrak S}_{n_{ij}}} [e_{i,m_{\sigma(1)}} , 
[e_{i,m_{\sigma(2)}} , \ldots , [ e_{i,m_{\sigma(n_{ij})}} , 
e_{j,n} ]]] = 0 \\
& \sum_{\sigma \in {\mathfrak S}_{n_{ij}}} [f_{i,m_{\sigma(1)}} , 
[f_{i,m_{\sigma(2)}} , \ldots , [ f_{i,m_{\sigma(n_{ij})}} , 
f_{j,n} ]]] = 0
\end{split}
\end{equation}
where $\alpha_{ij} \equiv \half A^{sym}_{ij}$ and $(A^{sym}_{ij})$ is the 
symmetrized Cartan matrix of $\fg$ properly normalized, as explained in 
section \ref{sect:212}.

Unfortunately the Hopf structure for this presentation of the Yangian 
$Y(\fg)$ is not explicit and no formula for the comultiplication of the 
generators $e_{i,n}$, $f_{i,n}$ and $h_{i,n}$ is known in the general case. 
In the case of $sl_{2}$, in terms of the generating functions (the index 
$i$ is omitted),
\begin{equation}
e^{+}(u) = \sum_{k \ge 0} e_{k} u^{-k-1} \quad f^{+}(u) = \sum_{k \ge 
0} f_{k} u^{-k-1} \quad h^{+}(u) = 1 + \sum_{k \ge 0} h_{k} u^{-k-1}
\end{equation}
the explicit formulae for the Hopf structure take the form \cite{Molev2001}
\begin{equation}
\label{eq:GopfDri1}
\begin{split}
& \Delta(e^{+}(u)) = e^{+}(u) \otimes 1 + \sum_{k=0}^{\infty} 
(-1)^{k} (f^{+}(u+1))^{k} h^{+}(u) \otimes (e^{+}(u))^{k+1} \\
& \Delta(f^{+}(u)) = 1 \otimes f^{+}(u) + \sum_{k=0}^{\infty} 
(-1)^{k} (f^{+}(u))^{k+1} \otimes h^{+}(u) (e^{+}(u+1))^{k} \\
& \Delta(h^{+}(u)) = \sum_{k=0}^{\infty} (-1)^{k} (k+1) 
(f^{+}(u+1))^{k} h^{+}(u) \otimes h^{+}(u) (e^{+}(u+1))^{k}
\end{split}
\end{equation}
for the coproduct,
\begin{equation}
\label{eq:GopfDri2}
\begin{split}
& S(e^{+}(u)) = - \big( h^{+}(u) + f^{+}(u+1) e^{+}(u) \big)^{-1} 
\, e^{+}(u) \\
& S(f^{+}(u)) = - f^{+}(u) \, \big( h^{+}(u) + f^{+}(u) e^{+}(u+1) 
\big)^{-1} \\
& S(h^{+}(u)) = \big( h^{+}(u) + f^{+}(u+1) e^{+}(u) \big)^{-1} \, 
\Big( 1 - f^{+}(u+1) \big( h^{+}(u) + f^{+}(u+1) e^{+}(u) 
\big)^{-1} e^{+}(u) \Big)
\end{split}
\end{equation}
for the antipode, and
\begin{equation}
\label{eq:GopfDri3}
\epsilon(e^{+}(u)) = 0 \qqquad \epsilon(f^{+}(u)) = 0 \qqquad 
\epsilon(h^{+}(u)) = 1
\end{equation}
for the counit.

\medskip

The link between the FRT formalism of the Yangian and the Drinfel'd second 
realization of the Yangian is given in the $sl_{2}$ case by the following 
Hopf isomorphism:
\begin{equation}
\label{eq:isoFRT}
\begin{split}
& e^{+}(u) \mapsto {T^{22}(u)}^{-1} \; T^{12}(u) \\
& f^{+}(u) \mapsto T^{21}(u) \; {T^{22}(u)}^{-1} \\
& h^{+}(u) \mapsto T^{11}(u) \; {T^{22}(u)}^{-1} - T^{21}(u) \; 
{T^{22}(u)}^{-1} \; T^{12}(u) \; {T^{22}(u)}^{-1} \\
\end{split}
\end{equation}
The corresponding result for the $sl_{N}$ case is rather cumbersome and is 
explicited in \cite{Io96}. Note that the Hopf structure (\ref{eq:GopfYang}) 
and the isomorphism (\ref{eq:isoFRT}) -- or its generalization -- allows 
one in principle to derive equations 
(\ref{eq:GopfDri1})--(\ref{eq:GopfDri3}). In practice this seems tractable 
in the $sl_{2}$ case only.

\medskip

\begin{remark}
Following a remark of Drinfel'd \cite{Dri86}, the commutation relations 
(\ref{eq:reldriyang}) of the Yangian can be obtained by the following 
construction, relying on the fact that $Y(\fg) \simeq A/\hbar A$ where 
$A$ is a suitable subsalgebra of $\cU_{q}(\widehat\fg)_{c=0} 
\underset{\CC[[\hbar]]}{\otimes} \CC((\hbar))$ and $q = e^\hbar$. 
Consider the Drinfel'd basis of $\cU_{q}(\widehat\fg)$ (see section 
\ref{sect:drinfeld}) with generators $\cK_{i,n}$, $\cX^{\pm}_{i,n}$ 
(with $1 \le i \le r$ and $n \in \ZZ$) satisfying the commutation 
relations (\ref{eq:reldri}) and set $C \equiv e^{\hbar c} = 1$. One 
defines for $n$ non-negative integer [the following calculations are 
due to A.I. Molev]
\begin{equation}
\begin{split}
& {\,\widetilde{\cX}^{\pm}_{i,n}}^{(\ell)} = 
\frac{(-1)^n}{\hbar^n} \sum_{k=0}^n (-1)^k \binom{n}{k} \; 
\cX^{\pm}_{i,k+\ell} \\
& {\,\widetilde{\cH}_{i,n}}^{(\ell)} = \frac{(-1)^n}{\hbar^n} 
\sum_{k=0}^n (-1)^k \binom{n}{k} \; \frac{\Psi^{+}_{i,k+\ell} - 
\Psi^{-}_{i,k+\ell}}{q-q^{-1}}
\end{split}
\end{equation}
which satisfy
\begin{equation}
{\,\widetilde{\cX}^{\pm}_{i,n}}^{(\ell)} = \hbar 
{\,\widetilde{\cX}^{\pm}_{i,n+1}}^{(\ell-1)} + 
{\,\widetilde{\cX}^{\pm}_{i,n}}^{(\ell-1)} \qquad \text{and} \qquad 
{\,\widetilde{\cH}_{i,n}}^{(\ell)} = \hbar 
{\,\widetilde{\cH}_{i,n+1}}^{(\ell-1)} + 
{\,\widetilde{\cH}_{i,n}}^{(\ell-1)}
\end{equation}
Then the commutation relations between the generators 
${\,\widetilde{\cX}^{\pm}_{i,n}}^{(0)}$ and 
${\,\widetilde{\cH}_{i,n}}^{(0)}$ ($n \in \NN$) in the quotient 
$A/\hbar A$ are equivalent to the commutation relations of the Yangian 
$Y(\fg)$ with the identification ${\,\widetilde{\cX}^{+}_{i,n}}^{(0)} 
\mapsto e_{i,n}$, ${\,\widetilde{\cX}^{-}_{i,n}}^{(0)} \mapsto f_{i,n}$ 
and ${\,\widetilde{\cH}_{i,n}}^{(0)} \mapsto h_{i,n}$. \finrmq
\end{remark}

\subsection{Double Yangians}

\subsubsection{Notion of quantum double}

Let $\fA$ and $\fA'$ be two finite-dimensional Hopf algebras and $\cR$ an 
invertible element of $\fA \otimes \fA'$ such that
\begin{equation}
\label{eq:hopfinv}
\begin{split}
& (\Delta^{\fA} \otimes \id) (\cR) = \cR_{13} \, \cR_{23} 
\qqquad (S^{\fA} \otimes \id) (\cR) = \cR^{-1} \\
& (\id \; \otimes \Delta^{\fA'}) (\cR) = \cR_{12} \, \cR_{13} 
\qqquad (\id \; \otimes S^{\fA'}) (\cR) = \cR^{-1}
\end{split}
\end{equation}
Then $\fA' \otimes \fA$ is also a Hopf algebra with coalgebra structure 
given by ($a' \otimes a \in \fA' \otimes \fA$)
\begin{align}
& \text{coproduct:} && \Delta(a' \otimes a) = \cR_{23} \; 
\Delta_{13}^{\fA'}(a') \; \Delta_{24}^{\fA}(a) \cR_{23}^{-1} &&&& 
\\[6pt] 
& \text{antipode:} && S(a' \otimes a) = \cR_{21}^{-1} \; \big( 
S^{\fA'}(a') \otimes S^{\fA}(a) \big) \; \cR_{21} &&&& \\[6pt] 
& \text{counit:} && \epsilon(a' \otimes a) = \epsilon^{\fA'}(a') \; 
\epsilon^{\fA}(a) &&&&
\end{align}
This Hopf algebra is denoted $\fA' \underset{\cR}{\otimes} \fA$.

\medskip

Let $\fA$ be a finite-dimensional Hopf algebra, $\fA_{\op}$ be the Hopf 
algebra with opposite multiplication $m_{\op} = m \circ \sigma$, 
$\fA^{\op}$ be the Hopf algebra with opposite comultiplication 
$\Delta^{\op} = \sigma \circ \Delta$ and $\fA^{*}$ be the dual Hopf algebra 
of $\fA$. One denotes by $\{a_{i}\}$ a basis of $\fA$ and by 
$\{a^{*}_{i}\}$ the dual basis of $\fA^{*}$. Consider the element $\cR \in 
\fA_{\op} \otimes \fA^{*}$ given by $\cR = \sum_{i} a_{i} \otimes 
a^{*}_{i}$. Then $\cR$ satisfies (\ref{eq:hopfinv}). It follows that 
$\fA^{*} \underset{\cR}{\otimes} \fA_{\op}$ is a Hopf algebra as described 
above.

\begin{proposition}
The algebra $\cD\fA = (\fA^{*} \underset{\cR}{\otimes} \fA_{\op})^{*}$ 
is a quasi-triangular Hopf algebra (isomorphic to $\fA \otimes 
(\fA^{*})^{\op}$ as a coalgebra), called the quantum double of $\fA$. 
The universal $\cR$ matrix of $\cD\fA$ is given by the canonical 
element of $\fA \otimes (\fA^{*})^{\op} \subset \cD\fA\otimes \cD\fA$ 
associated to the identity map $\fA \to \fA$.
\end{proposition}

\begin{remark}
In the case where $\fA$ is an infinite-dimensional Hopf algebra (which 
is the case we are hereafter interested in), one has to face with some 
difficulties \cite{ChaPre}. Indeed, when $\fA$ is finite-dimensional, 
the multiplication $\mu$ of $\fA$ induces the comultiplication on 
$\fA^{*}$ by $\mu^{*} : \fA^{*} \to \fA^{*} \otimes \fA^{*} \subset 
(\fA \otimes \fA)^{*}$. When $\fA$ is infinite-dimensional, it might 
happen that $\mu^{*}(\fA^{*}) \not\subset \fA^{*} \otimes \fA^{*}$, 
which prevents one from endowing the dual algebra $\fA^{*}$ with a 
canonical dual Hopf structure. One way to escape this problem is to 
introduce the notion of (non-degenerate) pairing of two Hopf algebras. 
A Hopf pairing $\langle \;,\; \rangle : \fA \otimes \fA' \to \CC$ 
between two Hopf algebras $\fA$ and $\fA'$ is a bilinear map such that 
for all $a_{1},a_{2} \in \fA$ and $a_{1}',a_{2}' \in \fA'$
\begin{equation}
\label{eq:pairing}
\begin{split}
& \langle \mu^{\fA}(a_{1} \otimes a_{2}),a' \rangle = \langle 
a_{1} \otimes a_{2},\Delta^{\fA'}(a') \rangle \qqquad \langle 
a,\mu^{\fA'}(a_{1}' \otimes a_{2}') \rangle = \langle 
\Delta^{\fA}(a),a_{1}' \otimes a_{2}' \rangle \\
& \langle S^{\fA}(a),a' \rangle = \langle a,S^{\fA'}(a') 
\rangle \qqquad \langle \epsilon^{\fA}(a),a' \rangle = \langle 
a,\iota^{\fA'}(a') \rangle \qqquad \langle \iota^{\fA}(a),a' 
\rangle = \langle a,\epsilon^{\fA'}(a') \rangle
\end{split}
\end{equation}
and $\langle a_{1} \otimes a_{2},a_{1}' \otimes a_{2}' \rangle = 
\langle a_{1},a_{1}' \rangle \langle a_{2},a_{2}' \rangle$. \\
When $\fA$ is finite-dimensional, eq. (\ref{eq:pairing}) is equivalent 
to have $\fA' = \fA^{*}$. When $\fA$ is infinite-dimensional, this is 
however a weaker statement. \\
In the case of a quantum double, it may be convenient to slightly 
modify these equations to work with the opposite comultiplication on 
the second Hopf algebra. \finrmq
\end{remark}

\subsubsection{Definition of the double Yangians}

Although the Yangian $Y(\fg)$ is a Hopf algebra, it is not 
quasi-triangular. In order to get a quasi-triangular Hopf algebra, one has 
to construct the quantum double $\dy{\fg}$ of the Yangian $Y(\fg)$. The 
double Yangian $\dy{\fg}$ is generated by the generators $e_{i,n}$, 
$f_{i,n}$ and $h_{i,n}$ where $i = 1,\ldots,r$ ($r$ is the rank of $\fg$) 
and $n \in \ZZ$, satisfying the relations (\ref{eq:reldriyang}). It is 
convenient to write the commutation relations of the double Yangian by 
introducing the following generating functionals:
\begin{equation}
\label{eq:genfun}
e_{i}^{\pm}(u) = \pm\sum_{\substack{k \ge 0 \\ k < 0}} e_{i,k} u^{-k-1} 
\;, \qquad 
f_{i}^{\pm}(u) = \pm\sum_{\substack{k \ge 0 \\ k < 0}} f_{i,k} u^{-k-1} 
\;, \qquad 
h_{i}^{\pm}(u) = 1 \pm\sum_{\substack{k \ge 0 \\ k < 0}} h_{i,k} u^{-k-1}
\end{equation}
and $e_{i}(u) = e_{i}^{+}(u) - e_{i}^{-}(u)$, $f_{i}(u) = f_{i}^{+}(u) - 
f_{i}^{-}(u)$. One gets
\begin{equation}
\begin{split}
& e_{i}(u) \, e_{j}(v) = \frac{u-v+\alpha_{ij}}{u-v-\alpha_{ij}} \; 
e_{i}(v) \, e_{j}(u) \\
& f_{i}(u) \, f_{j}(v) = \frac{u-v-\alpha_{ij}}{u-v+\alpha_{ij}} \; 
f_{i}(v) \, f_{j}(u) \\
& h_{i}^{\pm}(u) \, e_{j}(v) = 
\frac{u-v+\alpha_{ij}}{u-v-\alpha_{ij}} \; e_{j}(v) \, 
h_{i}^{\pm}(u) \\
& h_{i}^{\pm}(u) \, f_{j}(v) = 
\frac{u-v-\alpha_{ij}}{u-v+\alpha_{ij}} \; f_{j}(v) \, 
h_{i}^{\pm}(u) \\
& \big[ e_{i}(u) , f_{j}(v) \big] = \delta_{ij} \, \big( 
\delta(u-v) h_{i}^{+}(u) - \delta(u-v) h_{i}^{-}(v) \big)
\end{split}
\end{equation}
and for $i \ne j$ with $n_{ij} = 1-A_{ij}$
\begin{equation}
\label{eq:serreyang}
\begin{split}
& \sum_{\sigma \in {\mathfrak S}_{n_{ij}}} [e_{i}(u_{\sigma(1)}) , 
[e_{i}(u_{\sigma(2)}) , \ldots , [ e_{i}(u_{\sigma(n_{ij})}) , 
e_{j}(v) ]]] = 0 \\
& \sum_{\sigma \in {\mathfrak S}_{n_{ij}}} [f_{i}(u_{\sigma(1)}) , 
[f_{i}(u_{\sigma(2)}) , \ldots , [ f_{i}(u_{\sigma(n_{ij})}) , 
f_{j}(v) ]]] = 0
\end{split}
\end{equation}
where $\alpha_{ij} \equiv \half A^{sym}_{ij}$ and $\delta(u-v) \equiv 
\sum_{n \in \ZZ} u^{-n-1} v^n$.

\medskip

The double Yangian $\dy{\fg}$ is described as follows \cite{Kho95,KT96}. 
Let $Y^{\pm}$ be the algebras generated by the generating functionals 
$e_{i}^{\pm}(u)$, $f_{i}^{\pm}(u)$ and $h_{i}^{\pm}(u)$. One has $Y^{+} = 
Y(\fg) \subset \dy{\fg}$ and the dual with opposite comultiplication 
$Y(\fg)^{\op} \simeq Y^{-}$ (more precisely a suitable formal completion of 
it). The Hopf pairing between the generators of
$Y^{+}$ and those of $Y^{-}$ is given by
\begin{equation}
\langle e_{i}^{+}(u),f_{j}^{-}(v) \rangle = \frac{\delta_{ij}}{u-v} 
\qqquad \langle f_{i}^{+}(u),e_{j}^{-}(v) \rangle = 
\frac{\delta_{ij}}{u-v} \qqquad \langle h_{i}^{+}(u),h_{j}^{-}(v) 
\rangle = \frac{u-v+\alpha_{ij}}{u-v-\alpha_{ij}}
\end{equation}

\subsubsection{Universal R-matrix of the double Yangian $\dy{\fg}$}

\textsl{Universal R-matrix for $\dy{sl_{2}}$}

Being a quantum double, $\dy{\fg}$ is naturally endowed with a structure of 
quasi-triangular Hopf algebra and admits a universal R-matrix. In the case 
of $\dy{sl_{2}}$, the following decomposition of the universal R-matrix has 
been proved \cite{KT96}
\begin{equation}
\cR[\dy{sl_{2}}] = \cR_{E} \; \cR_{H} \; \cR_{F}
\end{equation}
where
\begin{align}
& \cR_{E} = \prod_{k \ge 0}^{\rightarrow} \exp(-e_{k} \otimes f_{-k-1}) 
\;, \qquad \qquad \cR_{F} = \prod_{k \ge 0}^{\leftarrow} \exp(-f_{k} 
\otimes e_{-k-1}) \;, 
\label{eq:RERF} \\
& \cR_{H} = \prod_{k \ge 0} \exp \left( \sum_{n \ge 0} \Big( 
-\frac{d}{du} \ln h_{+}(u) \Big)_{n} \otimes \Big( \ln h_{-}(x+2k+1) 
\Big)_{-n-1} \right)
\label{eq:RH}
\end{align}

\textsl{Universal R-matrix for $\dy{sl_{N}}$}

In the case of $\dy{sl_{N}}$, the same kind of decomposition of the 
universal R-matrix takes place
\begin{equation}
\cR[\dy{sl_{N}}] = \cR_{E} \; \cR_{H} \; \cR_{F}
\end{equation}
Unfortunately, although the factor $\cR_{H}$ has been rigorously derived, 
the factors $\cR_{E}$ and $\cR_{F}$ remain only conjectured \cite{KT96}. 
The factor $R_{H}$ is given by
\begin{equation}
\cR_{H} = \prod_{k \ge 0} \exp \Bigg( \sum_{i,j=1,\ldots,r} \sum_{n \ge 
0} \Big( -\frac{d}{du} \ln h_{i,+}(u) \Big)_{n} \otimes \Big( 
C_{ji}(T^{\half}) \ln h_{j,-} \big( x+N(k+\sfrac{1}{2}) \big) 
\Big)_{-n-1} \Bigg)
\end{equation}
Let $A^{sym}(q)$ be the quantum analogue of the the symmetrized Cartan 
matrix of $sl_{N}$, i.e. $A^{sym}_{ij}(q) = [A^{sym}_{ij}]_{q}$. The matrix 
$C(q)$ is defined by $C(q) = [N]_{q} \, (A^{sym}(q))^{-1}$. The operator $T 
= \exp(d/dx)$ is the shift operator, $T f(x) = f(x+1)$.

In order to write the expression of the factors $\cR_{E}$ and $\cR_{F}$, 
one needs to fix some notations. Let $E$ and $F$ be the unital subalgebras 
of $\dy{sl_{N}}$ generated respectively by the generators $e_{i,k}$ ($k \ge 
0$) and $f_{i,k}$ ($k \ge 0$). Let $\Pi$ and $\Pi^{+} \subset \Pi$ be the 
root and positive root systems of $sl_{N}$ and $\widehat{\Pi} = \{ \gamma + 
n\delta \,, m\delta \,|\, \gamma \in \Pi \,, n,m \in \ZZ \,, m \ne 0\}$ be 
the root system of $\widehat{sl}_{N}$. One defines the following two 
subsets of $\widehat{\Pi}$: $\Pi_{E} = \{\gamma+n\delta \,|\, 
\gamma\in\Pi^{+}, n\in\NN\}$ and $\Pi_{F} = \{-\gamma+n\delta \,|\, 
\gamma\in\Pi^{+}, n\in\NN\}$. $\Pi_{E}$ and $\Pi_{F}$ are equipped with two 
orderings $<_{E}$ and $<_{F}$ such that
\begin{align*}
& \xi_{1} <_{E} \xi_{1} + \xi_{2} <_{E} \xi_{2} \;\; \text{if} \;\; 
\xi_{1}, \xi_{2}, \xi_{1} + \xi_{2} \in \Pi_{E} \quad \text{and} \quad 
\gamma+n\delta <_{E} \gamma+m\delta \\
& \xi_{2} <_{F} \xi_{1} + \xi_{2} <_{F} \xi_{1} \;\; \text{if} \;\; 
\xi_{1}, \xi_{2}, \xi_{1} + \xi_{2} \in \Pi_{F}\quad \text{and} \quad 
-\gamma+n\delta >_{F} -\gamma+m\delta
\end{align*}
where $\gamma \in \Pi^{+}$ and $n < m$.

By induction, for any root $\xi = \gamma+n\delta \in \Pi_{E}$ (resp. $\xi = 
-\gamma+n\delta \in \Pi_{F}$), one constructs the roots generators $e_{\xi} 
\equiv e_{\gamma,n}$ (resp. $f_{\gamma,n}$) and $e_{-\xi} \equiv 
f_{\gamma,-n}$ (resp. $e_{\gamma,-n}$) by
\begin{equation}
e_{\xi} = [e_{\xi_{1}},e_{\xi_{2}}] \qquad \text{and} \qquad e_{-\xi} = 
[e_{-\xi_{2}},e_{-\xi_{1}}]
\end{equation}
where $\xi_{1} <_{E} \xi <_{E} \xi_{2}$ (resp. $\xi_{1} <_{F} \xi <_{F} 
\xi_{2}$) and $(\xi_{1},\xi_{2})$ is a minimal segment in the chosen 
ordering.

Then the factors $\cR_{E}$ and $\cR_{F}$ of the universal R-matrix of 
$\dy{sl_{N}}$ have been conjectured to be \cite{KT96}:
\begin{align}
& \cR_{E} = \prod_{\gamma\in\Delta_{+},n\in\NN}^{\rightarrow} 
\exp(-e_{\gamma,n} \otimes f_{\gamma,-n-1}) \\
& \cR_{F} = \prod_{\gamma\in\Delta_{+},n\in\NN}^{\rightarrow} 
\exp(-f_{\gamma,n} \otimes e_{\gamma,-n-1})
\end{align}
where the arrows correspond to ordered products according to $<_{E}$ and 
$<_{F}$ respectively.

Using the following $N$-dimensional evaluation representation with 
evaluation parameter $u$
\begin{align}
\label{eq:repdysln}
& \pi_{u}(e_{i,k}) = E_{i,i+1} \; (u+\sfrac{1}{2}(i-1))^k \nonumber \\
& \pi_{u}(f_{i,k}) = E_{i+1,i} \; (u+\sfrac{1}{2}(i-1))^k \\
& \pi_{u}(h_{i,k}) = (E_{i,i} - E_{i+1,i+1}) \; (u+\sfrac{1}{2}(i-1))^k 
\nonumber
\end{align}
the represented R-matrix of $\dy{sl_{N}}$ takes the Yang form up to a 
normalization factor
\begin{equation}
\label{eq:Rdysln}
R[\dy{sl_{N}}](u) = \rho(u) \left( \sum_{1 \le a \le N} E_{aa} \otimes 
E_{aa} + \frac{1}{u+1} \sum_{1 \le a \ne b \le N} \; ( u \, E_{aa} 
\otimes E_{bb} + E_{ab} \otimes E_{ba} ) \right)
\end{equation}
with
\begin{equation}
\label{eq:normdysln}
\rho(u) = \frac{\Gamma_{1}(u \,|\, N) \; \Gamma_{1}(u+N \,|\, 
N)}{\Gamma_{1}(u+1 \,|\, N) \; \Gamma_{1}(u+N-1 \,|\, N)}
\end{equation}

\subsubsection{Central extension $\dy{\fg}_{c}$ of the double Yangian 
$\dy{\fg}$}

The double Yangian $\dy{\fg}$ admits a central extension $\dy{\fg}_{c}$ 
\cite{Kho95}. It is generated by the generators $e_{i,n}$, $f_{i,n}$ and 
$h_{i,n}$ where $i = 1,\ldots,r$ ($r$ is the rank of $\fg$) and $n \in 
\ZZ$, a central element $c$ and a grading operator $d$. In terms of the 
generating functionals (\ref{eq:genfun}), the relations between the 
generators read
\begin{equation}
\begin{split}
& e_{i}(u) \, e_{j}(v) = \frac{u-v+\alpha_{ij}}{u-v-\alpha_{ij}} \; 
e_{i}(v) \, e_{j}(u) \\
& f_{i}(u) \, f_{j}(v) = \frac{u-v-\alpha_{ij}}{u-v+\alpha_{ij}} \; 
f_{i}(v) \, f_{j}(u) \\
& h_{i}^{\pm}(u) \, e_{j}(v) = 
\frac{u-v+\alpha_{ij}}{u-v-\alpha_{ij}} \; e_{j}(v) \, 
h_{i}^{\pm}(u) \\
& h_{i}^{+}(u) \, f_{j}(v) = 
\frac{u-v-\alpha_{ij}-\alpha_{ij}c}{u-v+\alpha_{ij}-\alpha_{ij}c} 
\; f_{j}(v) \, h_{i}^{+}(u) \\
& h_{i}^{-}(u) \, f_{j}(v) = 
\frac{u-v-\alpha_{ij}}{u-v+\alpha_{ij}} \; f_{j}(v) \, h_{i}^{-}(u) 
\\
& h_{i}^{+}(u) \, h_{j}^{-}(v) = 
\frac{u-v+\alpha_{ij}}{u-v-\alpha_{ij}} \; 
\frac{u-v-\alpha_{ij}-\alpha_{ij}c}{u-v+\alpha_{ij}-\alpha_{ij}c} 
\; h_{j}^{-}(v) \, h_{i}^{+}(u) \\
& \big[ e_{i}(u) , f_{j}(v) \big] = \delta_{ij} \, \big( 
\delta(u-v-\alpha_{ii}c) h_{i}^{+}(u) - \delta(u-v) h_{i}^{-}(v) 
\big)
\end{split}
\end{equation}
together with the Serre-type relations (\ref{eq:serreyang}). The action of 
the grading operator is given by
\begin{equation}
\big[ d , e_{i}(u) \big] = \frac{d}{du} \; e_{i}(u) \qqquad 
\big[ d , f_{i}(u) \big] = \frac{d}{du} \; f_{i}(u) \qqquad 
\big[ d , h_{i}^{\pm}(u) \big] = \frac{d}{du} \; h_{i}^{\pm}(u)
\end{equation}

\medskip

As the double Yangian $\dy{\fg}$, the central extension of the double 
Yangian $\dy{\fg}_{c}$ exhibits a structure of quantum double \cite{Kho95}. 
Let $\widehat{Y}^{+} = Y^{+} \otimes \CC[c]$ and $\widehat{Y}^{-}$ be the 
semi-direct product of $Y^{-}$ with the ring $\CC[[d]]$ of formal power 
series in $d$. $\widehat{Y}^{\pm}$ are Hopf algebras with $\Delta(c) = c 
\otimes 1 + 1 \otimes c$ and $\Delta(d) = d \otimes 1 + 1 \otimes d$. The 
Hopf paring on $Y^{+} \otimes Y^{-}$ is uniquely extended onto 
$\widehat{Y}^{+} \otimes \widehat{Y}^{-}$ by $\langle c,d \rangle = 1$.

\medskip

In the $sl_{2}$ case, the universal R-matrix of $\dy{sl_{2}}_{c}$ 
decomposes as
\begin{equation}
\label{eq:RDYsl2c}
\cR[\dy{sl_{2}}_{c}] = \cR_{E} \; \cR_{H} \; \exp(c \otimes d) \; 
\cR_{F}
\end{equation}
where $\cR_{E}$, $\cR_{F}$ and $\cR_{H}$ are given by (\ref{eq:RERF}) and 
(\ref{eq:RH}).

In the two-dimensional evaluation representation with evaluation parameter 
$u$, $\pi_{ev}(u)(e_{k}) = E_{12} \, u^k$, $\pi_{ev}(u)(f_{k}) = E_{21} \, 
u^k$, $\pi_{ev}(u)(h_{k}) = (E_{11}-E_{22}) \, u^k$, $\pi_{ev}(u)(c) = 0$, 
the evaluated R-matrix $R_{12}(u_{1}-u_{2}) = \big(\pi_{ev}(u_{1}) \otimes 
\pi_{ev}(u_{2})\big) \, \cR[\dy{sl_{2}}_{c}]$ is identified with 
(\ref{eq:Rdysln}) with the normalization factor (\ref{eq:normdysln}) 
(taking $N=2$). In the FRT formalism, $\dy{sl_{2}}_{c}$ is then defined by 
the relations
\begin{equation}
\begin{split}
R_{12}(u_{1}-u_{2}) \, L_{1}^{\pm}(u_{1}) \, L_{2}^{\pm}(u_{2}) &= 
L_{2}^{\pm}(u_{2}) \, L_{1}^{\pm}(u_{2}) \, R_{12}(u_{1}-u_{2}) \\
R_{12}(u_{1}-u_{2}-c) \, L_{1}^-(u_{1}) \, L_{2}^+(u_{2}) &= 
L_{2}^+(u_{2}) \, L_{1}^-(u_{2}) \, R_{12}(u_{1}-u_{2}) \\
\end{split}
\end{equation}
where $L^+(u) = \sum_{k \ge 0} L^+_{k} u^{-k}$ and $L^-(u) = \sum_{k \le 0} 
L^-_{k} u^{-k}$ are two $2 \!\times\! 2$ independent matrices, expressed in 
term of the Drinfel'd generators $e^{\pm}(u)$, $f^{\pm}(u)$ and 
$h^{\pm}(u)$ using a Gauss decomposition of the Lax matrices:
\begin{equation}
L^{\pm}(u) = \left(
\begin{array}{cc} 
1 & f^{\pm}(u^\mp) \\ 
0 & 1 \\ 
\end{array} 
\right) 
\left( 
\begin{array}{cc} 
k_1^{\pm}(u) & 0 \\ 
0 & k_2^{\pm}(u) \\ 
\end{array} 
\right) 
\left( 
\begin{array}{cc} 
1 & 0 \\ 
e^{\pm}(u) & 1 \\ 
\end{array} 
\right) 
\end{equation}
with $u^+ = u$ and $u^- = u-c$. Furthermore, $k_1^{\pm}(u) k_2^{\pm}(u-1) = 1$ 
and one defines $h^{\pm}(u)\equiv k_2^{\pm}(u)^{-1} k_1^{\pm}(u)$.

\subsection{Deformed double Yangians}

\subsubsection{Definition of the deformed double Yangian $\ddy{sl_{2}}{r}$.}

Consider the R-matrix of $\elpa{2}$, and perform the scaling limit $q \to 
1$, with $z = q^u$, $p = q^{2r}$, keeping $u$ and $r$ fixed. One obtains, 
up to a similarity transformation, the following R-matrix:
\begin{equation}
\label{eq:RDYr}
R[\ddy{sl_{2}}{r}](u,r) = \rho(u,r) \left(
\begin{array}{cccc}
1 & 0 & 0 & 0 \\
0 & \sfrac{\sin \pi u/r}{\sin \pi(u+1)/r} & \sfrac{\sin \pi/r}{\sin 
\pi(u+1)/r} & 0 \\
0 & \sfrac{\sin \pi/r}{\sin \pi(u+1)/r} & \sfrac{\sin \pi u/r}{\sin 
\pi(u+1)/r} & 0 \\
0 & 0 & 0 & 1 \\
\end{array}
\right)
\end{equation}
the normalization factor being
\begin{equation}
\label{eq:factnorm}
\rho(u,r) = \frac{\Gamma_{2}(r+1-u \,|\, r,2)^2 \, \Gamma_{2}(2+u \,|\, 
r,2) \, \Gamma_{2}(u \,|\, r,2)}{\Gamma_{2}(u+1 \,|\, r,2)^2 \, 
\Gamma_{2}(r-u \,|\, r,2) \, \Gamma_{2}(r+2-u \,|\, r,2)}
\end{equation} 
where $\Gamma_2$ is the multiple Gamma function of order 2 (see Appendix A).

Taking now $L(u) = \sum_{n \in \ZZ} L_{n} u^{-n}$, the deformed double 
Yangian $\ddy{sl_{2}}{r}$ is defined by
\begin{equation}
R_{12}(u_{1}-u_{2},r) \, L_{1}(u_{1}) \, L_{2}(u_{2}) = L_{2}(u_{2}) \, 
L_{1}(u_{2}) \, R_{12}(u_{1}-u_{2},r-c)
\end{equation}

\subsubsection{The deformed double Yangian $\ddy{sl_{2}}{r}$ is a QTQHA}

What is the status of this algebra? In fact, in the same way the quantum 
affine elliptic algebra $\elpa{2}$ appears as a Drinfel'd twist of the 
universal quantum affine algebra $\cU_{q}(\widehat{sl}_{2})$, we sketch below 
how the deformed double Yangian $\ddy{sl_{2}}{r}$ can be obtained as a 
Drinfel'd twist of the double Yangian $\dy{sl_{2}}_{c}$, promoting 
$\ddy{sl_{2}}{r}$ as a QTQHA. For this purpose, one must prove that the 
matrix (\ref{eq:RDYr}) is indeed an evaluation representation of a 
universal R-matrix obtained as a Drinfel'd twist of the universal R-matrix 
of $\dy{sl_{2}}$.

Introducing the notation 
\begin{equation}
\label{eq:Mbpm}
M(b^+,b^-) = \left(
\begin{array}{cccc}
1 & 0 & 0 & 0 \\
0 & \sfrac{1}{2}\,(b^++b^-) & \sfrac{1}{2}\,(b^+-b^-) & 0 \\[2mm]
0 & \sfrac{1}{2}\,(b^+-b^-) & \sfrac{1}{2}\,(b^++b^-) & 0 \\
0 & 0 & 0 & 1 \\
\end{array}
\right)
\end{equation}
one has obviously $M(a,b) M(a',b')= M(aa',bb')$ and $M(a,b)^{-1} = 
M(a^{-1},b^{-1})$. The R-matrix of $\dy{sl_{2}}$ can be written as 
$R[\dy{sl_{2}}](u) = \rho(u) M \left(1,\frac{u-1}{u+1} \right)$, and the 
R-matrix of $\ddy{sl_{2}}{r}$ takes the form $R[\ddy{sl_{2}}{r}](u) = 
\rho_r(u) M(b_r^+,b_r^-)$, with
\begin{align*}
b_r^+ &= \frac{\cos\frac{u-1}{2r}} {\cos\frac{u+1}{2r}} = 
\frac{\Gamma_1(r+u+1 \,|\, 2r) \Gamma_1(r-u-1 \,|\, 2r)} 
{\Gamma_1(r+u-1 \,|\, 2r) \Gamma_1(r-u+1 \,|\, 2r)} \\
b_r^- &= \frac{\sin\frac{u-1}{2r}} {\sin\frac{u+1}{2r}} = 
\frac{\Gamma_1(u+1 \,|\, 2r) \Gamma_1(2r-u-1 \,|\, 2r)} {\Gamma_1(u-1 
\,|\, 2r) \Gamma_1(2r-u+1 \,|\, 2r)} 
\end{align*}
The normalization factor of $\ddy{sl_{2}}{r}$ being rewritten as $\rho_r(u) 
= \rho_F(-u;r) \rho(u)$ $\rho_F(u;r)^{-1}$ with $\displaystyle \rho_F(u) = 
\frac{\Gamma_{2}(u+1+r \,|\, 2,r)^2} {\Gamma_{2}(u+r \,|\, 2,r) \; 
\Gamma_{2}(u+2+r \,|\, 2,r)}$, one obtains
\begin{equation}
R[\ddy{sl_{2}}{r}] = F_{21}(-u) R[\dy{sl_{2}}] F_{12}(u)^{-1}
\end{equation}
Using the notation (\ref{eq:Mbpm}), $F_{12}(u)$ is given by
\begin{equation}
F_{12}(u) = \rho_F(u) \, M \left( \frac{\Gamma_1(u+r-1 \,|\, 2r)} 
{\Gamma_1(u+r+1 \,|\, 2r)}, \frac{\Gamma_1(u+2r-1 \,|\, 2r)} 
{\Gamma_1(u+2r+1 \,|\, 2r)} \right)
\end{equation}
This twist-like matrix reads 
\begin{align}
F_{12}(u) &= \rho_F(u) \prod_{n=1}^\infty M \left( 1,\frac{u+1+2nr}{u 
-1+2n r} \right) M \left( \frac{u+1+(2n-1)r}{u-1+(2n-1)r},1 \right) 
\nonumber \\
&= \prod_{n=1}^\infty R(u+2nr)^{-1} \; \tau(R(u+(2n-1)r)^{-1})
\end{align}
where $\tau(M(a,b)) = M(b,a)$, $R$ is the R-matrix of $\dy{sl_{2}}_{c}$ and 
one uses the representation of $\rho_F(u)$ as an infinite product
\begin{equation}
  \rho_F(u) = \prod_{n=1}^\infty \rho(u+nr)^{-1} 
\end{equation}
The automorphism $\tau$ may be interpreted as the adjoint action of 
$(-1)^{\half h_{0}} \otimes 1$, where $h_{0}$ is the Cartan generator of 
$sl_{2} \subset \dy{sl_{2}}$, so that
\begin{align}
F_{12}(u) &= \prod_{n=1}^\infty R(u+2nr)^{-1} \Ad\Big( (-1)^{\half 
h_{0}} \otimes 1 \Big) R(u+(2n-1)r)^{-1} \nonumber \\
&= \prod_{n=1}^\infty \Ad\Big( (-1)^{\frac{n}{2} h_{0}} \otimes 1 
\Big) R(u+nr)^{-1}
\end{align}
Hence $F$ is solution of the difference equation
\begin{equation}
F_{12}(u) = \Ad\left( (-1)^{-\half h_{0}} \otimes 1 \right) F_{12}(u-r) 
\cdot R_{12}(u)
\end{equation}
Note that all the infinite products are logarithmically divergent. They are 
consistently regularized by the $\Gamma_1$ and $\Gamma_2$ functions. In 
particular, $\lim\limits_{r\rightarrow \infty}F = M(1,1) = \II_4$. \\
One imposes then the difference equation \emph{at the universal level}:
\begin{equation}
\label{eq:eqdiffuniv}
\cF_{12} = \Ad(\phi \otimes 1) (\cF_{12}) \cdot \cC
\end{equation}
where $\phi = (-1)^{-\half h_{0}} e^{-(r+c)d}$ and $\cC = e^{\half(c\otimes 
d + d\otimes c)} \cR$. A solution of (\ref{eq:eqdiffuniv}) is given by an 
infinite product
\begin{equation}
\label{eq:prodinf}
\cF_{12}(r) = \prod_{k \ge 1}^{\leftarrow} \cF_{k}(r) = \prod_{k \ge 
1}^{\leftarrow} \Ad(\phi^{-k} \otimes 1) (\cC_{12}^{-1})
\end{equation}
It can be proved along the same lines as eqs. (\ref{eq:DeltaidF}) to 
(\ref{eq:cocinduc}) that $\cF_{12}(r)$ satisfies the shifted cocycle 
condition. \\
Finally, one gets the following theorem \cite{AAFRR3}:
\begin{theorem}
\label{thm:yangr}
The Yangian-type Drinfel'd twist (\ref{eq:prodinf}) satisfies the 
shifted cocycle condition
\begin{equation*}
\cF_{12}(r) \, (\Delta \otimes \id) \cF(r) = \cF_{23}(r+c^{(1)}) \, 
(\id \otimes \Delta) \cF(r)
\end{equation*}
The deformed double Yangian $\ddy{sl_{2}}{r}$ is a quasi-triangular 
quasi-Hopf algebra with the universal R-matrix $\widetilde{\cR}(r) = 
\cF_{21}(r) \, \cR_{12} \, \cF_{12}^{-1}(r)$ where $\cR$ is the 
universal R-matrix of $\dy{sl_{2}}$. $\widetilde{\cR}(r)$ satisfies 
the shifted Yang--Baxter equation:
\begin{equation*}
\widetilde{\cR}_{12}(r+c^{(3)}) \, \widetilde{\cR}_{13}(r) \, 
\widetilde{\cR}_{23}(r+c^{(1)}) = \widetilde{\cR}_{23}(r) \, 
\widetilde{\cR}_{13}(r+c^{(2)}) \, \widetilde{\cR}_{12}(r)
\end{equation*}
\end{theorem}

\begin{remark}
Note that the evaluated R-matrices (\ref{eq:ruqa}) and (\ref{eq:RDYr}) 
are homothetical once the identification $q = e^{i\pi/r}$ and $z = 
e^{2i\pi u/r}$ is done. However they are used to construct distinct 
algebras which differ fundamentaly in their structure. The R-matrix 
(\ref{eq:ruqa}) is the evaluation of the universal R-matrix of the 
\emph{Hopf} algebra $\cU_{q}(\widehat{sl}_{2})$, while the R-matrix 
(\ref{eq:RDYr}) corresponds to the evaluation of the universal R-matrix 
of the \emph{quasi-Hopf} algebra $\ddy{sl_{2}}{r}$. In particular, the 
normalization factors of (\ref{eq:ruqa}) and (\ref{eq:RDYr}), which are 
different, are related roughly speaking to the contribution of the 
Cartan part of the corresponding universal R-matrices. \finrmq
\end{remark}

\subsubsection{Other presentation of $\ddy{sl_{2}}{r}$.}

As it was said before, the R-matrix $\overline{R}[\ddy{sl_{2}}{r}]$ 
obtained from the scaling limit of the R-matrix of $\elpa{2}$ differs from 
the R-matrix $R[\ddy{sl_{2}}{r}]$ given by Eq. (\ref{eq:RDYr}) by a 
similarity transformation. Indeed, one has
\begin{equation}
\overline{R}[\ddy{sl_{2}}{r}] = K_{21} \; R[\ddy{sl_{2}}{r}] \; 
K_{12}^{-1}
\label{eq:bleue}
\end{equation}
where
\begin{equation}
K = V \otimes V \qquad \mbox{with} \qquad V = \frac{1}{\sqrt{2}} \left( 
\begin{array}{rr} 1 & 1 \\ -1 & 1 \\ \end{array} \right)
\end{equation}
The Lax operators $L$ and $\overline{L}$ associated respectively to the 
R-matrices $R[\ddy{sl_{2}}{r}]$ and $\overline{R}[\ddy{sl_{2}}{r}]$ are 
connected by $\overline{L} = V L V^{-1}$, implying an isomorphism between 
the two corresponding Yangian structures.

One can identify $V$ with an evaluation representation of an element $g$: 
$V = \pi_{ev}(g)$ with $g = \exp\left(\frac{-\pi}{4}(f_0-e_0)\right)$. 
Since $e_0$ and $f_0$ lie in the undeformed Hopf subalgebra $\cU(sl_{2})$ 
of $\dy{sl_{2}}$, the coproduct of $g$ reads $\Delta(g) = g \otimes g$, so 
that
\begin{equation}
g_1 g_2 \, \Delta^\cF(g^{-1}) \cF = g_1 g_2 \, \cF g_1^{-1} g_2^{-1}
\label{eq:gG}
\end{equation}
where $\Delta$ is the coproduct of $\dy{sl_{2}}$ while $\Delta^\cF$ is the 
coproduct of $\ddy{sl_{2}}{r}$. \\
The two-cocycle $g_1 g_2 \Delta^\cF(g^{-1})$ is a coboundary (with respect 
to the coproduct $\Delta^\cF$). In representation, (\ref{eq:gG}) is equal 
to the scaling limit of the represented twist from $\cU_{q}(sl_{2})$ to 
$\elpa{2}$. \\
It follows that \cite{AAFRR3}
\begin{equation}
\overline{\cR}[\ddy{sl_{2}}{r}] \equiv g_1 g_2 \, 
\Delta_{21}^\cF(g^{-1}) \; \cR[\ddy{sl_{2}}{r}] \; \Delta_{12}^\cF(g) 
\, g_1^{-1} g_2^{-1} = g_1 g_2 \; \cR[\ddy{sl_{2}}{r}] \; g_1^{-1} 
g_2^{-1}
\end{equation}
satisfies the shifted Yang--Baxter equation of Theorem \ref{thm:yangr}.  \\
To recover (\ref{eq:bleue}), it suffices to use (\ref{eq:gG}) and to remark 
that $(\pi_{ev} \otimes \pi_{ev}) (g \otimes g)$ commutes with 
$R[\ddy{sl_{2}}{r}]$.

\subsubsection{Generalization to $\ddy{sl_{N}}{r}$}

\textbf{Definition of the deformed double Yangian $\ddy{sl_{N}}{r}$}

We now construct the deformed double Yangian $\ddy{sl_{N}}{r}$ 
\cite{AAFRo}. The R-matrix of the deformed double Yangian $\ddy{sl_{N}}{r}$ 
is obtained by taking the scaling limit ($q \to 1$ with $z = q^u$ and $p = 
q^{2r}$ and keeping $u,r$ fixed) of the R-matrix of the quantum affine 
elliptic algebra $\elpa{N}$. One gets
\begin{align}
\overline{R}(u,r) &= \rho_{DYr}(u) \times \nonumber \\
& \sum_{a,b,c=1}^{N} \frac{\displaystyle \sin\frac{\pi u}{r} \; 
\sin\frac{\pi}{r} \; \sin\frac{\pi}{Nr} (u+1+(b-a) r)} {\displaystyle N 
\sin\frac{\pi(u+1)}{r} \; \sin\frac{\pi}{Nr}(u+(b-c) r) \; 
\sin\frac{\pi}{Nr}(1-(a-c) r)} \; E_{ac} \otimes E_{b,a+b-c} 
\end{align}
the normalization factor $\rho_{DYr}(u)$ being defined by
\begin{equation}
\rho_{DYr}(u) = \frac{S_{2}(-u \,|\, r,N) \, S_{2}(1+u \,|\, r,N)} 
{S_{2}(u \,|\, r,N) \, S_{2}(1-u \,|\, r,N)}
\end{equation}
where $S_{2}(x \,|\, \omega_{1},\omega_{2})$ is Barnes' double sine 
function of periods $\omega_{1}$ and $\omega_{2}$ (see Appendix A). \\
It is possible to simplify this matrix by a similarity transformation. 
Defining $V_i^j = 1/\sqrt{N} \; \omega^{(i-1)j}$ where $\omega = 
\exp(2i\pi/N)$, the similar matrix $R = (V\otimes V) \overline{R} (V\otimes 
V)^{-1}$ has the following non-zero entries
\begin{alignat}{3}
\label{eq:matrixS}
& S_{aa}^{aa}(u) = \cot \frac{\pi u}{r} + \cot \frac{\pi}{r} && && 
\nonumber \\
& S_{ab}^{ab}(u) = \frac{e^{i\pi /r}}{\sin \frac{\pi}{r}} e^{- 2i\pi 
(b-a) /Nr} && \qqquad && \text{for} \;\; b-a\in\{1,...,N-1\} \\
& S_{ab}^{ba}(u) = \frac{e^{i\pi u/r}}{\sin \frac{\pi u}{r}} e^{- 2i\pi 
(b-a) u /Nr} && \qqquad && \text{for} \;\; b-a\in\{1,...,N-1\} 
\nonumber
\end{alignat}
where the matrix elements $R_{ij}^{kl}$ and $S_{ij}^{kl}$ are related by
\begin{equation}
\label{eq:matrixR}
R_{ij}^{kl} = - \rho_{DYr}(u) \; \frac{\displaystyle \sin\frac{\pi 
u}{r} \; \sin\frac{\pi}{r}} {\displaystyle \sin\frac{\pi(u+1)}{r}} \; 
S_{ij}^{kl}
\end{equation}

\textbf{The deformed double Yangian $\ddy{sl_{N}}{r}$ is a QTQHA}

Inspired by the expressions obtained in the $sl_{2}$ case, one postulates 
\cite{AAFR} the linear equation (\ref{eq:eqdiffuniv}) in $\dy{sl_{N}}_{c} 
\otimes \dy{sl_{N}}_{c}$ for the twist $\cF_{12}$, where $\phi = 
\omega^{-h_{0,\rho}} e^{-(r+c)d}$ and $\cC = e^{\frac{1}{2} (c \otimes d + 
d \otimes c)} \cR$. Equation (\ref{eq:eqdiffuniv}) can be solved by
\begin{equation}
\label{eq:FprodDYr}
\cF_{12}(r) = \prod_{k \ge 1}^{\leftarrow} \cF_{k}(r) = \prod_{k \ge 
1}^{\leftarrow} \Ad(\phi^{-k} \otimes 1) (\cC_{12}^{-1})
\end{equation}
The operator $d$ in the double Yangian $\dy{sl_{N}}_{c}$ satisfies 
$\displaystyle \big[ d,e_\alpha(u) \big] = \frac{d}{du} \; e_\alpha(u)$ for 
any root $\alpha$ and its coproduct is given by $\Delta(d) = d \otimes 1 + 
1 \otimes d$ (see \cite{Kho95}). \\
The generator $h_{0,\rho}$ of $\dy{sl_{N}}_{c}$ is such that
\begin{equation}
h_{0,\rho} e_\alpha(u) = e_\alpha(u)(h_{0,\rho}+(\rho,\alpha)) \;, 
\qquad h_{0,\rho} f_\alpha(u) = f_\alpha(u)(h_{0,\rho}-(\rho,\alpha )) 
\;, \qquad \big[ h_{0,\rho}, h_\alpha(u) \big] = 0
\end{equation}
and hence $\tau= \Ad\Big( \omega^{-h_{0,\rho}} \otimes \id \Big)$ is 
idempotent, since all the scalar products $(\rho,\alpha)$ are rational. \\
As in the $sl_{2}$ case, equations (\ref{eq:DeltaidF}) to 
(\ref{eq:cocinduc}) still hold and thus the infinite product 
(\ref{eq:FprodDYr}) satisfies the shifted cocycle relation. The twist 
$\cF(r)$ defines a QTQHA denoted $\ddy{sl_{N}}{r}$ with universal R-matrix 
$\cR[\ddy{sl_{N}}{r}](r) = \cF_{21}(r) \; \cR[\dy{sl_{N}}_{c}] \; 
\cF_{12}(r)^{-1}$.

\medskip

It remains to show that the corresponding evaluated R-matrix indeed 
coincides with the expression (\ref{eq:matrixR}). In the evaluation 
representation (\ref{eq:repdysln}), the linear equation 
(\ref{eq:eqdiffuniv}) takes the form
\begin{equation}
\label{eq:eqlinrep}
F(u) = (\phi \otimes 1)^{-1} F(u-r) (\phi \otimes 1) R(u)
\end{equation}
where $\phi = \text{diag}(\omega^{N-1},\omega^{N-2},\ldots,\omega,1)$. \\
The solution of (\ref{eq:eqlinrep}) is expressed in term of hypergeometric 
functions ${_{2}F_{1}}$:
\begin{align}
& F_{ab}^{ab}(u) = \frac{\Gamma\left(\frac{u-1}{r}+1\right)} 
{\Gamma\left(\frac{u}{r}+1\right)} \; 
\hypergeom{-\frac1r}{\frac{u-1}{r}+1}{\frac{u}{r}+1}{\omega^{b-a}} \\
& F_{ab}^{ba}(u) = -\frac{\omega^{b-a}}{r} \; 
\frac{\Gamma\left(\frac{u-1}{r}+1\right)} 
{\Gamma\left(\frac{u}{r}+2\right)} \; 
\hypergeom{-\frac1r+1}{\frac{u-1}{r}+1}{\frac{u}{r}+2}{\omega^{b-a}}
\end{align}
and $F_{aa}^{aa}(u) = 1$ (the indices $a \ne b$ run from 1 to $N$). \\
The twist $F(u)$ applied to the R-matrix of $\dy{sl_{N}}$, eq. 
(\ref{eq:Rdysln}), provides the R-matrix $\widetilde{R}(u)$ of the deformed 
double Yangian $\ddy{sl_{N}}{r}$, eq. (\ref{eq:matrixR}), the non-vanishing 
entries of which are expressed in terms of (\ref{eq:matrixS}). The proof 
follows by a direct computation using the properties of the hypergeometric 
functions $_{2}F_{1}$.

\medskip

This kind of construction can be generalized for any simple Lie algebra 
$\fg$. 

\subsection{Dynamical double Yangians}

\subsubsection{The dynamical double Yangian $\ddy{sl_{2}}{s}$}

Using the same kind of argument, it is possible to construct another 
deformation of the Yangian $\dy{sl_{2}}$, which can be considered as a 
dynamical deformation, denoted by $\ddy{sl_{2}}{s}$. The corresponding 
R-matrix reads:
\begin{equation}
R[\ddy{sl_{2}}{s}](u,s) = \rho(u) \left(
\begin{array}{cccc}
1 & 0 & 0 & 0 \\
0 & \sfrac{u}{u+1} & \sfrac{s+u}{s(u+1)} & 0 \\
0 & \sfrac{s-u}{s(u+1)} & \sfrac{s^2-1}{s^2} \, \sfrac{u}{u+1} & 0 
\\
0 & 0 & 0 & 1 \\
\end{array}
\right)
\end{equation}
$\rho(u)$ being the normalization factor of the R-matrix of $\dy{sl_{2}}$ 
and $s$ is a complex number. Actually this matrix can be obtained as a 
scaling limit of the quantum affine elliptic algebra $\elpb{2}$ with $p=0$.

The algebra $\ddy{sl_{2}}{s}$ is then defined by the RLL-type relations
\begin{equation}
\begin{split}
R_{12}(u_1-u_2,\lambda+h) \, & L_1(u_1,\lambda) \, 
L_2(u_2,\lambda+h^{(1)}) = \\
& L_2(u_2,\lambda) \, L_1(u_1,\lambda+h^{(2)}) \, 
R_{12}(u_1-u_2,\lambda)
\end{split}
\end{equation}

Let us prove that $\ddy{sl_{2}}{s}$ is a QTQHA obtained from $\dy{sl_{2}}$ 
by a Drinfel'd twist. To this aim, one uses the fact that $\cU(sl_{2})$ is 
a Hopf subalgebra of $\dy{sl_{2}}$.

One first constructs a twist from ${\cU}_{q}(sl_{2})$ to 
$\cB_{q,\lambda}(sl_{2})$ \cite{ABRR,Bab91,JKOS} as follows. Introduce 
$\phi = q^\chi = q^{\half h^2 + sh}$ where $s \in \CC$, and consider
\begin{equation}
\label{eq:twistbql2}
\cF = \prod_{k \ge 1}^{\leftarrow} \Ad (\phi^k \otimes 1) \; 
(\widehat{\cR}^{-1})
\end{equation}
where $\cR = q^{-\half h \otimes h} \; \widehat{\cR}$ is the universal 
R-matrix of ${\cU}_{q}(sl_{2})$. The twisted R-matrix $\widetilde{\cR}_{12} 
= \cF_{21} \; \cR_{12} \; \cF_{12}^{-1}$ defines the algebra 
$\cB_{q,\lambda}(sl_{2})$. The twist $\cF$ satisfies the following 
difference equation
\begin{equation}
\label{eq:eqlinbql}
\cF = \Ad(\phi^{-1} \otimes 1) (\cF) \cdot \widehat{\cR}
\end{equation}
the solution of which is precisely the infinite product (\ref{eq:twistbql2}).

Consider now the twist $\cF\,'$ obtained from (\ref{eq:twistbql2}) by 
taking the limit $q \to 1$ and keeping $s$ fixed. The consistency of this 
procedure follows from the Hopf algebra identification 
$\cU_{\hbar}(sl_{2})/(\hbar\;\cU_{\hbar}(sl_{2})) \simeq \cU(sl_{2})$ with 
$q = e^{\hbar}$. In this quotient, the twist $\cF\,'$ is given by the 
formula
\begin{equation}
\label{eq:twistdys}
\cF\,' = \sum_{n=0}^{\infty} \frac{1}{n!} \left( \prod_{k=0}^{n-1} [ 
(1+k-s)1 \otimes h ] \otimes 1 \right)^{-1} \; e^n \otimes f^n
\end{equation}
In fact, the scaling limit of the difference equation (\ref{eq:eqlinbql}) 
of the twist $\cF$ leads to the following equation for the twist $\cF\,'$, 
obtained as the first non-trivial order in $\hbar$ in (\ref{eq:eqlinbql}):
\begin{equation}
[ \chi \otimes 1 , \cF\,' ] = \cF\,' \widehat{r}
\end{equation}
where $\widehat{r}$ is given by $\widehat{\cR} = 1 \otimes 1 + \hbar 
\widehat{r} + o(\hbar)$. This equation has a unique solution expressed 
either by eq. (\ref{eq:twistdys}) or by the infinite product (it can be 
checked directly, order by order, that formulae (\ref{eq:twistdys}) and 
(\ref{eq:twistdysprod}) lead to the same expression in terms of the 
generators)
\begin{equation}
\label{eq:twistdysprod}
\cF\,' = \prod_{k \ge 1}^{\leftarrow} \Ad(\chi^{-k} \otimes 1) \Big( 1 
\otimes 1 + (\chi \otimes 1)^{-1} \widehat{r} \Big)
\end{equation}
Finally, the following theorem holds \cite{AAFR,Cla}:
\begin{theorem}
The dynamical Yangian-type Drinfel'd twist (\ref{eq:twistdysprod}) 
satisfies the shifted cocycle condition
\begin{equation*}
\cF\,'_{12}(s) \, (\Delta \otimes \id) \cF\,'(s) = \cF\,'_{23}(s + 
{h^\vee}^{(1)}) \, (\id \otimes \Delta) \cF\,'(s)
\end{equation*}
The dynamical double Yangian $\ddy{sl_{2}}{s}$ is a quasi-triangular 
quasi-Hopf algebra with the universal R-matrix $\widetilde{\cR}(s) = 
\cF\,'_{21}(s) \; \cR \; \cF\,'_{12}(s)$ where $\cR$ is the universal 
R-matrix of $\dy{sl_{2}}$. $\widetilde{\cR}(s)$ satisfies the 
dynamical Yang--Baxter equation:
\begin{equation*}
\widetilde{\cR}_{12}(s+{h^\vee}^{(3)}) \, \widetilde{\cR}_{13}(s) 
\, \widetilde{\cR}_{23}(s+{h^\vee}^{(1)}) = \widetilde{\cR}_{23}(s) 
\, \widetilde{\cR}_{13}(s+{h^\vee}^{(2)}) \, 
\widetilde{\cR}_{12}(s)
\end{equation*}
\end{theorem}
In that case, although the computation is rather tedious, the proof of the 
shifted cocycle condition can be done directly by using the sum formula 
(\ref{eq:twistdys}).

\subsubsection{Generalization to $\ddy{sl_{N}}{s}$}

The previous construction can be generalized to the case of $sl_{N}$ 
without difficulty \cite{AAFR}. It follows the same lines: \\
-- construction of the twist $\cF$ from ${\cU}_{q}(sl_{N})$ to 
$\cB_{q,\lambda}(sl_{N})$, \\
-- one takes the scaling limit of this twist consistently from the Hopf 
algebra identification 
\begin{equation}
\cU_{\hbar}(sl_{N})/(\hbar\;\cU_{\hbar}(sl_{N})) \simeq \cU(sl_{N})
\end{equation}
with $q = e^{\hbar}$, \\
-- the scaled twist $\cF\,'$ can also be obtained in its universal form as 
an infinite product given by (\ref{eq:twistdysprod}), the matrix 
$\widehat{r}$ being the classical limit of the matrix $\widehat{\cR}$ where 
$\cR = q^{- d_{ij} h_{i} \otimes h_{j}} \; \widehat{\cR}$ is the universal 
R-matrix of ${\cU}_{q}(sl_{N})$ and $d_{ij}$ is the inverse of the Cartan 
matrix of $sl_{N}$.

\medskip

In the fundamental representation for $sl_{N}$, the evaluated infinite 
product expression for $\cF\,'$ reads
\begin{equation}
F' = 1 \otimes 1 - \sum_{a<b} \frac{2}{x_{a}-x_{b}} \; E_{ab} \otimes 
E_{ba}
\end{equation}
Applied to the R-matrix of $\dy{sl_{N}}$ given in (\ref{eq:Rdysln}), this 
twist leads to the evaluated R-matrix of $\ddy{sl_{N}}{s}$:
\begin{align}
R[\ddy{sl_{N}}{s}](u,s) = \rho[\ddy{sl_{N}}{s}](u) & \Bigg( \sum_{a} 
E_{aa} \otimes E_{aa} + \sum_{a<b} \frac{u}{u+1} \; E_{aa} \otimes 
E_{bb} \nonumber \\
&+ \sum_{a>b} \left( 1 - \frac{4}{(x_{a}-x_{b})^2} \right) \; 
\frac{u}{u+1} \; E_{aa} \otimes E_{bb} \nonumber \\
&+ \sum_{a,b} \left( 1 + \frac{2u}{x_{a}-x_{b}} \right) \frac{1}{u+1} 
\; E_{ab} \otimes E_{ba} \Bigg)
\end{align}
the normalization factor being $\rho[\ddy{sl_{N}}{s}](u) = 
\rho[\dy{sl_{N}}](u)$, see eq. (\ref{eq:normdysln}).

\subsubsection{The deformed dynamical double Yangian $\ddy{sl_{2}}{r,s}$}

Considering now the scaling limit of the quantum affine elliptic algebra 
$\elpb{2}$, it is possible to construct a deformed dynamical double Yangian 
$\ddy{sl_{2}}{r,s}$ with the following R-matrix \cite{Cla}
\begin{equation}
\label{eq:Rdyrs}
R[\ddy{sl_{2}}{r,s}](u,r,s) = \rho(u,r) \left(
\begin{array}{cccc}
1 & 0 & 0 & 0 \\
0 & b(u) & c(u) & 0 \\
0 & c'(u) & b'(u) & 0 \\
0 & 0 & 0 & 1 \\
\end{array} 
\right)
\end{equation}
where
\begin{equation}
\begin{split}
& b(u) = \frac{\Gamma_{1}(r-s \,|\, r)^2} {\Gamma_{1}(r-s+1 \,|\, 
r) \, \Gamma_{1}(r-s-1 \,|\, r)} \; \frac{\sin \frac{\pi u}{r}} 
{\sin\frac{\pi(1+u)}{r}} \\
& c(u) = \frac{\sin\frac{\pi(s+u)}{r}}{\sin\frac{\pi s}{r}} \; 
\frac{\sin\frac{\pi}{r}}{\sin\frac{\pi(1+u)}{r}} \\
& b'(u) = \frac{\Gamma_{1}(s \,|\, r)^2} {\Gamma_{1}(s+1 \,|\, r) 
\, \Gamma_{1}(s-1 \,|\, r)} \; \frac{\sin\frac{\pi 
u}{r}}{\sin\frac{\pi(1+u)}{r}} \\
& c'(u) = \frac{\sin\frac{\pi(s-u)}{r}}{\sin\frac{\pi s}{r}} \; 
\frac{\sin\frac{\pi}{r}}{\sin\frac{\pi(1+u)}{r}}
\end{split}
\end{equation}
and the normalization factor $\rho(u,r)$ is given by eq. 
(\ref{eq:factnorm}). \\
The algebra $\ddy{sl_{2}}{r,s}$ is then defined by the relations
\begin{equation}
\begin{split}
R_{12}(u_1-u_2,\lambda+h) \, & L_1(u_1,\lambda) \, 
L_2(u_2,\lambda+h^{(1)}) = \\
& L_2(u_2,\lambda) \, L_1(u_1,\lambda+h^{(2)}) \, 
R_{12}(u_1-u_2,\lambda)
\end{split}
\end{equation}
In that case, although the algebra can be defined through RLL relations, 
the status of this algebra as a QTQHA, obtained by a Drinfel'd twist of 
$\dy{sl_{2}}_{c}$, remains still open. However the interest of this algebra 
is enforced by the following remark.

\begin{remark}
Considering the R-matrix of the $\ddy{sl_{2}}{r,s}$ algebra, eq. 
(\ref{eq:Rdyrs}), and taking the limit $s \to i\infty$, one gets a 
non-dynamical R-matrix given by
\begin{equation}
\label{eq:RDYrF}
R = \rho(u,r) \left(
\begin{array}{cccc}
1 & 0 & 0 & 0 \\
0 & \sfrac{\sin \pi u/r}{\sin \pi(u+1)/r} & e^{-i\pi u} \, 
\sfrac{\sin \pi/r}{\sin \pi(u+1)/r} & 0 \\
0 & e^{i\pi u} \, \sfrac{\sin \pi/r}{\sin \pi(u+1)/r} & 
\sfrac{\sin \pi u/r}{\sin \pi(u+1)/r} & 0 \\
0 & 0 & 0 & 1 \\
\end{array}
\right)
\end{equation}
where $\rho(u,r)$ is the same as above. \\
In fact, this matrix can be obtained from the R-matrix (\ref{eq:RDYr}) 
by a similarity transformation. This similarity transformation can be 
constructed as a Drinfel'd twist at the universal level. Considering 
$g' = \exp(h_{1}/2r)$, one defines the following shifted 
\emph{coboundary}
\begin{equation}
\cK_{12}(r) = g'(r) \otimes g'(r+c^{(1)}) \; \Delta^\cF({g'}^{-1})
\end{equation}
where $\Delta^\cF$ is the coproduct of $\ddy{sl_{2}}{r}$. It obeys a 
shifted cocycle condition
\begin{equation}
\cK_{12}(r) \; (\Delta^\cF \otimes \id) \cK(r) = 
\cK_{23}(r+c^{(1)}) \; (\id \otimes \Delta^{\cF'}) \cK(r)
\end{equation}
as a consequence of
\begin{equation}
(\Delta^\cF \otimes \id) \; \Delta^\cF({g'}^{-1}) = (\id \otimes 
\Delta^{\cF'}) \; \Delta^\cF({g'}^{-1})
\end{equation}
with $\cF'_{23}(r) = \cF_{23}(r+c^{(1)})$. This last equation is the 
quasi-coassociativity property for the coproduct $\Delta^\cF$. \\
Finally, the matrix $\cK_{21}(r) \, \cR[\ddy{sl_{2}}{r}] \, 
\cK_{12}^{-1}(r)$ satisfies the shifted Yang--Baxter equation, and the 
corresponding evaluated R-matrix coincides with (\ref{eq:RDYrF}). 
\finrmq
\end{remark}

\bigskip

\textbf{Acknowledgements:} The author would like to thank the organizers of 
the School for their kind invitation and support. This work greatly 
benefitted from a careful reading of the manuscript and helpful comments of 
D. Arnaudon and E. Ragoucy. The author warmfully thanks them.

\clearpage

\markboth{Appendix}{Appendix}
\appendix
\section*{Appendix}
\section{Notations}
\setcounter{equation}{0}

\subsubsection*{Jacobi theta functions}

Let $\HH = \{ z\in\CC \,\vert\, \text{Im} z > 0 \}$ be the upper half-plane 
and $\Lambda_\tau = \{ \lambda_1\,\tau + \lambda_2 \,\vert\, 
\lambda_1,\lambda_2 \in \ZZ \,, \tau \in \HH \}$ the lattice with basis 
$(1,\tau)$ in the complex plane. One denotes the congruence ring modulo 
$N$ by $\ZZ_N \equiv \ZZ/N\ZZ$ with basis $\{0,1,\dots,N-1\}$. One sets 
$\omega = e^{2i\pi/N}$. Finally, for any pairs 
$\gamma=(\gamma_1,\gamma_2)$ and $\lambda=(\lambda_1,\lambda_2)$ of 
numbers, we define the (skew-symmetric) pairing 
$\langle\gamma,\lambda\rangle \equiv \gamma_1\lambda_2 - 
\gamma_2\lambda_1$.

\medskip

\noindent
One defines the Jacobi theta functions with rational characteristics
$\gamma=(\gamma_1,\gamma_2) \in \sfrac{1}{N} \ZZ \times \sfrac{1}{N} \ZZ$ by
\begin{equation}
\vartheta\car{\gamma_1}{\gamma_2}(\xi,\tau) = \sum_{m \in \ZZ} 
\exp\big(i\pi(m+\gamma_1)^2\tau + 2i\pi(m+\gamma_1)(\xi+\gamma_2)\big) 
\end{equation}
The functions $\vartheta\car{\gamma_1}{\gamma_2}(\xi,\tau)$ satisfy the
following shift properties
\begin{align}
& \vartheta\car{\gamma_1+\lambda_1}{\gamma_2+\lambda_2}(\xi,\tau) = 
\exp(2i\pi \gamma_1\lambda_2) \,\, 
\vartheta\car{\gamma_1}{\gamma_2}(\xi,\tau) \\
& \vartheta\car{\gamma_1}{\gamma_2}(\xi+\lambda_1\tau+\lambda_2,\tau) 
= \exp(-i\pi\lambda_1^2\tau-2i\pi\lambda_1\xi) \, 
\exp(2i\pi\langle\gamma,\lambda\rangle) \, 
\vartheta\car{\gamma_1}{\gamma_2}(\xi,\tau) 
\end{align}
where $\gamma=(\gamma_1,\gamma_2) \in \sfrac{1}{N}\ZZ \times 
\sfrac{1}{N}\ZZ$ and $\lambda=(\lambda_1,\lambda_2) \in \ZZ \times \ZZ$. 
\\
Moreover, for arbitrary $\lambda=(\lambda_1,\lambda_2)$ (not necessarily 
integers), one has the following shift exchange
\begin{equation}
\vartheta\car{\gamma_1}{\gamma_2}(\xi+\lambda_1\tau+\lambda_2,\tau) = 
\exp\big(-i\pi\lambda_1^2\tau-2i\pi\lambda_1(\xi+\gamma_2+\lambda_2)\big) 
\, \vartheta\car{\gamma_1+\lambda_1}{\gamma_2+\lambda_2}(\xi,\tau)
\end{equation}

\medskip

Consider the usual Jacobi theta function
\begin{equation}
\Theta_p(z) = (z;p)_\infty \, (pz^{-1};p)_\infty \, (p;p)_\infty
\end{equation}
where the infinite multiple products are defined by
\begin{equation}
(z;p_1,\dots,p_m)_\infty = \prod_{n_i \ge 0} (1-z p_1^{n_1} \dots 
p_m^{n_m})
\end{equation}
It satisfies $\Theta_p(pz) = \Theta_p(z^{-1}) = -z^{-1} \Theta_p(z)$. \\
The Jacobi theta functions with rational characteristics 
$(\gamma_1,\gamma_2) \in \sfrac{1}{N} \ZZ \times \sfrac{1}{N} \ZZ$ can be 
expressed in terms of the $\Theta_{p}$ function as
\begin{equation}
\vartheta\car{\gamma_1}{\gamma_2}(\xi,\tau) = (-1)^{2\gamma_1\gamma_2} 
\, p^{\frac{1}{2}\gamma_1^2} \, z^{2\gamma_1} \, 
\Theta_{p}(-e^{2i\pi\gamma_2} p^{\gamma_1+\frac{1}{2}} z^2) 
\end{equation}
where $p = e^{2i\pi\tau}$ and $z = e^{i\pi \xi}$.

\subsubsection*{Multiple Gamma and Sine functions}

$\Gamma_r$ is the multiple Gamma function of order $r$ given by
\begin{equation}
\Gamma_{r}(x \,|\, \omega_{1},\dots,\omega_{r}) = \exp \left( \left. 
\frac{\partial}{\partial s} \; \zeta_{r}(x,s \,|\, 
\omega_{1},\dots,\omega_{r}) \right\vert_{s=0} \right)
\end{equation}
where
\begin{equation}
\zeta_{r}(x,s \,|\, \omega_{1},\dots,\omega_{r}) = 
\sum_{n_{1},\dots,n_{r} \ge 0} 
(x+n_{1}\omega_{1}+\dots+n_{r}\omega_{r})^{-s}
\end{equation}
is the multiple zeta function. In particular $\Gamma_1(x\,|\,\omega_1) = 
\frac{\displaystyle\omega_1^{x/\omega_1}} {\displaystyle \sqrt{2\pi 
\omega_1}} \Gamma\left(\frac{\displaystyle 
x}{\displaystyle\omega_1}\right)$. It has the following property
\begin{equation}
\frac{\Gamma_r(x+\omega_i\,|\,\omega_1,\cdots \omega_r)} 
{\Gamma_r(x\,|\,\omega_1,\cdots \omega_r)} = 
\frac{1}{\Gamma_{r-1}(x\,|\,\omega_1,\cdots \omega_{i-1}, 
\omega_{i+1}\cdots,\omega_r)} 
\end{equation}
Multiple sine functions of order $r$ are defined by \cite{Bar1901,JM96}
\begin{equation}
S_r(x\,|\,\omega_1,\cdots \omega_r) = \Gamma_r(x\,|\,\omega_1,\cdots 
\omega_r)^{-1} \Gamma_r(\omega_1+\cdots+\omega_r - x\,|\,\omega_1,\cdots 
\omega_r)^{(-1)^r} 
\end{equation}
They satisfy
\begin{equation}
\frac{S_r(x+\omega_i\,|\,\omega_1,\cdots \omega_r)} 
{S_r(x\,|\,\omega_1,\cdots \omega_r)} = 
\frac{1}{S_{r-1}(x\,|\,\omega_1,\cdots \omega_{i-1}, 
\omega_{i+1}\cdots,\omega_r)}
\end{equation}
In particular $S_1(x\,|\,\omega_1) = 2\sin\left(\frac{\displaystyle \pi 
x}{\displaystyle \omega_1} \right)$ and Barnes' double sine function of 
periods $\omega_{1}$ and $\omega_{2}$ is given by
\begin{equation}
S_{2}(x \,|\, \omega_{1},\omega_{2}) = 
\frac{\Gamma_{2}(\omega_{1}+\omega_{2}-x \,|\, \omega_{1},\omega_{2})} 
{\Gamma_{2}(x \,|\, \omega_{1},\omega_{2})}
\end{equation}
One has the following properties
\begin{align}
& \frac{S_{2}(x+\omega_{1} \,|\, \omega_{1},\omega_{2})}{S_{2}(x \,|\, 
\omega_{1},\omega_{2})} = \frac{1}{\displaystyle 2\sin\frac{\pi 
x}{\omega_{2}}} \\
& S_{2}(x \,|\, \omega_{1},\omega_{2}) S_{2}(-x \,|\, 
\omega_{1},\omega_{2}) = -4 \sin\frac{\pi x}{\omega_{1}} \sin\frac{\pi 
x}{\omega_{2}}
\end{align}

\section{Deformation of an algebra}
\setcounter{equation}{0}

\begin{definition}
Let $\cA$ be a Lie algebra over the field $\CC$ with bracket 
$[\,\cdot\,,\,\cdot\,]$. A $p$-cochain is a $p$-linear skew-symmetric 
map $C_{p}: \wedge^p \cA \to \cA$. The Chevalley coboundary operator 
$\partial$ maps a $p$-cochain to $(p+1)$-cochains as
\begin{align}
\partial C_{p}(u_{0},\ldots,u_{p}) &= \sum_{i=0}^p (-1)^i \; 
\big[u_{i} , C_{p}(u_{0}, \ldots, \widehat{u_{i}}, \ldots, u_{p}) 
\big] \nonumber \\
&+ \sum_{0 \le i<j \le p} (-1)^{i+j} \; C_{p} \big([u_{i},u_{j}], 
u_{0}, \ldots, \widehat{u_{i}}, \ldots, \widehat{u_{j}}, \ldots, 
u_{p} \big) 
\end{align}
where $\widehat{u_{i}}$ means that $u_{i}$ is omitted.
\end{definition}
It can be checked that $\partial$ satisfies $\partial^2 C_{p} = 0$. A 
$p$-cochain $C_{p}$ is called a $p$-cocycle if $\partial C_{p} = 0$, i.e. 
is an element of $\text{Ker} \partial$. The space of $p$-cocycles is 
denoted $\cZ^p(\cA,\cA)$. A $p$-cochain $C_{p}$ is called a $p$-coboundary 
if $C_{p} = \partial C_{p-1}$, i.e. is an element of $\text{Im} \partial$. 
The space of $p$-coboundaries is denoted $\cB^p(\cA,\cA)$. Since 
$\partial^2 = 0$, one has $\cB^p(\cA,\cA) \subset \cZ^p(\cA,\cA)$.
\begin{definition}
Let $\cB^p(\cA,\cA)$ and $\cZ^p(\cA,\cA)$ be the spaces of 
$p$-coboundaries and $p$-cocycles with respect to the Chevalley 
coboundary operator $\partial$. The quotient $\cH^p(\cA,\cA) = 
\cZ^p(\cA,\cA) / \cB^p(\cA,\cA)$ is called the $p$-th Chevalley 
cohomology space.
\end{definition}
The space $\cH^p(\cA,\cA)$ describes the non-trivial $p$-cocycles, i.e. 
cocycles which are not coboundaries. One has $\cH^*(\cA,\cA) \equiv 
\oplus_{p} \cH^p(\cA,\cA) = \text{Ker} \partial/\text{Im} \partial$.

\medskip

The theory of Chevalley cohomology is intimately related to the deformation 
of a Lie algebra. This last notion has been precisely defined by 
Gerstenhaber \cite{Ger64,GS92,Stern}.
\begin{definition}
Let $\cA$ be Lie algebra over the field $\CC$ and denote by 
$\CC[[\hbar]]$ the ring of formal series in the parameter $\hbar$. A 
deformation of $\cA$ is an algebra $\cA_{\hbar}$ over $\CC[[\hbar]]$ 
such that $\cA_{\hbar}/\hbar\cA_{\hbar} \approx \cA$. Two deformations 
$\cA_{\hbar}$ and $\cA'_{\hbar}$ are equivalent if they are isomorphic 
over $\CC[[\hbar]]$. A deformation $\cA_{\hbar}$ is trivial if it is 
isomorphic to the original algebra $\cA$ (considered as a 
$\CC[\hbar]]$-algebra).
\end{definition}

Let $\{u_{i}\}$ be the set of generators of a Lie algebra $\cA$ with 
commutation relations
\begin{equation}
[u_{i},u_{j}] = f_{ij}^{k} \, u_{k}
\end{equation}
Consider the deformation $\cA_{\hbar}$ of the enveloping algebra of $\cA$: 
there exists a new bracket $[\,\cdot\,,\,\cdot\,]_{\hbar}$ such that
\begin{equation}
[u_{i},u_{j}]_{\hbar} = f_{ij}^{k} \, u_{k} + \sum_{p=1}^{\infty} 
\hbar^p \varphi_{p}(u_{i},u_{j})
\end{equation}
where the $\varphi_{p}$ are antisymmetric bilinear maps taking values in 
$\cA$, i.e. Chevalley 2-cochains. Imposing to the new bracket 
$[\,\cdot\,,\,\cdot\,]_{\hbar}$ to satisfy the Jacobi identity, one gets 
the following constraints ($\partial$ denoting the Chevalley coboundary 
operator):
\begin{align}
& \partial \varphi_{1} = 0 \\
& \partial \varphi_{p} = \sum_{r+s=p} \Big( 
\varphi_{r}\big(\varphi_{s}(u_{i},u_{j}),u_{k}\big) + 
\varphi_{r}\big(\varphi_{s}(u_{j},u_{k}),u_{i}\big) + 
\varphi_{r}\big(\varphi_{s}(u_{k},u_{i}),u_{j}\big) \Big) \equiv 
\psi_{p}(u_{i},u_{j},u_{k}) \nonumber \\
& \hspace*{14cm} (p>1)
\label{eq:phin}
\end{align}
Therefore $\varphi_{1}$ is a 2-cocycle, while $\varphi_{p}$ is determined 
by the $\varphi_{r}$ with $r<p$, up to 2-cocycles. Now, if equation 
(\ref{eq:phin}) is satisfied up to order $p$ (i.e. the deformation is 
consistent up to order $p$), one can show by a direct calculation that the 
r.h.s. of (\ref{eq:phin}) at order $p+1$ is a 3-cocycle: $\partial 
\psi_{p+1} = 0$. If one wants equation (\ref{eq:phin}) to be satisfied at 
order $p+1$ (i.e. the deformation extends to order $p+1$), this 3-cocycle 
must be indeed a 3-coboundary: $\psi_{p+1} = \partial \varphi_{p+1}$. It 
follows:
\begin{proposition}
The third Chevalley cohomology space $\cH^3(\cA,\cA)$ classifies the 
obstructions to deformations of a Lie algebra. In particular, if 
$\cH^3(\cA,\cA)$ is the null space, any 2-cocycle of $\cZ^2(\cA,\cA)$ 
leads to a deformation of $\cA$.
\end{proposition}
The question is now to classify all \emph{non-equivalent} deformations. Let 
$\cA_{\hbar}$ and $\cA'_{\hbar}$ be two deformations of the Lie algebra 
$\cA$. $\cA_{\hbar}$ and $\cA'_{\hbar}$ are equivalent if there is an 
isomorphism $\cI = 1 + \sum_{p=1}^{\infty} \hbar^p \cI_{p}$ of $\cA$ such 
that $\cI([u_{i},u_{j}]_{\hbar}) = [\cI(u_{i}),\cI(u_{j})]'_{\hbar}$.

A deformation at order 1 is trivial if the 2-cocycle $\varphi_{1}$ is in 
fact a 2-coboundary. Using similar arguments as above, if two deformations 
are equivalent up to order $p$, they differ only by 2-cocycles which are 
2-coboundaries. If one wants to extend the equivalence at order $p+1$, this 
imposes that the 2-cocycle at order $p+1$ is also a 2-coboundary. It 
follows:
\begin{proposition}
The obstructions to equivalence between deformations of a Lie algebra 
$\cA$ lie in the Chevalley cohomology space $\cH^2(\cA,\cA)$. In 
particular, if $\cH^2(\cA,\cA)$ is the null space, all deformations are 
trivial.
\end{proposition}

\medskip

In the case of Hopf algebras, one has to deal also with the deformation of 
the coproduct. If $\cA$ is a Hopf algebra with coproduct $\Delta$, a 
deformation $\cA_{\hbar}$ is endowed with a deformed coproduct 
$\Delta_{\hbar} = \Delta + \sum_{p=1}^{\infty} \hbar^p \Delta_{p}$ where 
$\Delta_{p}$ are homomorphisms from $\cA$ to $\cA \otimes \cA$. As before, 
one has to introduce suitable cohomologies (related to the algebra and the 
coalgebra structures). See \cite{ChaPre,Stern} and references therein for 
more details.


\cleardoublepage


\baselineskip=15pt
\markboth{Bibliography}{Bibliography}

\begin{thebibliography}{10}

\bibitem{ABF}
G.~Andrews, R.J. Baxter, and P.J. Forrester, \emph{Eight-vertex {SOS} model and
  generalized {R}ogers--{R}amanujan type identities}, J. Stat. Phys.
  \textbf{35} (1984), 193--266.

\bibitem{AACFR}
D.~Arnaudon, J.~Avan, N.~Cramp\'{e}, L.~Frappat, and E.~Ragoucy,
  \emph{{R}-matrix presentation for (super)-{Y}angians}, e-print arXiv.org
  \texttt{math.QA/0111325}.

\bibitem{AAFR}
D.~Arnaudon, J.~Avan, L.~Frappat, and E.~Ragoucy, \emph{Yangian and quantum
  universal solutions of {G}ervais--{N}eveu--{F}elder equations}, to appear in
  Commun. Math. Phys., e-print arXiv.org \texttt{math.QA/0104181}.

\bibitem{AAFRR3}
D.~Arnaudon, J.~Avan, L.~Frappat, E.~Ragoucy, and M.~Rossi, \emph{On the
  quasi-{H}opf structure of deformed double {Y}angians}, Lett. Math. Phys.
  \textbf{51} (2000), 193--204, e-print arXiv.org \texttt{math.QA/0001034}.

\bibitem{Cla}
D.~Arnaudon, J.~Avan, L.~Frappat, E.~Ragoucy, and M.~Rossi, \emph{Towards a
  cladistics of double {Y}angians and elliptic algebras}, J. Phys. A (Math.
  Gen.) \textbf{33} (2000), 6279--6309, e-print arXiv.org
  \texttt{math.QA/9906189}.

\bibitem{AAFRo}
D.~Arnaudon, J.~Avan, L.~Frappat, and M.~Rossi, \emph{Deformed double {Y}angian
  structures}, Rev. Math. Phys. \textbf{12} (2000), 945--963, e-print arXiv.org
  \texttt{math.QA/9905100}.

\bibitem{ABRR}
D.~Arnaudon, E.~Buffenoir, E.~Ragoucy, and Ph. Roche, \emph{Universal solutions
  of quantum dynamical {Y}ang--{B}axter equations}, Lett. Math. Phys.
  \textbf{44} (1998), 201--214, e-print arXiv.org \texttt{q-alg/9712037}.

\bibitem{ABB96}
J.~Avan, O.~Babelon, and E.~Billey, \emph{The {G}ervais--{N}eveu--{F}elder
  equation and quantum {C}alogero--{M}oser systems}, Commun. Math. Phys.
  \textbf{178} (1996), 281--299, e-print arXiv.org \texttt{hep-th/9505091}.

\bibitem{AFRS3}
J.~Avan, L.~Frappat, M.~Rossi, and P.~Sorba, \emph{Deformed ${W}_{N}$ algebras
  from elliptic $sl_{N}$ algebras}, Commun. Math. Phys. \textbf{199} (1999),
  697--728, e-print arXiv.org \texttt{math.QA/9801105}.

\bibitem{AFRS5}
J.~Avan, L.~Frappat, M.~Rossi, and P.~Sorba, \emph{Universal construction of
  q-deformed {W}-algebras}, Commun. Math. Phys. \textbf{202} (1999), 445--461,
  e-print arXiv.org \texttt{math.QA/9807048}.

\bibitem{SKAO96}
H.~Awata, H.~Kubo, S.~Odake, and J.~Shiraishi, \emph{A quantum deformation of
  the {V}irasoro algebra and the {M}acdonald symmetric functions}, Lett. Math.
  Phys. \textbf{38} (1996), 33--51, e-print arXiv.org \texttt{q-alg/9507034}.

\bibitem{AKOS96}
H.~Awata, H.~Kubo, S.~Odake, and J.~Shiraishi, \emph{Quantum ${W}_{N}$ algebras
  and {M}acdonald polynomials}, Commun. Math. Phys. \textbf{179} (1996),
  401--415, e-print arXiv.org \texttt{q-alg/9508011}.

\bibitem{Bab91}
O.~Babelon, \emph{Universal exchange algebra for {B}loch waves and {L}iouville
  theory}, Commun. Math. Phys. \textbf{139} (1991), 619--643.

\bibitem{BBB96}
O.~Babelon, D.~Bernard, and E.~Billey, \emph{A quasi-{H}opf algebra
  interpretation of quantum $3-j$ and $6-j$ symbols and difference equations},
  Phys. Lett. B \textbf{375} (1996), 89--97, e-print arXiv.org
  \texttt{q-alg/9511019}.

\bibitem{Bar1901}
E.W. Barnes, \emph{The theory of the double gamma function}, Philos. Trans.
  Roy. Soc. A \textbf{196} (1901), 265--388.

\bibitem{Bax72}
R.J. Baxter, \emph{Partition function of the eight-vertex lattice model}, Ann.
  Phys. \textbf{70} (1972), 193--228.

\bibitem{Bax82}
R.J. Baxter, \emph{Exactly solved models in statistical mechanics}, Academic
  Press, London, 1982.

\bibitem{Bel81}
A.A. Belavin, \emph{Dynamical symmetry of integrable quantum systems}, Nucl.
  Phys. B \textbf{180} (1981), 189--200.

\bibitem{BL93}
D.~Bernard and A.~LeClair, \emph{The quantum double in integrable quantum field
  theory}, Nucl. Phys. B \textbf{399} (1993), 709--748.

\bibitem{BR99}
E.~Buffenoir and Ph. Roche, \emph{Harmonic analysis on the quantum {L}orentz
  group}, Commun. Math. Phys. \textbf{207} (1999), 499--555, e-print arXiv.org
  \texttt{q-alg/9710022}.

\bibitem{ChaPre}
V.~Chari and A.~Pressley, \emph{A guide to quantum groups}, Cambridge
  University Press, 1994.

\bibitem{ChCh}
D.V. Chudnovsky and G.V. Chudnovsky, \emph{Completely {X}-symmetric
  ${S}$-matrices corresponding to theta functions}, Phys. Lett. A \textbf{81}
  (1981), 105.

\bibitem{DJMO}
E.~Date, M.~Jimbo, T.~Miwa, and M.Okado, \emph{Fusion of the eight-vertex {SOS}
  model}, Lett. Math. Phys. \textbf{12} (1986), 209--215.

\bibitem{DF93}
J.~Ding and I.B. Frenkel, \emph{Isomorphism of two realizations of quantum
  affine algebra ${U}_q(\widehat{gl(n)})$}, Commun. Math. Phys. \textbf{156}
  (1993), 277--300.

\bibitem{Dri85}
V.G. Drinfeld, \emph{Hopf algebras and the quantum {Y}ang--{B}axter equation},
  Sov. Math. Dokl. \textbf{32} (1985), 254--258.

\bibitem{Dri86}
V.G. Drinfeld, \emph{Quantum groups}, Proc. ICM-1986, Berkeley, California
  (Academic Press, ed.), vol.~1, 1986, pp.~798--820.

\bibitem{Dri88}
V.G. Drinfeld, \emph{A new realization of {Y}angians and quantized affine
  algebras}, Sov. Math. Dokl. \textbf{36} (1988), 212--216.

\bibitem{Dri90}
V.G. Drinfeld, \emph{Quasi-{H}opf algebras}, Leningrad Math. Journ. \textbf{1}
  (1990), 1419--1457.

\bibitem{EF97}
B.~Enriquez and G.~Felder, \emph{Elliptic quantum groups
  ${E}_{\tau,\eta}(sl_2)$ and quasi-{H}opf algebras}, Commun. Math. Phys.
  \textbf{195} (1998), 651--689, e-print arXiv.org \texttt{q-alg/9703018}.

\bibitem{FRT}
L.D. Faddeev, N.Yu. Reshetikhin, and L.A. Takhtajan, \emph{Quantization of
  {L}ie groups and {L}ie algebras}, Leningrad Math. Journ. \textbf{1} (1990),
  193--225.

\bibitem{FF96}
B.~Feigin and E.~Frenkel, \emph{Quantum {W}-algebras and elliptic algebras},
  Commun. Math. Phys. \textbf{178} (1996), 653--678, e-print arXiv.org
  \texttt{q-alg/9508009}.

\bibitem{Fel94}
G.~Felder, \emph{Elliptic quantum groups}, Proc. ICMP Paris, 1994, e-print
  arXiv.org \texttt{hep-th/9412207}.

\bibitem{FIJKMY}
O.~Foda, K.~Iohara, M.~Jimbo, R.~Kedem, T.~Miwa, and H.~Yan, \emph{An elliptic
  quantum algebra for $\widehat{sl}_2$}, Lett. Math. Phys. \textbf{32} (1994),
  259--268, e-print arXiv.org \texttt{hep-th/9403094}.

\bibitem{FIJKMY95}
O.~Foda, K.~Iohara, M.~Jimbo, R.~Kedem, T.~Miwa, and H.~Yan, \emph{Notes on
  highest weight modules of the elliptic algebra ${A}_{q,p}(\widehat{sl}_2)$},
  Prog. Theor. Phys. Suppl. \textbf{118} (1995), 1--34, e-print arXiv.org
  \texttt{hep-th/9405058}.

\bibitem{FR96}
E.~Frenkel and N.Yu. Reshetikhin, \emph{Quantum affine algebras and
  deformations of the {V}irasoro and {W}-algebras}, Commun. Math. Phys.
  \textbf{178} (1996), 237--264, e-print arXiv.org \texttt{q-alg/9505025}.

\bibitem{Fro97b}
C.~Fr{\o}nsdal, \emph{Generalization and exact deformations of quantum groups},
  Publication RIMS Kyoto University \textbf{33} (1997), 91--149, e-print
  arXiv.org \texttt{q-alg/9606020}.

\bibitem{Fro97a}
C.~Fr{\o}nsdal, \emph{Quasi-{H}opf deformations of quantum groups}, Lett. Math.
  Phys. \textbf{40} (1997), 117--134, e-print arXiv.org \texttt{q-alg/9611028}.

\bibitem{Ger64}
M.~Gerstenhaber, \emph{On the deformation of rings and algebras}, Ann. Math.
  \textbf{79} (1964), 59--103.

\bibitem{GS92}
M.~Gerstenhaber and S.D. Schack, \emph{Algebras, bialgebras, quantum groups and
  algebraic deformations}, Contemporary Mathematics \textbf{134} (1992),
  51--92, In ``Deformation theory and quantum groups with applications to
  mathematical physics''.

\bibitem{GN84}
J.L. Gervais and A.~Neveu, \emph{Novel triangle relation and absence of
  tachyons in {L}iouville string field theory}, Nucl. Phys. B \textbf{238}
  (1984), 125--141.

\bibitem{IIJMNT}
M.~Idzumi, K.~Iohara, M.~Jimbo, T.~Miwa, T.~Nakashima, and T.~Tokihiro,
  \emph{Quantum affine symmetry in vertex models}, Int. J. Mod. Phys.
  \textbf{A8} (1993), 1479--1511, e-print arXiv.org \texttt{hep-th/9208066}.

\bibitem{Io96}
K.~Iohara, \emph{Bosonic representations of {Y}angian double ${DY}_{\hbar}(g)$
  with $g = gl_{N},sl_{N}$}, J. Phys. A (Math. Gen.) \textbf{29} (1996),
  4593--4621, e-print arXiv.org \texttt{q-alg/9603033}.

\bibitem{Jim85}
M.~Jimbo, \emph{A q-difference analogue of ${U}(g)$ and the {Y}ang--{B}axter
  equation}, Lett. Math. Phys. \textbf{10} (1985), 63--69.

\bibitem{Jim86}
M.~Jimbo, \emph{Quantum {R}-matrix for the generalized {T}oda system}, Commun.
  Math. Phys. \textbf{102} (1986), 537--547.

\bibitem{JKKMW}
M.~Jimbo, R.~Kedem, H.~Konno, T.~Miwa, and R.~Weston, \emph{Difference
  equations in spin chains with a boundary}, Nucl. Phys. B \textbf{448} (1995),
  429--456, e-print arXiv.org \texttt{hep-th/9502060}.

\bibitem{JKM}
M.~Jimbo, H.~Konno, and T.~Miwa, \emph{Massless {XXZ} model and degeneration of
  the elliptic algebra ${A}_{q,p}(\widehat{sl}_2)$}, Math. Phys. Studies
  \textbf{20} (1997), 117--138, e-print arXiv.org \texttt{hep-th/9610079}.

\bibitem{JKOS2}
M.~Jimbo, H.~Konno, S.~Odake, and J.~Shiraishi, \emph{Elliptic algebra
  ${U}_{q,p}(\widehat{sl}_2)$: {D}rinfel'd currents and vertex operators},
  Commun. Math. Phys. \textbf{199} (1999), 605--647, e-print arXiv.org
  \texttt{math.QA/9802002}.

\bibitem{JKOS}
M.~Jimbo, H.~Konno, S.~Odake, and J.~Shiraishi, \emph{Quasi-{H}opf twistors for
  elliptic quantum groups}, Transformation Groups \textbf{4} (1999), 303--327,
  e-print arXiv.org \texttt{q-alg/9712029}.

\bibitem{JKMO88}
M.~Jimbo, A.~Kuniba, T.~Miwa, and M.~Okado, \emph{The ${A}({N})^{(1)}$ face
  models}, Commun. Math. Phys. \textbf{119} (1988), 543--565.

\bibitem{JM96}
M.~Jimbo and T.~Miwa, \emph{Quantum {KZ} equation with $|q|=1$ and correlation
  functions of the {XXZ} model in the gapless regime}, J. Phys. A (Math. Gen).
  \textbf{29} (1996), 2923--2958, e-print arXiv.org \texttt{hep-th/9601135}.

\bibitem{JMO88}
M.~Jimbo, T.~Miwa, and M.~Okado, \emph{Solvable lattice models related to the
  vector representation of classical simple {L}ie algebras}, Commun. Math.
  Phys. \textbf{116} (1988), 507--525.

\bibitem{JS97}
B.~Jur\v{c}o and P.~Schupp, \emph{{AKS} scheme for face and
  {C}alogero--{M}oser--{S}utherland type models}, e-print arXiv.org
  \texttt{solv-int/9710006}.

\bibitem{Kho95}
S.M. Khoroshkin, \emph{Central extension of the {Y}angian double}, Collection
  SMF, 7{\`e}me rencontre du contact franco-belge en alg{\`e}bre, Reims, 1995,
  e-print arXiv.org \texttt{q-alg/9602031}.

\bibitem{KLP}
S.M. Khoroshkin, D.~Lebedev, and S.~Pakuliak, \emph{Elliptic algebra
  ${A}_{q,p}(\widehat{sl}_2)$ in the scaling limit}, Commun. Math. Phys.
  \textbf{190} (1998), 597--627, e-print arXiv.org \texttt{q-alg/9702002}.

\bibitem{KLCP}
S.M. Khoroshkin, A.~LeClair, and S.~Pakuliak, \emph{Angular quantization of the
  {S}ine--{G}ordon model at the free fermion point}, Adv. Theor. Math. Phys.
  \textbf{3} (1999), 1227--1287, e-print arXiv.org \texttt{hep-th/9904082}.

\bibitem{KT91}
S.M. Khoroshkin and V.N. Tolstoy, \emph{Universal {R}-matrix for quantized
  (super)algebras}, Commun. Math. Phys. \textbf{141} (1991), 599--617.

\bibitem{KT92}
S.M. Khoroshkin and V.N. Tolstoy, \emph{The uniqueness theorem for the
  universal {R}-matrix}, Lett. Math. Phys. \textbf{24} (1992), 231--244.

\bibitem{KT93}
S.M. Khoroshkin and V.N. Tolstoy, \emph{The {C}artan-{W}eyl basis and the
  universal {R}-matrix for quantum {K}ac--{M}oody algebras and superalgebras},
  Quantum symmetries, 1993, H. Doebner and V. Dobrev eds.

\bibitem{KT94}
S.M. Khoroshkin and V.N. Tolstoy, \emph{Twisting of quantum (super)algebras.
  {C}onnection of {D}rinfeld's and {C}artan--{W}eyl realizations for quantum
  affine algebras},  (1994), e-print arXiv.org \texttt{hep-th/9404036}.

\bibitem{KT96}
S.M. Khoroshkin and V.N. Tolstoy, \emph{Yangian double and rational
  {R}-matrix}, Lett. Math. Phys. \textbf{36} (1996), 373--402, e-print
  arXiv.org \texttt{hep-th/9406194}.

\bibitem{Kon97}
H.~Konno, \emph{Degeneration of the elliptic algebra $\elpa{2}$ and form
  factors in the {S}ine-{G}ordon theory}, to appear in the {CRM} series in
  mathematical physics (Springer Verlag, ed.), 1996, e-print arXiv.org
  \texttt{hep-th/9701034}.

\bibitem{Kon98}
H.~Konno, \emph{An elliptic algebra ${U}_{q,p}(\widehat{sl}_2)$ and the fusion
  {RSOS} model}, Commun. Math. Phys. \textbf{195} (1998), 373--403, e-print
  arXiv.org \texttt{q-alg/9709013}.

\bibitem{KuRe81}
P.P. Kulish and N.Yu. Reshetikhin, \emph{The quantum linear problem for the
  sine--gordon equation and higher representations}, Zap. Nauchn. Sem. LOMI
  \textbf{101} (1981), 101--110.

\bibitem{KuSkl82}
P.P. Kulish and E.V. Sklyanin, \emph{On the solutions of the {Yang}--{B}axter
  equation}, J. Sov. Math. \textbf{19} (1982), 1596--1620.

\bibitem{Molev2001}
A.I. Molev, \emph{Yangians and their applications}, Handbook of {A}lgebra, vol.
  3, Elsevier, to appear.

\bibitem{MNO96}
A.I. Molev, M.~Nazarov, and G.~Olshanski\u{\i}, \emph{Yangians and classical
  {L}ie algebras}, Russ. Math. Surveys \textbf{51} (1996), 205--282, e-print
  arXiv.org \texttt{hep-th/9409025}.

\bibitem{Skl80}
E.K. Sklyanin, \emph{Quantum version of the method of inverse scattering
  problem}, J. Sov. Math. \textbf{19} (1982), 1546--1596, translated from Zap.
  Nauchn. Sem. \textbf{95} (1980), 55--128.

\bibitem{Skl82}
E.K. Sklyanin, \emph{Some algebraic structures connected with the
  {Y}ang--{B}axter equation}, Funct. Anal. Appl. \textbf{16} (1982), 263--270.

\bibitem{Skl83}
E.K. Sklyanin, \emph{Some algebraic structures connected with the
  {Y}ang--{B}axter equation. {R}epresentations of quantum algebras}, Funct.
  Anal. Appl. \textbf{17} (1983), 273--284.

\bibitem{Skl85}
E.K. Sklyanin, \emph{On an algebra generated by quadratic relations}, Uspekhi
  Mat. Nauk. \textbf{40} (1985), 214.

\bibitem{SF78}
E.K. Sklyanin and L.D. Faddeev, \emph{Quantum mechanical approach to completely
  integrable field theory models}, Sov. Phys. Dokl. \textbf{23} (1978),
  902--904.

\bibitem{SFT80}
E.K. Sklyanin, L.D. Faddeev, and L.A. Takhtajan, \emph{The quantum inverse
  problem method}, Theor. Math. Phys. \textbf{40} (1980), 688--706, translated
  from Teor. Mat. Fiz. \textbf{40} (1979), 194--220.

\bibitem{Stern}
D.~Sternheimer, \emph{Deformation quantization: twenty years after},
  International conference on particles, fields and gravitation, Lodz, Poland,
  1998, e-print arXiv.org \texttt{math.QA/9809056}.

\bibitem{Yang67}
C.N. Yang, \emph{Some exact results for the many-body problem in one dimension
  with repulsive delta-function interaction}, Phys. Rev. Lett. \textbf{19}
  (1967), 1312--1314.

\end{thebibliography}
\end{document}